\RequirePackage{fix-cm}

\documentclass[twocolumn]{autart}    

\usepackage{multirow}
\usepackage{eurosym}
\usepackage{pifont}
\usepackage[utf8]{inputenc}
\usepackage{mathtools}
\usepackage{amssymb}
\usepackage{amsfonts}
\usepackage{graphicx}
\usepackage{anyfontsize}
\usepackage{subcaption}
\usepackage[numbers]{natbib}
\usepackage{nicefrac}
\usepackage{svg}
\usepackage[per-mode=fraction]{siunitx}
\usepackage{paralist}
\usepackage{textcomp,gensymb}
\graphicspath{{Figures}}
\usepackage{float}
\usepackage[english]{babel}
\usepackage{multicol}
\usepackage{wrapfig}
\usepackage{pdfpages}
\usepackage{xcolor}
\usepackage{fancyhdr}
\usepackage{algorithm}
\usepackage{algpseudocode}
\usepackage{enumitem}
\usepackage{hyperref}

\newcommand{\cmark}{\ding{51}}
\newcommand{\xmark}{\ding{55}}

\newcommand{\qeq}{\bar{q}}

\newcommand{\EigN}{\mathcal{M}}
\newcommand{\R}{\mathbb{R}}                                

\newcommand{\M}{\mathcal{M}}                                
\newcommand{\Q}{\mathcal{Q}}                                
\newcommand{\X}{\mathfrak{X}}                               
\newcommand{\ext}{\text{d}}                                 

\newcommand{\coord}{q}                                      

\newcommand{\CoordC}{Z}
\newcommand{\coordC}{z}

\newcommand{\state}{x}
\newcommand{\eqstate}{\bar{x}}

\newcommand{\stateH}{\zeta}

 \colorlet{myellow}{green!10!orange}

\begin{document}

\begin{frontmatter}
\title{Symmetric Lyapunov Subcenter Manifolds for Periodic Regulation of Mechanical Systems\thanksref{footnoteinfo}} 

\thanks[footnoteinfo]{This work is supported by the Advanced Grant M-Runner (Grant Agreement No. 835284) by the European Research Council (ERC). \\
This work is supported by the Advanced Grant PortWings (Grant Agreement No. 787675) by the ERC.}

\author[Twente]{Yannik P. Wotte}\ead{y.p.wotte@utwente.nl},    
\author[TUM,DLR]{Arne Sachtler}\ead{arne.sachtler@dlr.de},               
\author[TUM,DLR]{Alin Albu-Sch\"{a}ffer}\ead{alin.albu-schaeffer@tum.de},  
\author[Twente]{Stefano Stramigioli}\ead{s.stramigioli@utwente.nl},
\author[Delft,DLR]{Cosimo Della Santina}\ead{cosimodellasantina@gmail.com}

\address[Twente]{Robotics and Mechatronics, University of Twente (UT), Drienerlolaan 5, 7522 NB Enschede, Netherlands}  
\address[TUM]{Technical University of Munich (TUM), Munich, Germany}             
\address[DLR]{Institute of Robotics and Mechatronics, German Aerospace Center (DLR), Oberpfaffenhofen, Germany}
\address[Delft]{Cognitive Robotics Department, Delft University of Technology (TU Delft), Netherlands}        


\begin{abstract} 
Multi-body mechanical systems have rich internal dynamics, whose solutions can be exploited as energy-efficient control targets. Yet, solutions non-trivially depend on system parameters, obscuring feasible properties for use as target trajectories. 
For periodic regulation tasks in robotics applications, we investigate properties of nonlinear oscillations collected in Lyapunov subcenter manifolds (LSMs) of conservative mechanical systems (CMs). 
Using a time-symmetry of CMs, it is shown that mild non-resonance conditions guarantee that LSMs exclusively consist of oscillations between two points of zero velocity. The existence of a unique generator is proven, which is a connected, 1D manifold that collects these points of zero velocity for a given LSM.  
Furthermore, it is shown that an additional spatial symmetry provides LSMs with yet stronger properties of Rosenberg manifolds. Here all oscillations pass through a unique equilibrium configuration, which can be favorable for control applications.
These theoretical results are numerically confirmed on two mechanical systems: a double pendulum and a 5-link pendulum.
\end{abstract}

\end{frontmatter}





\section{Introduction}\label{s1}

\begin{figure*}[h]
    \centering
    \def\svgwidth{1.3\linewidth}\scriptsize
    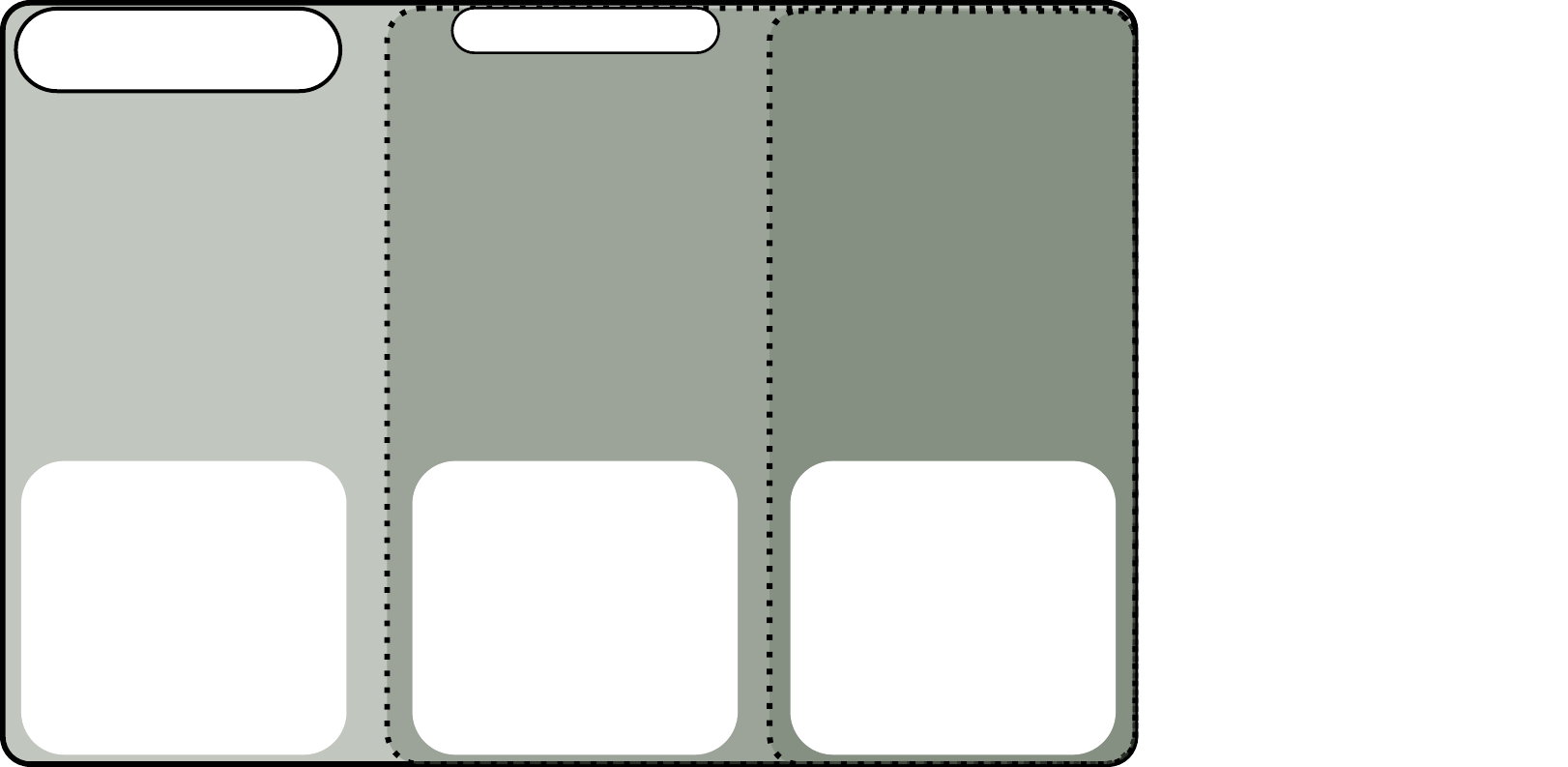
    \caption{Summary of the article: for conservative mechanical systems with configuration $\state \in \Q$, momentum $P\in T^*_\state \Q$, and projecting all data to $\Q$ for ease of visualization. Lyapunov subcenter manifolds (panels \textbf{a} and \textbf{d}) are families of general periodic oscillations $(\state(t),P(t))$ springing from an equilibrium $(\eqstate,0)$. We prove conditions for LSMs to have stronger properties: Theorems~\ref{thm:LSM_mechanics_geometric} and~\ref{thm:Eigenmanifold_Condition} make it highly common for LSMs to become weak Eigenmanifolds (panels \textbf{b} and \textbf{e}) that collect periodic brake trajectories (oscillating between brake points $(\state,0)$), and Eigenmanifolds that collect periodic brake trajectories whose configuration trajectory does not self-intersect (called geometric eigenmodes). Theorem~\ref{Theorem:1} shows that yet stronger conditions turn LSMs into (weak) Rosenberg manifolds (panels \textbf{c} and \textbf{f}), where all modal configurations pass through $\eqstate$. In both (weak) Eigenmanifolds and (weak) Rosenberg manifolds, brake points are collected on a connected 1D submanifold that we call the generator.}
    \label{fig:overview}
\end{figure*}

Practical advances in robotics are pushing the need for understanding robust features of the nonlinear dynamics of multibody mechanical systems~\citep*{Wensing2023,Noemie2023,Baddoo2023}. For example, a key challenge on the frontier of robotics lies in enabling robots to move efficiently while performing periodic motions (such as repetitive industrial operations and locomotion), thus avoiding motor torque limits and mitigating limited battery capacity~\citep*{Tong2024,Ficht2021,Taheri2023,Soori2023,Dovgopolik2023}.\\ 
%
%
A growing body of literature addresses this challenge by exciting internal robot dynamics for low-torque control implementations, with modal analysis serving as the basis for these works~\citep*{Santina2021,raff2022connecting,Cheng2024,Alora2023,Vu2024}. %
Given a system with state $\zeta$, control-input $u$ representing a motor torque, current or voltage and assuming control-affine dynamics $$\dot{\zeta} = f(\zeta) + g(\zeta) u\,.$$ Then the most energy-efficient control-input is simply $u = 0$: when tracking target trajectories that are natural motions of $f(\zeta)$, this corresponds to the steady-state control-input, forming the rationale behind exciting internal robot dynamics. For periodic regulation tasks, (periodic) natural motions from a dissipation-free model serve as target trajectories for the real system, where stabilizing controllers compensate for small dissipative effects but not the much larger inertial terms~\citep{Sachtler2024,Coelho2022}.\\
To use natural motions as control-targets in practice, optimization of robots and their internal dynamics is essential~\citep*{Wotte2023}. Yet, solutions and their properties can change drastically when parameters are changed. Robust properties of natural motions of the nonlinear and often chaotic open-loop dynamics are required~\citep*{AlbuSchaeffer2023}, forming the center of this investigation.\\
Robots are usually multibody systems whose configuration is modeled to evolve on Riemannian manifolds with curvature~\citep*{Bullo2004}, and general theoretical insights rely on a differential geometric formalism. Inisghts for such systems are challenging to derive by direct application of classic modal theory from nonlinear mechanics~\citep*{Rosenberg1966,Jezequel1991,Vakakis1996,Touz2006}, leading instead to postulated geometric definitions~\citep*{AlbuSchaeffer2020}.\\ 
%
In this article we go back to the foundations of modal theory in nonlinear dynamics, starting from Lyapunov subcenter manifold (LSM) theory to describe the families of oscillations that can be expected in multibody dynamical systems, and to derive robust, geometric properties (Figure~\ref{fig:overview}).\\ 
Given a dynamical system $\dot{\stateH} = f(\stateH)$ with a conserved quantity $H(\stateH)$ (\textit{i.e.}, $\dot{H} = 0$), LSMs are two-dimensional submanifolds of the state-space that collect periodic trajectories, and originate at equilibria where the linearized system has imaginary eigenvalues~\citep*{Kelley1967}. LSMs have been heavily investigated over the past decades~\citep*{Jain2022,Stoychev2023,Touze2025} and have found important applications to disparate scenarios~\citep*{Axas2024,Debeurre2024,Pozzi2024}. LSMs are well-studied: there are results on existence, uniqueness, differentiability~\citep*{Sijbrand1985,Llave1997}, and persistence under dissipative disturbances~\citep*{Llave2018} in the form of spectral submanifolds.\\
The theory of LSMs directly applies to conservative mechanical systems (CMs), \textit{i.e.}, including conservative multi-body mechanical systems, for which the Hamiltonian 
\begin{equation}\label{eq:abstract_hamiltonian}
    H(\state,P) = \frac{1}{2} M(\state)^{-1}(P,P) + V(\state)  
\end{equation} 
is the conserved quantity~--~where $\state \in \Q,P \in T^*_\state\Q$ denote abstract configuration and momentum variables, $M(\state)$ 
the symmetric, positive-definite inertia tensor, and $V:\Q \rightarrow \mathbb{R}$ the potential energy. It is well-known that unique LSMs exist around equilibria of CMs (\textit{i.e.}, minima of $V$) when the eigenvalues of the linearized system fulfill a non-resonance condition. However, high-level properties of the periodic trajectories collected in LSMs of CMs are apriori undetermined in this framework, see Figure~\ref{fig:overview} a \& d.\\ 
Our contribution utilizes the time-symmetry 
inherent to CMs to show that unique LSMs collect families of periodic brake trajectories, which are trajectories whose configuration oscillates back and forth between points with zero momentum. Uniqueness of LSMs is a robust feature over large ranges of parameters, and we argue that guaranteed brake trajectories make LSMs promising for e.g., optimization of target trajectories for pick-and-place tasks. In line with related literature, 
we refer to these families as \textit{Eigenmanifolds}, see Figure~\ref{fig:overview} b \& e.\\ 
Furthermore, we show that an additional spatial symmetry 
gives LSMs the stronger properties of Rosenberg manifolds~\citep*{AlbuSchaeffer2020}, see Figure~\ref{fig:overview} c \& f. 
In these structures all oscillation pass through the equilibrium configuration, which 
makes it easy to switch between modes by impulsive control actions~\citep*{Santina2021a}. We derive conditions on the inertia-tensor and potential energy function of CMs, for the existing LSMs to be Rosenberg manifolds. In local coordinates $q,p$ with equilibrium at $q = 0$, inertia matrix $M(q)\in\mathbb{R}^{n\times n}$ and potential energy $V(q)\in\mathbb{R}$, these conditions can be enforced as $M(q) = M(-q)$ and $V(q) = V(-q)$. This makes it possible to constrain Eigenmanifolds into Rosenberg manifolds at a design stage, enabling stricter properties of feasible target trajectories.\\
In summary, our contributions are:
\begin{enumerate}
    \item Theorem~\ref{thm:symmetry_propagation_system_to_LSM} on properties of LSMs induced by symmetries of dynamic systems, and the propagation of symmetries to individual trajectories in LSMs
    \item 
    Theorem~\ref{thm:Eigenmanifold_Condition}, stating that unique LSMs collect periodic brake trajectories and admit generators, which lends itself to a convenient parametrization of LSMs.
    \item Theorem~\ref{Theorem:1}, showing that conditions $M(q) = M(-q)$ and $V(q) = V(-q)$ let (weak) Eigenmanifolds in CMs exhibit the stronger properties of (weak) Rosenberg manifolds.
\end{enumerate}
The rest of the paper is organized as follows.
%
%
Sec.~\ref{ch-EM1:sec:background} introduces background theory on CMs, LSMs and discrete symmetries. 
Sec.~\ref{ch-EM1:sec:LSMs_in_CMS} presents conditions on $M(q)$ and $V(q)$ that determine existence and uniqueness of LSMs in CMs, and shows conditions under which self-symmetric LSMs consist of self-symmetric periodic motions. 
Sec.~\ref{ch-EM1:sec:theorems} derives the main theoretical contributions of this article, properties of self-symmetric LSMs and self-symmetric periodic motions in CMs. For periodic orbits, this motivates our definitions of geometric eigenmodes and geometric Rosenberg modes. Unique LSMs of CMs are proven to be collections of geometric eigenmodes, and to admit a generator, by which we justify the definition of Eigenmanifolds. A sufficient condition is derived for Eigenmanifolds to fulfill the stronger property that all trajectories pass through the equilibrium configuration, which aligns with the nonlinear mechanics definition of Rosenberg manifolds.  \\
Sec.~\ref{ch-EM1:sec:examples} verifies the results in numerical examples of CMs, highlighting that the existence of Eigenmanifolds and Rosenberg manifolds follows from rules implemented at a design stage. Supplementary code is found at~\cite{data-Link}.\\
We end with a conclusion in Sec.~\ref{ch-EM1:sec:conclusions}. Proofs are given in Appendix~\ref{ap:additional_proofs}, and Appendix B 
(only in the supplementary material~\cite{arxiV-Link}) contains background on 
choice of coordinates.

\subsection{Related literature}
Exciting and exploiting efficient nonlinear oscillations is an open challenge in nonlinear control theory~\cite{moyron2021orbital,romero2021orbital,saetre2021robust}, with important ramifications in robotics' locomotion and execution of industrial operations~\cite{garofalo2018passive,remy2018passive-gaits,kashiri2018overview,rosa2021topological,balderas2021performing}, and biomechanics~\cite{hatton2021geometry,seyfarth2021template,biewener2022stability}. 
To this end, Eigenmanifolds have been recently introduced (see \citep[Sec. 7]{AlbuSchaeffer2020}) as an extension of modal analysis to multi-body systems, in a way that naturally lends itself to control applications: design of mechanical systems~\cite{Ding2024} and shaping of Eigenmanifolds via the systems potential energy~\cite{Sachtler2022,Wotte2023} has been explored. Eigenmanifolds of interest are stabilized in simulation~\cite{Santina2021b,Santina2021c}, on real systems~\cite{Bjelonic2022}, can be stabilized in an energy-efficient manner~\cite{Coelho2022} and allow the swing-up of weakly actuated systems to high energies~\cite{Sachtler2024}. Eigenmanifolds have also been applied in classifying the accuracy of reduced order models in soft robotics~\cite{Santina2021d,Pustina2024}. \\
%
%
Eigenmanifold theory has its precursors in the study of nonlinear normal modes~\cite{Mikhlin2004,Kerschen2009,Mikhlin2023} for nonlinear mechanical systems with constant inertia tensors. Rosenberg modes are well-studied for conservative mechanical systems~\cite{Rosenberg1966,Peeters2009,Avramov2013}, and their symmetry properties are known by the condition that $V(q) = V(-q)$. Strict normal modes are a special case of Eigenmanifolds that exist under compatibility of the potential gradient with the  geodesic flow~\cite{AlbuStramigioli2021}. 

%
Within the nonlinear dynamics literature, LSMs are long studied~\cite{Kelley1967,Sijbrand1985,Llave1997} and their persistence in the limit of small dissipation was shown~\cite{Llave2018}, but this has not been connected to insights about Eigenmanifolds in multi-body dynamics. LSMs have also been applied to the model order reduction of partial differential equations~\cite{Vanderbauwhede1992}.\\
A general treatment of discrete symmetries in dynamical systems is given by Michael Field~\cite{Field2007}. Large classes of symmetries were also studied in the context of periodic orbits and nonlinear normal modes~\cite{Mishra1974,Chechin2000}. However, these works largely treat symmetries on vector and Banach spaces, and do not directly apply to Eigenmanifold theory. 
{An excellent overview to time-reversal symmetries in Hamiltonian systems is given by Lamb \& Roberts~\cite{Lamb1998}, covering many properties of direct relevance to Eigenmanifold theory, such as the existence of time-symmetric LSMs.}
Alomair \& Montaldi~\cite{Alomair2017} investigate the interplay of a time-reversal and a spatial symmetry of Hamiltonian systems in the context of periodic orbits, {proving the existence of a discrete number of doubly symmetric periodic orbits around symmetric equilibria. However, both of~\cite{Lamb1998,Alomair2017} leave properties of such symmetric orbits and their families implicit.}  

\section{Preliminaries and problem statement}\label{ch-EM1:sec:background}

This section introduces the Hamiltonian description of conservative mechanical systems, Lyapunov subcenter manifolds, and discrete symmetries. 
It concludes with a formal problem statement. 

\subsection{Notation}

\begin{table*}[t]
    \centering
    \caption{Table of Symbols and Notation}
    \begin{tabular}{c|c|c}
        Symbol & Description & Defined in \\
         \hline 
         \vspace{-0.15cm}
         $\Q$ & configuration manifold & Sec. \ref{ch-EM1:ssec:background_preliminaries} \\
         \vspace{-0.15cm}
         $T^*\Q$ & cotangent bundle & Sec. \ref{ch-EM1:ssec:background_preliminaries} \\
         \vspace{-0.15cm}
         $\EigN \subset T^*\Q$ & LSM, Eigenmanifold or Rosenberg manifold & Secs. \ref{ch-EM1:ssec:background_preliminaries}, \ref{ch-EM1:ssec: LSM} \\
         \vspace{-0.15cm}
         $\mathcal{R} \subset \EigN$ & generator of an Eigenmanifold or Rosenberg manifold & Sec. \ref{ch-EM1:ssec:background_preliminaries} \\
         \vspace{-0.15cm}
         $\state \in \Q,\,P \in T^*_\state\Q,\, \zeta = (\state,P) \in T^*\Q$ & abstract configuration and momentum variables & Sec. 
         \ref{ch-EM1:ssec:background_preliminaries} \\
         \vspace{-0.15cm}
         $q,p \in \R^n,\, \coordC = (q,p) \in \R^{2n}$ & configuration and momentum in a canonical chart& Sec. 
         \ref{ch-EM1:ssec:background_preliminaries} \\
         \vspace{-0.15cm}
         $H:T^*\Q \rightarrow \R$ & Hamiltonian & Sec. \ref{s1}\\
         \vspace{-0.15cm}
         $X_H \in \X(T^*\Q)$& abstract Hamiltonian vector field & Sec. \ref{ch-EM1:ssec:background_preliminaries} \\
         \vspace{-0.15cm}
         $\lambda_0 \in \mathbb{C}$ & Eigenvalue of linearized system & Sec. \ref{ch-EM1:ssec: LSM} \\
         \vspace{-0.15cm}
         $E_0 \in T_{(\eqstate,0)} T^*\Q$ & Eigenspace of linearized system & Sec. \ref{ch-EM1:ssec: LSM} \\
         \vspace{-0.15cm}
         $\phi_\gamma:[0,1] \rightarrow \Q$ & immersion or immersion of geometric Eigenmode $\stateH$ & Sec. \ref{ch-EM1:ssec:background_preliminaries} \\
         \vspace{-0.15cm}
         $\tau \in \R$ & time-scaling & Sec. \ref{ch-EM1:ssec:background_on_symmetries} \\
         \vspace{-0.15cm}
         $\sigma:\M\rightarrow\M$ & symmetry map & Sec. \ref{ch-EM1:ssec:background_on_symmetries} \\
         \vspace{-0.15cm}
         $S:\R^n \rightarrow \R^n$ & coordinate version of symmetry map & Sec. \ref{ch-EM1:ssec:background_on_symmetries} \\
         \vspace{-0.15cm}
         $\text{Fix}(\sigma)$ & fixed set of a symmetry & Sec. \ref{ch-EM1:ssec:background_on_symmetries} \\
         $\Psi^t_f:\M\rightarrow\M $ & Flow of $f \in \X(\M)$ & Sec. \ref{ch-EM1:ssec:background_on_symmetries}
    \end{tabular}
    \label{tab:notation}
\end{table*}

Smooth manifolds are denoted $\M, \mathcal{N},\mathcal{R}, \mathcal{Q}$. At a point $\state\in \M$, $T_{\state}\M$ denotes the tangent space of $\M$ at $\state$, and $T^*_{\state}\M$ denotes the respective cotangent space. The tangent bundle over $\M$ is $T\M$, and the cotangent bundle is $T^*\M$. The set of vector fields over $\M$ is $\X(\M)$. For a smooth map $\sigma:\M\rightarrow\mathcal{N}$ between smooth manifolds $\M,\mathcal{N}$ we write $\sigma \in C^\infty(\M,\mathcal{N})$. The pushforward of $\sigma$ is $\sigma_*:T_\state\M\rightarrow T_{\sigma(\state)}\mathcal{N}$, and the pullback is $\sigma^*:T^*_{\sigma(\state)}\mathcal{N}\rightarrow T^*_\state\M$. The image of $\sigma$ is denoted $\sigma(\M):=\{\sigma(\state) \;|\; \state\in\M\} \subseteq \mathcal{N}$.\\
%
Additional symbols introduced in the paper are summarized in Table \ref{tab:notation}.

\subsection{Conservative mechanical systems}\label{ch-EM1:ssec:background_preliminaries}

We describe CMs with configuration space $\Q$ by their Hamiltonian dynamics on the cotangent bundle $T^*\Q$. Here a Hamiltonian $H:\M\rightarrow\R$ canonically determines a Hamiltonian vector field $X_H \in \X(T^*\Q)$ (for details on the construction, see e.g., ~\citep{Bullo2004}). Let $\stateH :\R\rightarrow T^*\Q$ denote a solution to the abstract Hamiltonian dynamics
\begin{equation}\label{eq:abstract_hamiltonian_vector_field}
    \dot{\stateH } = X_H(\stateH )\,,\; \stateH (0) = \stateH \,.
\end{equation}
The system~\eqref{eq:abstract_hamiltonian_vector_field} is a CM when $H$ is given by~\eqref{eq:abstract_hamiltonian}.
The corresponding flow of $X_H$ is
\begin{equation}\label{eq:flow}
    \Psi_{X_H}^t:T^*\Q \rightarrow T^*\Q\,,\; \stateH(t) = \Psi_{X_H}^t(\stateH_0)\,, 
\end{equation}
where $\Psi_{X_H}^t$ is defined such that $\stateH(t) = \Psi_{X_H}^t(\stateH_0)$ solves~\eqref{eq:abstract_hamiltonian_vector_field}. It is assumed that $X_H$ is complete. The set of points on $\stateH(t)$ is 
\begin{equation}
   \stateH(\mathbb{R}) := \{\stateH(t)\;|\; t \in \R \}\,. 
\end{equation}
%
In a canonical chart (see Appendix B.1 
in the supplementary material~\cite{arxiV-Link}) 
with coordinates $q,p \in \R^n$, the Hamiltonian dynamics of a CM follow from~\eqref{eq:abstract_hamiltonian},~\eqref{eq:abstract_hamiltonian_vector_field} as:
\begin{align}\label{eq:hamiltonian_ode_q}
    \Dot{q} &= \frac{\partial H}{\partial p}(q,p) = M(q)^{-1}p\,, \\ 
    \Dot{p} &= -\frac{\partial H}{\partial q}(q,p) = -\frac{\partial}{\partial q} \left(\frac{1}{2} p^\top M(q)^{-1} p\right) - \frac{\partial V}{\partial q}(q) \,.
    \label{eq:hamiltonian_ode_p}
\end{align}
We are interested in families of periodic oscillations springing from an isolated, stable equilibrium of~\eqref{eq:hamiltonian_ode_q},~\eqref{eq:hamiltonian_ode_p} (\textit{i.e.}, a point $(\bar{\coord},\bar{p})$ with $\bar{p} = 0$, $\frac{\partial V}{\partial q}(\bar{q}) = 0$, and $\frac{\partial^2 V}{\partial q \partial q}(\bar{q}) > 0$). In a coordinate-free setting we denote this equilibrium as $\bar{\stateH} = (\eqstate,0)$, and in a chart-based setting we choose coordinates such that the equilibrium is $\bar{\coordC} = (\bar{\coord},0) = (0,0)$, with explicit exceptions only in Sec.~\ref{ch-EM1:sec:examples}.

\subsection{Lyapunov subcenter manifolds}\label{ch-EM1:ssec: LSM}
Lyapunov subcenter manifolds (LSMs) are briefly introduced, focusing on applications to systems of the type~\eqref{eq:hamiltonian_ode_q},~\eqref{eq:hamiltonian_ode_p}. For a complete treatment, refer to~\citep{Llave1997,Llave2018} and references therein.

Given a dynamical system in local coordinates $\coordC\in \mathbb{R}^{2n}$
\begin{equation}\label{eq:LSM_ODE}
    \dot{\coordC} = A \coordC + B(\coordC) \,,
\end{equation}
where $A\in \mathbb{R}^{2n \times 2n}$, $B:\mathbb{R}^{2n} \rightarrow \mathbb{R}^{2n}$ is analytic and $B(0) = 0$, $\frac{\partial B}{\partial \coordC}(0) = 0$.

Denote as $(\lambda_1,\hdots, \lambda_{2n})$ the eigenvalues of $A$, and define the linear eigenspace associated with $\lambda_k$ 
\begin{equation}
    E_k := \ker(\lambda_k I-A) \subset \mathbb{C}^{2n}\,.
\end{equation}
 
We repeat \citep[Ass. 2.1 \& Thm. 2.5]{Llave2018} on the existence and uniqueness of LSMs: 

\begin{thm}[Lyapunov subcenter manifolds]\label{thm:LSM}
If it holds that 
\begin{enumerate}[label=(\roman*)]
    \itemsep0em 
    \item A is diagonalizable.
    \item A has a pair of conjugate eigenvalues $\pm i \omega_0$ with zero real part.
    \item The remaining $2n-2$ eigenvalues $\lambda_k$ of A, with $1\leq k \leq 2n-2$, are \textit{non-resonant} with the eigenvalues $\pm i \omega_0$, \textit{i.e.}, it holds for all $k$ that 
    \begin{equation}
        \frac{\lambda_k}{i\omega_0} \notin \mathbb{Z}\,.
    \end{equation}
    \item The system~\eqref{eq:LSM_ODE} has a conserved quantity $H(\coordC(t))$ along integral curves $\coordC(t)$, \textit{i.e.}, 
    \begin{equation}
        \frac{d}{dt} H(\coordC(t)) = 0\,.
    \end{equation}
    Further, $H(0) = 0$, $(\frac{\partial}{\partial \coordC} H)(0) = 0$, and \begin{equation*} Y^\top \left(\left(\frac{\partial^2}{\partial \coordC \partial \coordC^T} H\right)(0)\right)Y > 0 \;, \forall Y \in E_0,\end{equation*}with $E_0$ the eigenspace associated with $\pm i\omega_0$.
\end{enumerate}
Then there exists, locally, a unique 2D submanifold $\EigN\subseteq\mathbb{R}^{2n}$ of periodic trajectories tangent to the eigenspace $E_0$ of $\pm i \omega_0$. 
\end{thm}
This 2D submanifold $\EigN$ is called the Lyapunov subcenter manifold associated with the eigenspace $E_0$ at $z = 0$. Apart from periodicity, there are no conditions on trajectories collected within an LSM (cf. Figure~\ref{fig:overview}a). When the non-resonance condition is not fulfilled, Lyapunov subcenter manifolds can still be found, but they are not unique~\cite{Sijbrand1985}.\\
It follows by repeated application of Theorem~\ref{thm:LSM}:
\begin{cor}\label{cor:n-LSMs}
    If $d$ of the $n$ pairs of eigenvalues $\pm i \omega_k$ of $A$, $1\leq k \leq n$, are mutually non-resonant, and non-resonant with the remaining eigenvalues, then there exist, locally, at least $d$ \textit{unique} families of periodic orbits, with a unique family being tangent to a given eigenspace $E_k$ of $A$. 
\end{cor}
Note that Corollary~\ref{cor:n-LSMs} provides a lower bound since the conditions of Theorem~\ref{thm:LSM} are sufficient, but not necessary.

\subsection{Symmetries of CMs}\label{ch-EM1:ssec:background_on_symmetries}
Discrete symmetries of Hamiltonian dynamics are introduced, alongside a well-known time symmetry of CMs. For a general treatment of symmetries of Hamiltonian systems see e.g.,~\cite{Lamb1998,Alomair2017}.

\begin{defn}[Discrete Symmetry]\label{def:symmetry}
    A symmetry $(\sigma,\tau)$ of the Hamiltonian dynamics $X_H \in \X(T^*\Q)$ consists of a diffeomorphism $\sigma:T^*\Q\rightarrow T^*\Q$ and a scalar $\tau \in \{-1, 1\}$ such that 
    \begin{equation}
        X_H\circ\sigma = \tau \sigma_* X_H \,. \label{eq:symmetry_condition}
    \end{equation}
\end{defn}
\noindent
It is a well-known result that symmetries map solutions to solutions:
\begin{align}
    \Psi^{t}_{X_H} \big(\sigma(\stateH)\big) &= \sigma\big(\Psi^{\tau t}_{X_H}(\stateH)\big) \,, \label{eq:symmetry_mapping}
\end{align}
\textit{i.e.}, for any solution $\stateH(t)$ of the dynamics $X_H$, also $\sigma(\stateH(\tau t))$ is a solution.

The fixed orbit subset $\text{Fix}(\sigma) \subset T^*\Q$ is defined as the union of orbits of ${X_H}$ that are fixed under a symmetry $(\sigma,\tau)$ 
\begin{equation}
    \text{Fix}(\sigma) := \{ \stateH \in \Psi^{\R}_{X_H}(\stateH_0) \; | \; \Psi^{\R}_{X_H}(\stateH_0) = \sigma \Psi^{\R}_{X_H}(\stateH_0)  \}\,.
\end{equation}

In a local chart $(U,\CoordC)$ with $\zeta,\sigma(\zeta) \in U \subset T^*\Q$ and $\CoordC:U\subset T^*\Q \rightarrow\mathbb{R}^{2n}$, the vector field $X_H$ is represented by components $f^i(\coordC) \in \R$, $i \in \{1,\cdots,n\}$, 
and the symmetry as a map $S:\mathbb{R}^{2n}\rightarrow\mathbb{R}^{2n}$ with $S = \CoordC\circ\sigma\circ \CoordC^{-1}$. 
Then condition~\eqref{eq:symmetry_condition} reads 
\begin{equation}
    f^i(S(\coordC)) = \tau \frac{\partial S^i}{\partial \coordC^j} f^j(\coordC) \,.
\end{equation}
If this holds, then for any solution $\coordC(t)$ of the dynamics $\dot{\coordC}^i = f^i(\coordC)$, also $S\big(\coordC(\tau t)\big)$ is a solution.

The equations of motion of a CM with Hamiltonian dynamics~\eqref{eq:abstract_hamiltonian_vector_field} are subject to the time-reversal symmetry $(\state, P, t) \rightarrow (\state, -P, -t)$. As shown by~\cite{Lamb1998}, this follows {more generally for any Hamiltonian system when} $H(\state,P) = H(\state,-P)$. 
Following Definition~\ref{def:symmetry}, this is written as $({\sigma}_{1},\tau_{1})$, with
\begin{equation}\label{eq:tr_sym}
    {\sigma}_{1}\big((\state,P)\big) = (\state, -P)\,, \quad \tau_{1} = -1 \,.
\end{equation}
As a consequence of this symmetry, for any solution $\stateH(t) = (\state(t),P(t))$ of~\eqref{eq:abstract_hamiltonian_vector_field} with Hamiltonian~\eqref{eq:abstract_hamiltonian}, also $\sigma_1(\stateH(\tau_1 t)) = (\state(-t),-P(-t))$ is a solution of~\eqref{eq:abstract_hamiltonian_vector_field}. 

\subsection{Problem statement}
In this article, we aim to answer the following questions:
\vspace{2mm}
\begin{enumerate}
    \item How do the conditions for existence and uniqueness of LSMs translate to conditions on the inertia tensor $M(q)$ and the potential energy $V(q)$? 
    \item What geometric properties of periodic trajectories in LSMs follow from the symmetry~\eqref{eq:tr_sym}?
    \item What other symmetries, corresponding to conditions on $M(q)$ and $V(q)$, result in desirable properties of LSMs?
\end{enumerate}
\section{Lyapunov subcenter manifolds in conservative mechanical systems}\label{ch-EM1:sec:LSMs_in_CMS}

\subsection{Conditions on inertia tensor and potential energy}
We present Theorem~\ref{thm:LSM_mechanics_geometric}, containing conditions on the inertia tensor $M(q)$ and the potential energy $V(q)$ that guarantee the existence and uniqueness of LSMs at an equilibrium $(\eqstate,0)$. This is a minor contribution through its geometric (\textit{i.e.}, chart-free) formalism.

\begin{thm}[LSMs in mechanical systems]\label{thm:LSM_mechanics_geometric}
    Given the Hamiltonian dynamics~\eqref{eq:abstract_hamiltonian_vector_field} on $T^*\Q$ with analytic Hamiltonian~\eqref{eq:abstract_hamiltonian}. If there is a point $\eqstate \in \Q$ where it holds that 
    \begin{enumerate}[label=(\roman*)]
    \itemsep0em
        \item $\text{d}V(\eqstate) = 0$ and $\text{Hess}(V)(\eqstate) > 0$ 
        \item The $(0,2)$-tensor $M(\eqstate)$ is positive definite and symmetric 
        \item The $(1,1)$-tensor $M(\eqstate)^{-1}\text{Hess}(V)(\eqstate)$ is diagonalizable 
        \item The roots $\pm\lambda_0$ of eigenvalues $\lambda_0^2 = -\omega_0^2$ of $M(\eqstate)^{-1}\text{Hess}(V)(\eqstate)$ 
        are \textit{non-resonant} with the remaining roots $\pm\lambda_k$, $ 1\leq k \leq n-1$, in the sense that 
        \begin{equation}
            \frac{\lambda_k}{\lambda_0} \notin \mathbb{Z}\,.
        \end{equation}
    \end{enumerate}
    Then there exists, locally, a unique 2D submanifold $\EigN\subseteq T^*\Q$ of periodic trajectories which is tangent to the eigenspace $E_0 = D_0 \oplus M(\eqstate)D_0 \subset T_{(\bar{\stateH})}T^*\Q$, 
    where $D_0 = \ker(\omega_0^2 I- M(\eqstate)^{-1}\text{Hess}(V)(\eqstate))$. 
\end{thm}
A
\begin{pf}
    The proof relies on showing that the conditions of Theorem~\ref{thm:LSM_mechanics_geometric} are equivalent to those of Theorem~\ref{thm:LSM}. 
    
    First, express the Hamiltonian dynamics of a conservative mechanical system~\eqref{eq:abstract_hamiltonian_vector_field},~\eqref{eq:hamiltonian_ode_p} in canonical coordinates $(q,p)$ such that the equilibrium $(\eqstate,0)$ is mapped to $q = 0$, $p = 0$. Further write $\coordC = (q,p) \in \R^{2n}$. 
    The dynamics then take the form 
    \begin{align}
    \Dot{\coordC} &:= f(\coordC) = \begin{pmatrix}
                    M(q)^{-1}p \\
                    -\frac{\partial}{\partial q} (p^\top M(q)^{-1} p) - \frac{\partial V}{\partial q}(q)
                \end{pmatrix}  \,,
    \end{align}
    which can be rewritten as 
    \begin{equation}
        \Dot{\coordC} = A \coordC + B(\coordC)\,,
    \end{equation}
    with 
    \begin{align}
        A &= \frac{\partial f}{\partial \coordC}\big((0,0)\big) = \begin{bmatrix}
            0 & M(0)^{-1} \\ -\frac{\partial^2 V}{\partial q \partial q}(0) & 0 
        \end{bmatrix} \,, \\
        B(\coordC) &= f(\coordC) - A \coordC \,.
    \end{align}
    The dynamics are of the required form~\eqref{eq:LSM_ODE} for Theorem~\ref{thm:LSM}, since $B(0) = 0$, and $\frac{\partial B}{\partial \coordC}(0) = 0$. $B(\coordC)$ is also analytic, given that $M(q),V(q)$ are analytic.
    This does not rely on the choice of a particular chart: the components of $A$ correspond to the components of a $(1,1)$-tensor on $T^*\Q$, given that $\text{d}V(\eqstate) = 0$, which follows from $(0,0)$ being an equilibrium of the dynamics~\eqref{eq:hamiltonian_ode_q},~\eqref{eq:hamiltonian_ode_p}. 
    
    The eigenvalues $\lambda$ of $A = \begin{bmatrix}
        0 & B \\ C & 0 
    \end{bmatrix}$ are the solution to 
    \begin{equation}
        \det(\lambda I_{2n} - A) = \det(\lambda^2 I_{n} - BC) = 0 \,,
    \end{equation}
    such that $A$ is diagonalizable since the $(1,1)$-tensor $BC = -M(0)^{-1} \frac{\partial^2 V}{\partial q \partial q}(0)$ is diagonalizable. Thus, condition $(i)$ of Theorem~\ref{thm:LSM} is fulfilled.
    
    Further, $\lambda^2 = -\omega^2$ for some positive $\omega \in \R_+$, as $M(0)^{-1} \frac{\partial^2 V}{\partial q \partial q}(0)$ is a positive definite matrix collecting the components of the positive definite tensor $M(\eqstate)^{-1}\text{Hess}(V)(\eqstate)$, and eigenvalues $-\omega_j^2$ of $-M(0)^{-1} \frac{\partial^2 V}{\partial q \partial q}(0)$ correspond to pairs of eigenvalues  
    \begin{equation}
        \lambda_{j,\pm} = \pm i \omega_j
    \end{equation}
    of $A$. Thus, condition $(ii)$ of Theorem~\ref{thm:LSM} is fulfilled.

    The eigenvalues $\lambda_{0,\pm}$ are non-resonant with the remaining $2n - 2$ eigenvalues $\lambda_{k,\pm}$ of $A$ by assumption,
    so condition $(iii)$ of Theorem~\ref{thm:LSM} is fulfilled.

    Finally, a conserved quantity $H\big(\coordC(t)\big)$ fulfilling condition $(iv)$ of Theorem~\ref{thm:LSM} is given by the Hamiltonian $H\big(\coordC(t)\big) = \frac{1}{2} p^\top M(q)^{-1} p + V(q)$\,.
\end{pf}

The conditions of Theorem~\ref{thm:LSM_mechanics_geometric} are chart-invariant: the equilibrium is guaranteed to be of the form $(\eqstate,0)$ in any canonical chart, $M(q)$ is positive definite and symmetric in any chart, and eigenvalues of the component matrix $M(0)^{-1}\frac{\partial^2 V}{\partial q \partial q}(0)$ of the $(1,1)$-tensor $M(\eqstate)^{-1}\text{Hess}(V)(\eqstate)$ are chart-invariant. 


\subsection{Propagation of symmetries}

Theorem~\ref{thm:symmetry_propagation_system_to_LSM} presented hereafter, is a novel result that is essential to derive results in the remainder of this section. Given a symmetry $(\sigma,\tau)$ of $X_H$, it presents conditions for an LSM $\EigN$ and trajectories contained within it to be fixed sets of the symmetry.

\begin{thm}[Symmetric LSMs]\label{thm:symmetry_propagation_system_to_LSM}
    Given dynamics~\eqref{eq:abstract_hamiltonian_vector_field} where~\eqref{eq:abstract_hamiltonian} fulfills the conditions of Theorem~\ref{thm:LSM_mechanics_geometric}, and further given a symmetry $(\sigma,\tau)$ of $X_H$ that fulfills Definition~\ref{def:symmetry}. Denote by $\EigN$ the Lyapunov subcenter manifold tangent to\footnote{To be precise, with $\i:\EigN\rightarrow T^*\Q$ the natural embedding of $\EigN$ into $T^*\Q$, then $\EigN$ is tangent to $E_0$ if $\i_*T_{(\eqstate,0)}\EigN \subset E_0$.}
    the Eigenspace $E_0$. If it holds that 
    \begin{equation}
        \sigma((\eqstate,0)) = (\eqstate,0)\,,
    \end{equation}
    and
    \begin{equation}\label{condition:symmetry_propagation_system_to_LSM}
        E_0 = \sigma_* E_0\,,
    \end{equation}
    then it also holds that 
    \begin{equation}
        \EigN = \sigma \big(\EigN\big)\,,
    \end{equation}
    and any evolution $\stateH(\mathbb{R}) \in \EigN$ satisfies 
    \begin{equation}
        \stateH(\mathbb{R}) = \sigma(\stateH(\mathbb{R}))\,.
    \end{equation}
\end{thm}

\begin{pf}
    Begin by showing that $\EigN = \sigma(\EigN)$. The presence of the symmetry $\sigma$ implies that both $\EigN$ and $\EigN' := \sigma(\EigN)$ are LSMs, and condition $\sigma(0) = 0$ implies that $0 \in \EigN$ and $0 \in \EigN'$, \textit{i.e.}, both are LSMs about the equilibrium $\coordC = 0$. Further, ${\sigma}_*E_0 = E_0$ implies that both $\EigN$ and $\EigN'$ are tangent to the eigenspace $E_0$. Theorem~\ref{thm:LSM_mechanics_geometric} holds, such that $\EigN$ is the unique LSM tangent to $E_0$, and hence it must be that $\EigN = \EigN'$.

    To show that also $\stateH(\R) = \sigma(\stateH(\R))$, introduce coordinates $(U,X)$ on $\EigN$, \textit{i.e.}, $U \subset \EigN$ is a neighborhood of $\coordC = 0$ and $X:\EigN \rightarrow \R^2$. Further, choose polar coordinates $X$ and a time scaling $\tau:\EigN \rightarrow \mathbb{R}_+$ such that $X(\stateH(t)) = (r, t)$ with $r$ constant for a given oscillation $\stateH \in \EigN$. 
    Then $\sigma:\EigN \rightarrow\EigN$ translates into coordinates as
    \begin{equation}
        S(r,t) = (\alpha(r),t + c)\,,
    \end{equation}
    Since it must map solutions to solutions.
    The condition 
    \begin{equation}
        \sigma\circ \sigma = \text{id}_{\EigN}\,, \quad \sigma(0) = 0
    \end{equation}
    translates into coordinates as 
    \begin{equation}
        S(S(r,\theta)) = (r,\theta)\,, \quad \Phi(0,\theta) = 0\,,
    \end{equation}
    or in terms of $\alpha$
    \begin{equation}
        \alpha(\alpha(r)) = r\,, \quad \alpha(0) = 0\,,
    \end{equation}
    With the additional constraint that $\alpha(r) \geq 0$, this has the unique continuous solution $\alpha(r) = r$. Thus
    \begin{equation}
        S(r,t) = (r, t+c)\,,
    \end{equation}
    necessarily maps solutions into themselves, \textit{i.e.}, $\sigma(\stateH(\R)) = \stateH(\R)$, completing the proof.
\end{pf}


\begin{rem}
    In a local chart, condition~\eqref{condition:symmetry_propagation_system_to_LSM} requires that $E_0$ is an eigenspace of the Jacobian $\frac{\partial S}{\partial \coordC}_{|\coordC=0}$.
\end{rem}

\section{Properties of Lyapunov subcenter manifolds induced by symmetry}\label{ch-EM1:sec:theorems}

This section presents the main results, combining Theorems~\ref{thm:LSM_mechanics_geometric},~\ref{thm:symmetry_propagation_system_to_LSM} with the symmetry~\eqref{eq:tr_sym} to derive properties of self-symmetric LSMs in CMs. It further presents a general class of spatial symmetries that are equivalent to conditions on $M(q),V(q)$ and that lead to desirable control-theoretic properties of LSMs in CMs.

\subsection{Time-symmetric Lyapunov subcenter manifolds}

Properties of time-symmetric Lyapunov subcenter manifold are collected in a theorem:
\begin{thm}[Eigenmanifolds]\label{thm:Eigenmanifold_Condition}
    Given the conditions of Theorem~\ref{thm:LSM_mechanics_geometric}, and with $\pi(\state,p) = \state$ 
    the projection of $T^*\Q$ to $\Q$. Then the unique LSM $\EigN\subset T^*\Q$ associated with the eigenspace $E_0 \subset T_{(\bar{q},0)}T^*\Q$ satisfies the following properties: 
    \begin{enumerate}
        \item Immersion: For any trajectory $\zeta(t) \in \EigN$ there is an immersion\footnote{The interval $[0,1]$ is a manifold with boundary, for which the notion of a immersion is not immediately defined. An immersion $\sigma:\mathcal{P}\rightarrow\mathcal{N}$ for a manifold $\mathcal{P}$ with boundary $\partial\mathcal{P}$ is defined such that $\sigma_{|\partial \mathcal{P}}$ is an immersion of $\partial \mathcal{P}$. For $\mathcal{P} = [0,1]$, the boundary is $\partial\mathcal{P} = \{0,1\}$. Since this boundary is 0-dimensional it has no tangent-space, any map with injective differential on $(0,1)$ can be defined as an immersion for the interval $[0,1]$.}
        $\phi_{\zeta}:[0,1]\rightarrow \Q$ such that 
        \begin{equation}\label{eq:immersion-eigenmode}
            \phi_{\zeta}([0,1]) = \pi(\zeta(\R))
        \end{equation}
        \item Embedding: For a neighborhood $\EigN'\subset \EigN$ around the equilibrium, $\phi_{\zeta}:[0,1]\rightarrow \Q$ is an embedding. When $\dim \Q = 2$, then $\EigN'= \EigN$.
        \item Brake points: A trajectory $\zeta(t) \in \EigN$ encounters two points with $p(t) = 0$ over a period of oscillation $T_{\zeta}$, at the configurations $\phi_{\zeta}(0)$ and $\phi_{\zeta}(1)$, and they are a time $T_{\zeta}/2$ apart.
        \item Generator: The configurations $\phi_{\zeta}(0)$ and $\phi_{\zeta}(1)$ of the $\zeta\in\EigN$ lie on a connected, 1D-submanifold $\mathcal{R}\subset \EigN$ that passes through the equilibrium.
    \end{enumerate}
\end{thm}

\begin{pf}
    See Appendix~\ref{ap-ssec:time-symmetry}. 
\end{pf}

The existence of an immersion $\phi_{\zeta}:[0,1]\rightarrow \Q$ satisfying~\eqref{eq:immersion-eigenmode} means that configurations $\state(t)$ evolve on a 1D line that may self-intersect. When $\phi_{\zeta}:[0,1]\rightarrow \Q$ is an embedding, then the evolution of the configuration $\state(t)$ does not self-intersect.\\ 
The existence of points with $P(t) = 0$ has immediate implications for application, suggesting that motions collected on LSMs of CMs are interesting for e.g., pick-and-place-like tasks where a manipulator comes to a full stop at two points, but not immediately apply to locomotion, where the system is in continuous motion. This property also allows to reuse insights from numerical continuation methods of nonlinear normal modes~\cite{Kerschen2009} which search for new periodic motions only adjusting the initial configuration $\state_0$ and period $T_\zeta$, but not the initial momentum $P_0$. This directly generalizes to the numerical continuation of unique LSMs in any type of CM.
Finally, the fact that configurations $\phi_{\zeta}(0)$ and $\phi_{\zeta}(1)$ lie on a connected, 1D-submanifold $\mathcal{R}\subset \EigN$ is convenient for parameterizing motions on time-symmetric LSMs in terms of these 1D-submanifolds. \\
We make a number of definitions inspired by Theorem~\ref{thm:Eigenmanifold_Condition}. The definitions closely reflect standard definitions~\cite{Gluck1983},~\cite{Kerschen2009}, recent definitions in~\citep[Sec. 7]{AlbuSchaeffer2020}. Compared to~\citep[Sec. 7]{AlbuSchaeffer2020}, they follow from first principles. 

\begin{defn}[Weak Geometric Eigenmodes]\label{def:eigenmode}
   A \textbf{geometric eigenmode} is a $T$-periodic oscillation $\stateH:\mathbb{R}\rightarrow T^*\Q$, such that 
   $\pi (\stateH(\mathbb{R}))\cong [0,1]$, \textit{i.e.}, there is an embedding $\phi_{\stateH}: [0,1] \rightarrow \Q$ such that $\phi_{\stateH}([0,1]) = \pi (\stateH(\mathbb{R}))$.
   Instead, $\stateH(t)$ is a \textbf{weak geometric eigenmode} if $\phi_{\stateH}: [0,1] \rightarrow \Q$ is an immersion and $\phi_{\stateH}(0)\neq \phi_{\stateH}(1)$. 
   The points $\phi_{\stateH}(0)$ and $\phi_{\stateH}(1)$ are called \textbf{brake points}~\cite{Gluck1983}.
\end{defn}

\begin{defn}[Weak Eigenmanifold]\label{def:eigenmanifold}
An LSM $\EigN\subseteq T^*Q$ is an \textbf{(weak) Eigenmanifold} if all trajectories contained within it are \textbf{(weak) geometric eigenmodes}.
\end{defn}

\begin{defn}[Generator]
    The set $\mathcal{R} \subset \EigN$ of brake points is called the generator of $\EigN$.
\end{defn}
The generator satisfies that $\EigN = \Psi^\R_{X_H}(\mathcal{R})$, \textit{i.e.}, it generates $\EigN$ by forward evolution of the dynamics $X_H$.\\
Prior definitions of Eigenmanifolds (\citep[Sec. 7]{AlbuSchaeffer2020}) explicitly assumed the existence of the generator. It follows from Theorem~\ref{thm:Eigenmanifold_Condition} that this assumption is not necessary, since the generator exists for any family of weak geometric eigenmodes.\\
A lower bound for the number of Eigenmanifolds follows from Corollary~\ref{cor:n-LSMs} and Theorem~\ref{thm:Eigenmanifold_Condition}:
\begin{cor}[Number of Eigenmanifolds]\label{cr:n-Eigenmanifolds}
     Consider a system of the type~\eqref{eq:hamiltonian_ode_q},~\eqref{eq:hamiltonian_ode_p}, with an equilbrium $\bar{q}$ a minimum of the potential $V(q)$. If $m \leq n$ of the $n$ eigenvalues $\omega^2_k$ of the $(1,1)$-tensor $M(\bar{q})^{-1}\frac{\partial^2 V}{\partial q^T \partial q}(\bar{q})$ are mutually non-resonant, and non-resonant with the remaining eigenvalues, then there exist, locally, at least $m$ unique Eigenmanifolds. 
\end{cor}

\subsection{Time and spatially symmetric Lyapunov subcenter manifolds}
Here it is shown that subjecting the Lyapunov subcenter manifold to an additional spatial symmetry causes all trajectories to pass through the equilibrium configuration.\\
With $\bar{\state} = \arg\min V(\state) \in \Q$ the equilibrium configuration of~\eqref{eq:abstract_hamiltonian_vector_field}, let $\varphi:\Q\rightarrow\Q$ be a diffeomorphism such that:
\begin{align}
    \varphi(\bar{\state}) &= \bar{\state} \,, \label{eq:eq-sym-1} \\
    \varphi\circ\varphi &= \text{id}_{\Q} \,, \label{eq:eq-sym-2}\\
    \varphi_*(\bar{\state}) &= -\text{id}_{T_{\bar{\state}}\Q} \label{eq:eq-sym-3}\,.
\end{align}
Denote a class of equivariant symmetries by $\big(\state,P,t\big) \rightarrow \big(\varphi(\state),(\varphi^{-1})^*P,t\big)$, and write it as $(\sigma_2,\tau_2)$ with 
\begin{equation}\label{eq:eq_sym}
    \sigma_2((\state,P)) = (\varphi(\state), (\varphi^{-1})^*P)\,, \quad \tau_2 = 1\,. 
\end{equation}
Symmetries of this class apply to the CM~\eqref{eq:abstract_hamiltonian_vector_field} with Hamiltonian~\eqref{eq:abstract_hamiltonian} when $M = \varphi^* M$ and $V = V \circ \varphi$, as is shown in Appendix~\ref{ch-EM1:ap:symmetry}.\\
%
A local chart exists (see Appendix~\ref{ch-EM1:ap:symmetry}) that maps $\bar{\state}$ to $q = 0$, and where the symmetry reads $(q,p,t)\rightarrow(-q,-p,t)$, such that $(\sigma_2,\tau_2)$ can be written as $(S_2,\tau_2)$ with
\begin{equation} \label{eq:equi_sym} 
    S_2((q,p)) = (- q, - p)\,, \quad \tau_2 = 1 \,.
\end{equation}
Refer to these coordinates as $\varphi$-equivariant coordinates. In these coordinates the conditions $M = \varphi^* M$ and $V = V \circ \varphi$ read 
\begin{align}
    M(q) &= M(-q)\,, \\
    V(q) &= V(-q)\,.
\end{align}

\begin{thm}[Rosenberg manifolds]\label{Theorem:1}
Given the conditions of Theorem~\ref{thm:LSM_mechanics_geometric}, and additionally $\varphi:\Q\rightarrow\Q$ satisfying~\eqref{eq:eq-sym-1} to~\eqref{eq:eq-sym-3} and
\begin{enumerate}[label=(\roman*)]
    \item $V = V\circ\varphi$ and $M = \varphi^* M$
    \item The equilibrium $\eqstate$ is the only fixed point of $\varphi$ such that $V(\eqstate) \leq V(Q)$ for $Q \in \EigN$
\end{enumerate}
then a unique LSM $\EigN\subset T^*\Q$ associated with the eigenspace $E_0 \subset T_{(\bar{q},0)}T^*\Q$ has, in addition to the properties in Theorem~\ref{thm:Eigenmanifold_Condition}, the properties  
\begin{enumerate}
    \item Equilibrium: all trajectories $\zeta(t) \in \EigN$ pass through the equilibrium configuration $\eqstate \in \pi(\zeta(\R))$.
    \item Time to equilibrium: with $T_{\zeta}$ the period of $\zeta(t)$, then the points with $\pi(\zeta(t)) = \eqstate$ are a time $T_{\zeta}/2$ apart and a time $T_{\zeta}/4$ from the brake points. 
\end{enumerate}
\end{thm}
\begin{pf}
See Appendix~\ref{proof:rosenberg_manifold}.
\end{pf}
%
%
We make two definitions based on Theorem~\ref{Theorem:1}:
\begin{defn}[Weak Geometric Rosenberg Mode]\label{def:rosenberg_mode}
    A $T$-periodic oscillation $\stateH:\mathbb{R}\rightarrow T^*\Q$ is a \textbf{(weak) geometric Rosenberg mode} with respect to an equilibrium $\eqstate$ if it is both a (weak) geometric eigenmode, and $\eqstate \in \pi_Q (\stateH(\mathbb{R}))$, \textit{i.e.}, if $\stateH$ passes through the equilibrium configuration $\eqstate$ for some momentum $P$.
\end{defn}
%
\begin{defn}[Weak Rosenberg Manifold]\label{def:rosenberg_geo}
	An LSM $\EigN\subseteq T^*Q$ is a \textbf{(weak) Rosenberg manifold} if all trajectories contained within it are (weak) Rosenberg modes.
\end{defn}

\begin{rem}
Definitions~\ref{def:rosenberg_mode} and~\ref{def:rosenberg_geo} 
are equivalent to Definition 8 in~\cite{AlbuSchaeffer2020}, but disagree with other definitions in the literature:
Rosenberg~\cite{Rosenberg1966} defined Rosenberg modes as modes for which the coordinates $q_1(t),\cdots,q_n(t)$ oscillate in unison. Oscillation in unison is ill-defined in a coordinate-free setting: such coordinates can be found for any geometric eigenmode, and the property is not invariant under changes of coordinates. Definition~\ref{def:rosenberg_mode} remedies this and generalizes~\cite{Rosenberg1966}:  oscillation in unison implies that modes pass through the unique equilibrium at $q = 0$, and motion through the equilibrium configuration is a coordinate independent feature.
Peeters et al.~\cite{Peeters2009} defined extended Rosenberg modes, which allow for coordinates to oscillate asynchronously. This definition can be read to include arbitrary periodic modes, geometric eigenmodes and extended Rosenberg modes, making it non-specific. Definition~\ref{def:rosenberg_mode} excludes arbitrary periodic oscillations and makes geometric Rosenberg modes a subset of geometric eigenmodes.
\end{rem}

\begin{figure}[t]
  \centering
  \begin{tabular}{c c}
  \begin{subfigure}[t]{.5\columnwidth}
      \def\svgwidth{12em}\scriptsize
\begingroup%
  \makeatletter%
  \providecommand\color[2][]{%
    \errmessage{(Inkscape) Color is used for the text in Inkscape, but the package 'color.sty' is not loaded}%
    \renewcommand\color[2][]{}%
  }%
  \providecommand\transparent[1]{%
    \errmessage{(Inkscape) Transparency is used (non-zero) for the text in Inkscape, but the package 'transparent.sty' is not loaded}%
    \renewcommand\transparent[1]{}%
  }%
  \providecommand\rotatebox[2]{#2}%
  \newcommand*\fsize{\dimexpr\f@size pt\relax}%
  \newcommand*\lineheight[1]{\fontsize{\fsize}{#1\fsize}\selectfont}%
  \ifx\svgwidth\undefined%
    \setlength{\unitlength}{944.04483727bp}%
    \ifx\svgscale\undefined%
      \relax%
    \else%
      \setlength{\unitlength}{\unitlength * \real{\svgscale}}%
    \fi%
  \else%
    \setlength{\unitlength}{\svgwidth}%
  \fi%
  \global\let\svgwidth\undefined%
  \global\let\svgscale\undefined%
  \makeatother%
  \begin{picture}(1,0.61046935)%
    \lineheight{1}%
    \setlength\tabcolsep{0pt}%
    \put(0.72027186,0.0046705){\color[rgb]{0,0,0}\makebox(0,0)[lt]{\lineheight{1.25}\smash{\begin{tabular}[t]{l}$\Q$\end{tabular}}}}%
    \put(0,0){\includegraphics[width=\unitlength,page=1]{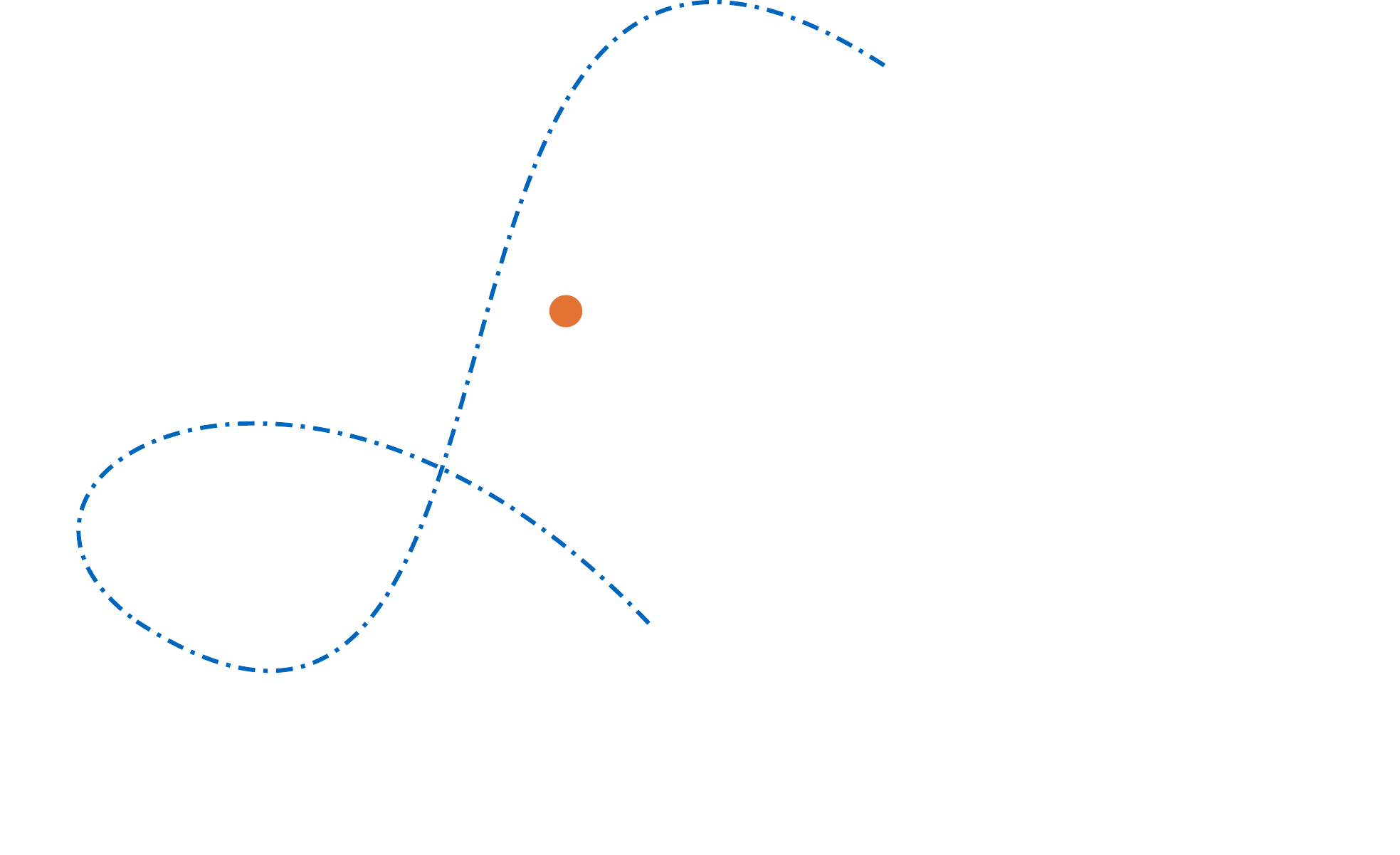}}%
    \put(0.41427136,0.3193896){\color[rgb]{0,0,0}\makebox(0,0)[lt]{\lineheight{1.25}\smash{\begin{tabular}[t]{l}$\eqstate$\end{tabular}}}}%
    \put(0,0){\includegraphics[width=\unitlength,page=2]{weak_geometric_eigenmode_svg-tex.pdf}}%
    \put(0.6582308,0.54286977){\color[rgb]{0,0,0}\makebox(0,0)[lt]{\lineheight{1.25}\smash{\begin{tabular}[t]{l}$\left(\state, 0\right)$\end{tabular}}}}%
    \put(-0.00263783,0.4822284){\color[rgb]{0,0,0}\makebox(0,0)[lt]{\lineheight{1.25}\smash{\begin{tabular}[t]{l}$\left(\state(t), P(t)\right)$\end{tabular}}}}%
  \end{picture}%
\endgroup%

      \subcaption{Weak geometric eigenmode}\label{fig:weak_geometric_eigenmode}
  \end{subfigure}
  &
  \begin{subfigure}[t]{.4\columnwidth}
      \def\svgwidth{12em}\scriptsize
\begingroup%
  \makeatletter%
  \providecommand\color[2][]{%
    \errmessage{(Inkscape) Color is used for the text in Inkscape, but the package 'color.sty' is not loaded}%
    \renewcommand\color[2][]{}%
  }%
  \providecommand\transparent[1]{%
    \errmessage{(Inkscape) Transparency is used (non-zero) for the text in Inkscape, but the package 'transparent.sty' is not loaded}%
    \renewcommand\transparent[1]{}%
  }%
  \providecommand\rotatebox[2]{#2}%
  \newcommand*\fsize{\dimexpr\f@size pt\relax}%
  \newcommand*\lineheight[1]{\fontsize{\fsize}{#1\fsize}\selectfont}%
  \ifx\svgwidth\undefined%
    \setlength{\unitlength}{814.84531862bp}%
    \ifx\svgscale\undefined%
      \relax%
    \else%
      \setlength{\unitlength}{\unitlength * \real{\svgscale}}%
    \fi%
  \else%
    \setlength{\unitlength}{\svgwidth}%
  \fi%
  \global\let\svgwidth\undefined%
  \global\let\svgscale\undefined%
  \makeatother%
  \begin{picture}(1,0.80302781)%
    \lineheight{1}%
    \setlength\tabcolsep{0pt}%
    \put(0.81461626,0.00541105){\color[rgb]{0,0,0}\makebox(0,0)[lt]{\lineheight{1.25}\smash{\begin{tabular}[t]{l}$\Q$\end{tabular}}}}%
    \put(0,0){\includegraphics[width=\unitlength,page=1]{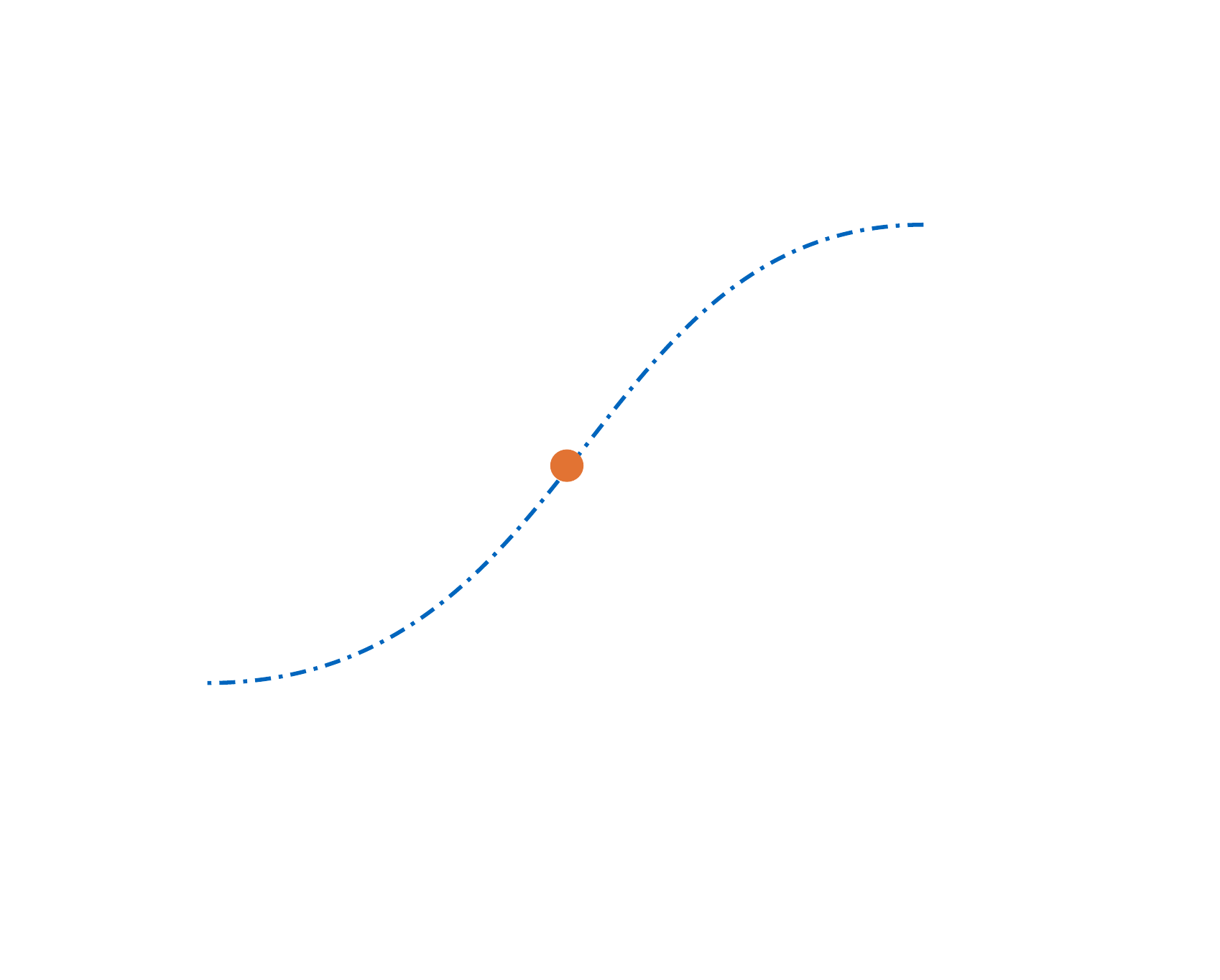}}%
    \put(0.47639627,0.33776164){\color[rgb]{0,0,0}\makebox(0,0)[lt]{\lineheight{1.25}\smash{\begin{tabular}[t]{l}$\eqstate$\end{tabular}}}}%
    \put(0,0){\includegraphics[width=\unitlength,page=2]{weak_geometric_rosenbergmode2_svg-tex.pdf}}%
    \put(0.60404086,0.77505576){\color[rgb]{0,0,0}\makebox(0,0)[lt]{\lineheight{1.25}\smash{\begin{tabular}[t]{l}$\left(\state, 0\right)$\end{tabular}}}}%
    \put(-0.00305607,0.0316292){\color[rgb]{0,0,0}\makebox(0,0)[lt]{\lineheight{1.25}\smash{\begin{tabular}[t]{l}$\left(\varphi(\state), 0\right)$\end{tabular}}}}%
    \put(0.03139839,0.52052151){\color[rgb]{0,0,0}\makebox(0,0)[lt]{\lineheight{1.25}\smash{\begin{tabular}[t]{l}$\left(\state(t), P(t)\right)$\end{tabular}}}}%
  \end{picture}%
\endgroup%

      \subcaption{Weak geometric Rosenberg mode}\label{fig:weak_geometric_rosenberg_mode}
  \end{subfigure}
  \end{tabular}
  \caption{\small Example of a weak geometric eigenmode in panel (a) and a $\varphi$-equivariant weak geometric Rosenberg mode in panel (b). Extending their non-weak counterparts, these allow for self-intersections in configuration space. } 
  \label{fig:weak_geometric_modes}
\end{figure}

\section{Examples}\label{ch-EM1:sec:examples}
Here, aspects of Theorem~\ref{thm:Eigenmanifold_Condition}, Corollary~\ref{cr:n-Eigenmanifolds} and Theorem~\ref{Theorem:1} are highlighted numerically. Existence of Eigenmanifolds, the number of Eigenmanifolds and the existence of Rosenberg manifolds are presented for examples with increasing complexity. 

\subsection{2-DoF systems}
\begin{figure}[t]
  \centering
  \begin{tabular}{c c}
  \begin{subfigure}[t]{.235\textwidth}
    \def\svgwidth{10em}
\begingroup%
  \makeatletter%
  \providecommand\color[2][]{%
    \errmessage{(Inkscape) Color is used for the text in Inkscape, but the package 'color.sty' is not loaded}%
    \renewcommand\color[2][]{}%
  }%
  \providecommand\transparent[1]{%
    \errmessage{(Inkscape) Transparency is used (non-zero) for the text in Inkscape, but the package 'transparent.sty' is not loaded}%
    \renewcommand\transparent[1]{}%
  }%
  \providecommand\rotatebox[2]{#2}%
  \newcommand*\fsize{\dimexpr\f@size pt\relax}%
  \newcommand*\lineheight[1]{\fontsize{\fsize}{#1\fsize}\selectfont}%
  \ifx\svgwidth\undefined%
    \setlength{\unitlength}{534.81315827bp}%
    \ifx\svgscale\undefined%
      \relax%
    \else%
      \setlength{\unitlength}{\unitlength * \real{\svgscale}}%
    \fi%
  \else%
    \setlength{\unitlength}{\svgwidth}%
  \fi%
  \global\let\svgwidth\undefined%
  \global\let\svgscale\undefined%
  \makeatother%
  \begin{picture}(1,1.03548634)%
    \lineheight{1}%
    \setlength\tabcolsep{0pt}%
    \put(0,0){\includegraphics[width=\unitlength,page=1]{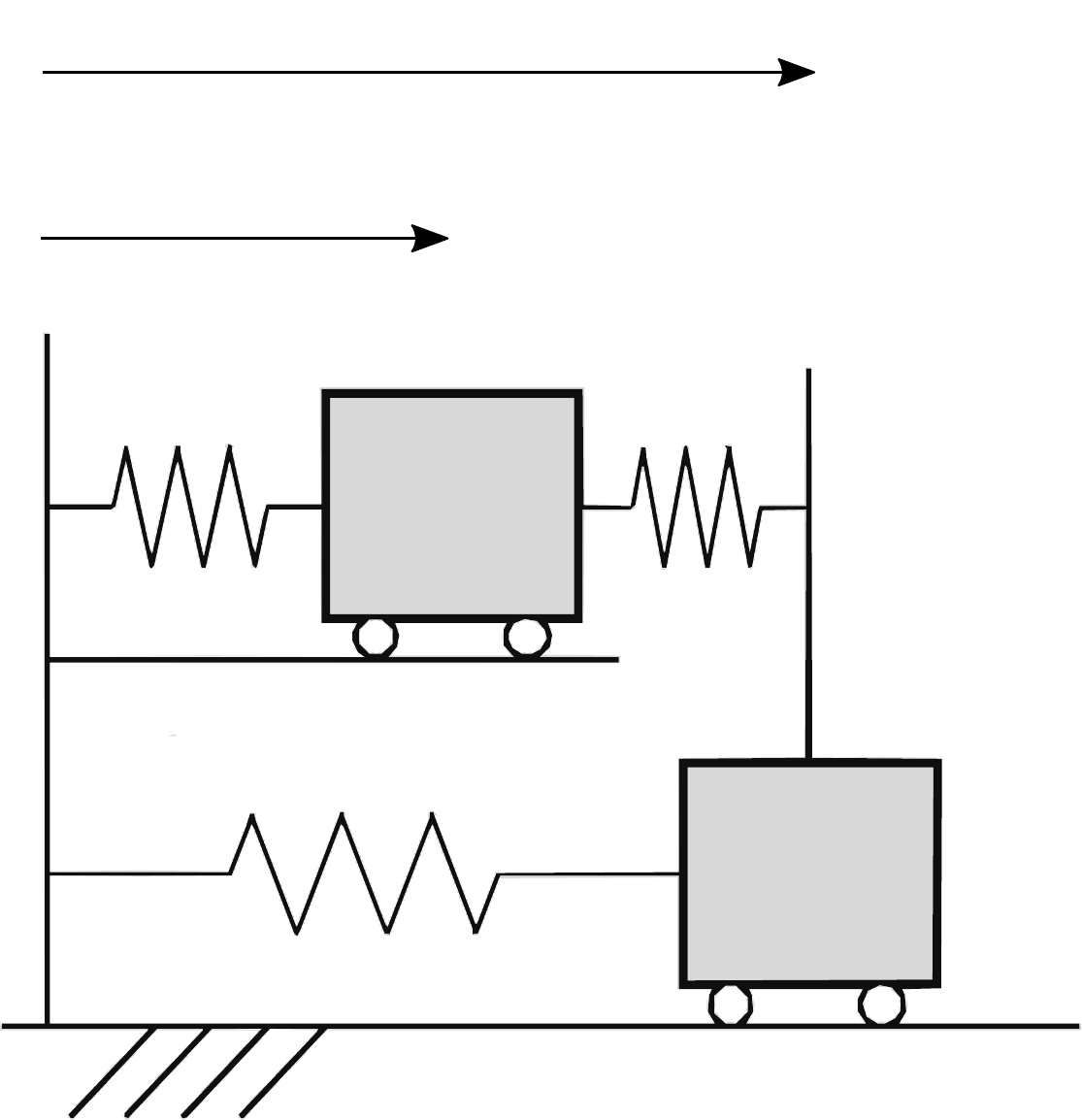}}%
    \put(0.35590997,0.54380976){\color[rgb]{0,0,0}\makebox(0,0)[lt]{\lineheight{1.25}\smash{\begin{tabular}[t]{l}$m$\end{tabular}}}}%
    \put(0.69499934,0.2083318){\color[rgb]{0,0,0}\makebox(0,0)[lt]{\lineheight{1.25}\smash{\begin{tabular}[t]{l}$m$\end{tabular}}}}%
    \put(0.12059427,0.85244649){\color[rgb]{0,0,0}\makebox(0,0)[lt]{\lineheight{1.25}\smash{\begin{tabular}[t]{l}$q_1$\end{tabular}}}}%
    \put(0.28405065,0.99655241){\color[rgb]{0,0,0}\makebox(0,0)[lt]{\lineheight{1.25}\smash{\begin{tabular}[t]{l}$q_2$\end{tabular}}}}%
  \end{picture}%
\endgroup%

      \subcaption{Euclidean system}\label{fig:eucl_sys}
  \end{subfigure}
  &
  \begin{subfigure}[t]{.24\textwidth}
    \def\svgwidth{10em}
\begingroup%
  \makeatletter%
  \providecommand\color[2][]{%
    \errmessage{(Inkscape) Color is used for the text in Inkscape, but the package 'color.sty' is not loaded}%
    \renewcommand\color[2][]{}%
  }%
  \providecommand\transparent[1]{%
    \errmessage{(Inkscape) Transparency is used (non-zero) for the text in Inkscape, but the package 'transparent.sty' is not loaded}%
    \renewcommand\transparent[1]{}%
  }%
  \providecommand\rotatebox[2]{#2}%
  \newcommand*\fsize{\dimexpr\f@size pt\relax}%
  \newcommand*\lineheight[1]{\fontsize{\fsize}{#1\fsize}\selectfont}%
  \ifx\svgwidth\undefined%
    \setlength{\unitlength}{467.82552928bp}%
    \ifx\svgscale\undefined%
      \relax%
    \else%
      \setlength{\unitlength}{\unitlength * \real{\svgscale}}%
    \fi%
  \else%
    \setlength{\unitlength}{\svgwidth}%
  \fi%
  \global\let\svgwidth\undefined%
  \global\let\svgscale\undefined%
  \makeatother%
  \begin{picture}(1,1.16607588)%
    \lineheight{1}%
    \setlength\tabcolsep{0pt}%
    \put(0,0){\includegraphics[width=\unitlength,page=1]{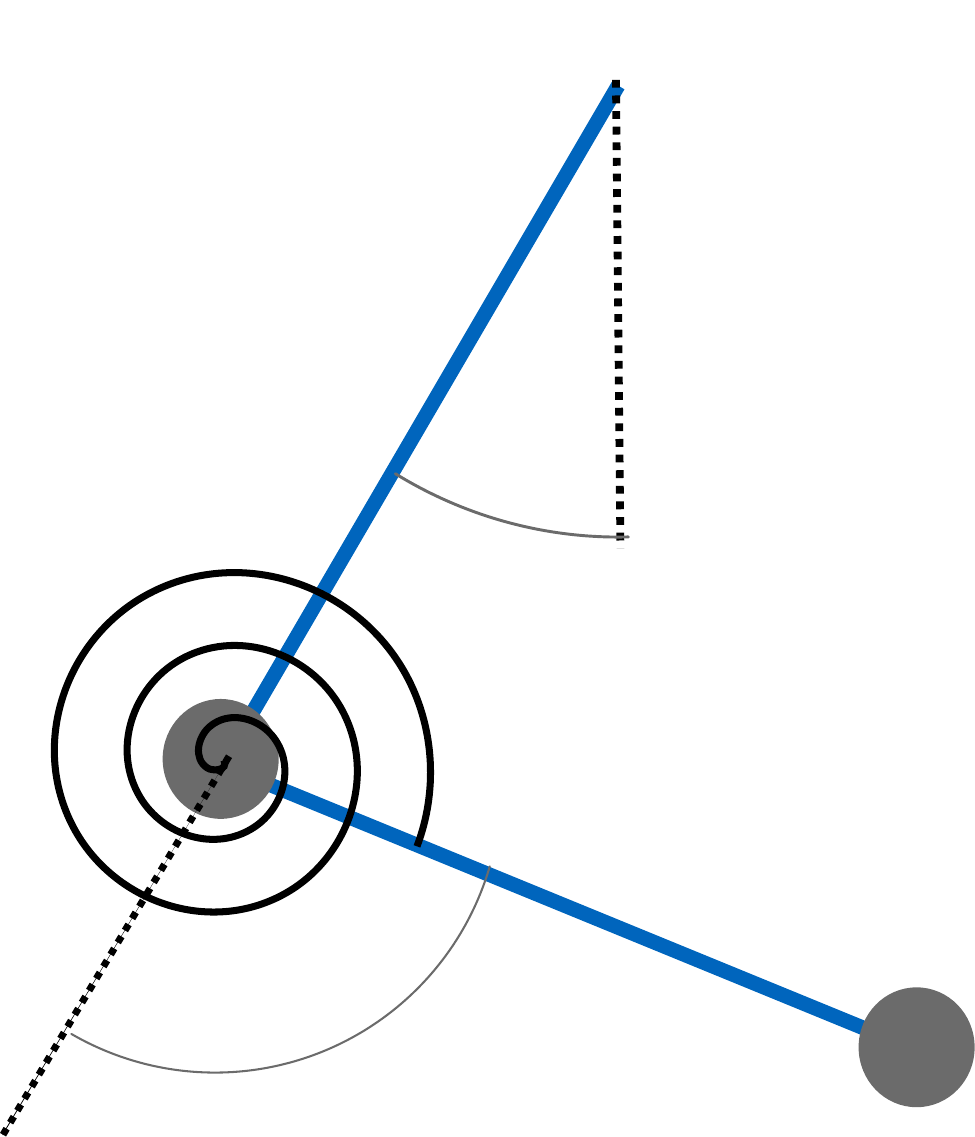}}%
    \put(0.50140437,0.68361838){\color[rgb]{0,0,0}\makebox(0,0)[lt]{\lineheight{1.25}\smash{\begin{tabular}[t]{l}$q_1$\end{tabular}}}}%
    \put(0.24033561,0.13840837){\color[rgb]{0,0,0}\makebox(0,0)[lt]{\lineheight{1.25}\smash{\begin{tabular}[t]{l}$q_2$\end{tabular}}}}%
    \put(0,0){\includegraphics[width=\unitlength,page=2]{DP_svg-tex.pdf}}%
  \end{picture}%
\endgroup%

      \subcaption{Non-Euclidean system}\label{fig:double_pendulum}
  \end{subfigure}
  \end{tabular}
  \caption{\small Examples with $\dim(\Q) = 2$: the system of two masses in panel (a) has a Euclidean configuration space, $\R^2$ with a constant inertia tensor. The double pendulum in panel (b) has a non-Euclidean configuration space, $T^2$ with a non-zero curvature, where the associated inertia tensor is non-constant in any choice of coordinates.} 
  \label{fig:example_systems}
\end{figure}

We begin by describing our class of example systems: coupled masses with a linear motion profile (Fig.~\ref{fig:eucl_sys}) and a double pendulum (Fig.~\ref{fig:double_pendulum}), exerted to different potential fields. 
The coupled masses in Fig.~\ref{fig:eucl_sys} evolve on a Euclidean configuration manifold $\Q_C = \R^2$ with curvature-free metric tensor 
\begin{equation}
    M_C(q) =     
    \begin{bmatrix}
    m & 0 \\
    0 & m
    \end{bmatrix} \,.
\end{equation}
The double pendulum in Fig.~\ref{fig:double_pendulum} evolves on a non-Euclidean configuration manifold $\Q_{DP} = T^2$ with metric tensor 
\begin{equation} \label{s4eq_MDP}
    M_{DP}(q) = 
    \begin{bmatrix}
    I+3 d^2 m \big(1 +2 \cos(q_2)) & d^2 m \big(1+\cos(q_2)\big) \\ d^2 m \big(1+\cos(q_2)\big) & I+d^2 m
    \end{bmatrix} \,.
\end{equation} 
Note that both $M_C(q)$ and $M_{DP}(q)$ fulfill condition (i) of Theorem~\ref{Theorem:1}, i.e., $M_C(q) = M_C(-q)$ and $M_{DP}(q) = M_{DP}(-q)$.

\begin{figure}[t]
    \centering
      \begin{tabular}{c c c}
        \begin{subfigure}[b]{.2\columnwidth}
          \def\svgwidth{5.5em}
\begingroup%
  \makeatletter%
  \providecommand\color[2][]{%
    \errmessage{(Inkscape) Color is used for the text in Inkscape, but the package 'color.sty' is not loaded}%
    \renewcommand\color[2][]{}%
  }%
  \providecommand\transparent[1]{%
    \errmessage{(Inkscape) Transparency is used (non-zero) for the text in Inkscape, but the package 'transparent.sty' is not loaded}%
    \renewcommand\transparent[1]{}%
  }%
  \providecommand\rotatebox[2]{#2}%
  \newcommand*\fsize{\dimexpr\f@size pt\relax}%
  \newcommand*\lineheight[1]{\fontsize{\fsize}{#1\fsize}\selectfont}%
  \ifx\svgwidth\undefined%
    \setlength{\unitlength}{440.98596693bp}%
    \ifx\svgscale\undefined%
      \relax%
    \else%
      \setlength{\unitlength}{\unitlength * \real{\svgscale}}%
    \fi%
  \else%
    \setlength{\unitlength}{\svgwidth}%
  \fi%
  \global\let\svgwidth\undefined%
  \global\let\svgscale\undefined%
  \makeatother%
  \begin{picture}(1,1.8084928)%
    \lineheight{1}%
    \setlength\tabcolsep{0pt}%
    \put(0,0){\includegraphics[width=\unitlength,page=1]{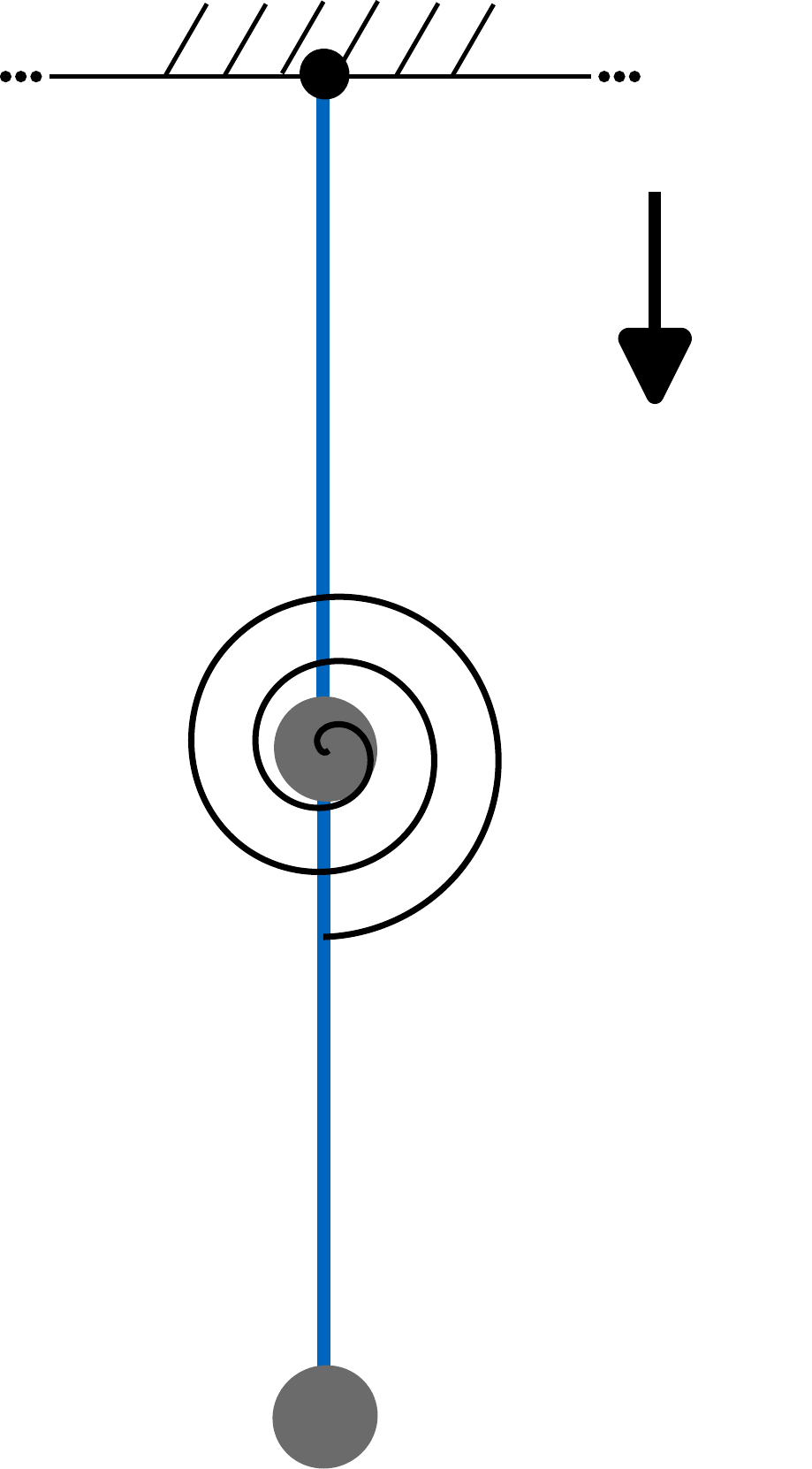}}%
    \put(0.66749544,0.84411032){\color[rgb]{0,0,0}\makebox(0,0)[lt]{\lineheight{1.25}\smash{\begin{tabular}[t]{l}$k$\end{tabular}}}}%
    \put(0.85171943,1.45437031){\color[rgb]{0,0,0}\makebox(0,0)[lt]{\lineheight{1.25}\smash{\begin{tabular}[t]{l}$g$\end{tabular}}}}%
  \end{picture}%
\endgroup%

          \subcaption{$V_{s1}$}\label{fig:Vs1}
        \end{subfigure}
        &
        \begin{subfigure}[b]{.3\columnwidth}
          \def\svgwidth{6em}
\begingroup%
  \makeatletter%
  \providecommand\color[2][]{%
    \errmessage{(Inkscape) Color is used for the text in Inkscape, but the package 'color.sty' is not loaded}%
    \renewcommand\color[2][]{}%
  }%
  \providecommand\transparent[1]{%
    \errmessage{(Inkscape) Transparency is used (non-zero) for the text in Inkscape, but the package 'transparent.sty' is not loaded}%
    \renewcommand\transparent[1]{}%
  }%
  \providecommand\rotatebox[2]{#2}%
  \newcommand*\fsize{\dimexpr\f@size pt\relax}%
  \newcommand*\lineheight[1]{\fontsize{\fsize}{#1\fsize}\selectfont}%
  \ifx\svgwidth\undefined%
    \setlength{\unitlength}{566.32637609bp}%
    \ifx\svgscale\undefined%
      \relax%
    \else%
      \setlength{\unitlength}{\unitlength * \real{\svgscale}}%
    \fi%
  \else%
    \setlength{\unitlength}{\svgwidth}%
  \fi%
  \global\let\svgwidth\undefined%
  \global\let\svgscale\undefined%
  \makeatother%
  \begin{picture}(1,0.93421354)%
    \lineheight{1}%
    \setlength\tabcolsep{0pt}%
    \put(0,0){\includegraphics[width=\unitlength,page=1]{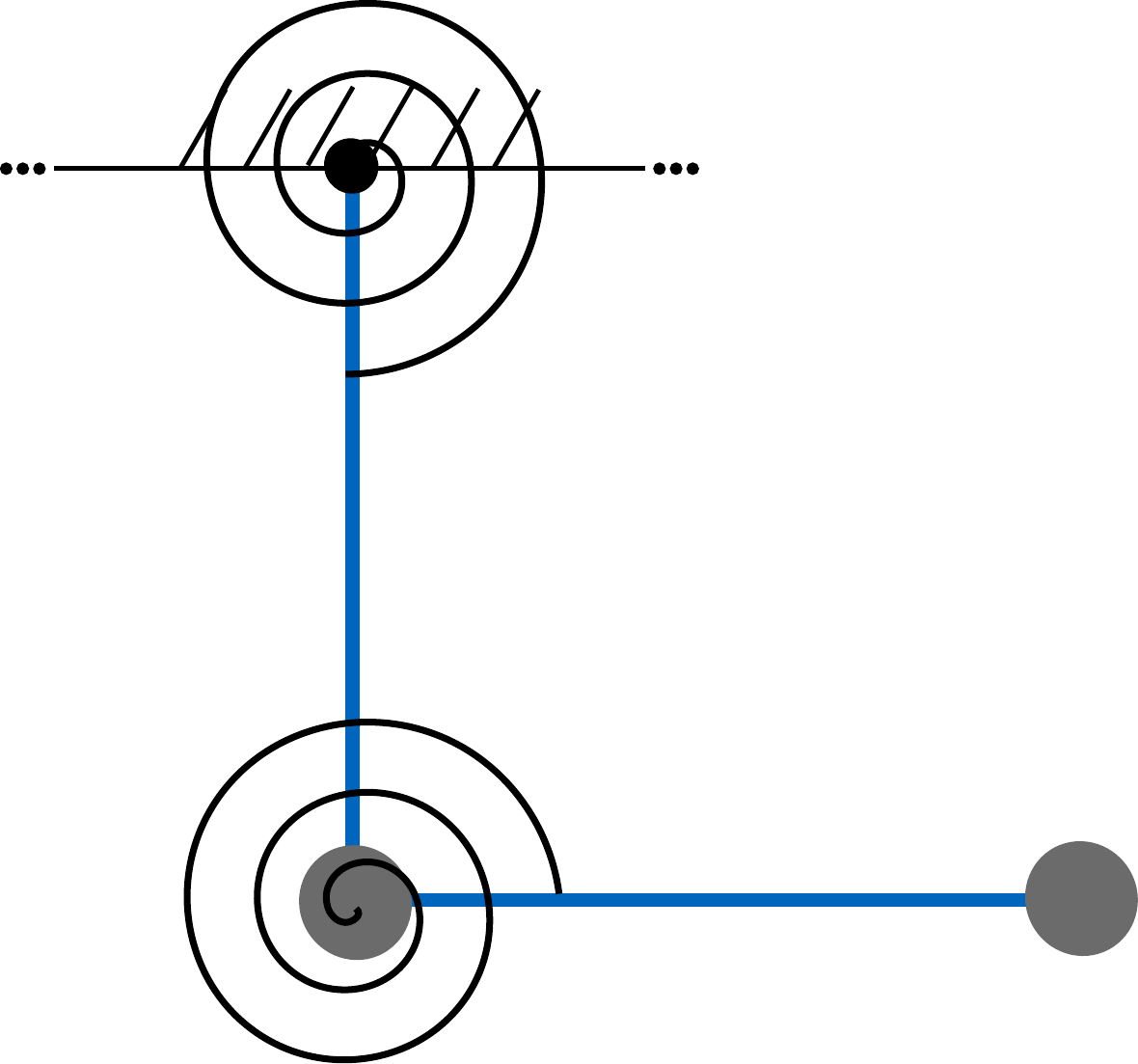}}%
    \put(0.46442547,0.26793141){\color[rgb]{0,0,0}\makebox(0,0)[lt]{\lineheight{1.25}\smash{\begin{tabular}[t]{l}$k$\end{tabular}}}}%
    \put(0.46520532,0.62438315){\color[rgb]{0,0,0}\makebox(0,0)[lt]{\lineheight{1.25}\smash{\begin{tabular}[t]{l}$k$\end{tabular}}}}%
  \end{picture}%
\endgroup%

          \vspace{1.55cm}
          \subcaption{$V_{s2}$}\label{fig:Vs2}
        \end{subfigure}
        &
        \begin{subfigure}[b]{.3\columnwidth}
          \def\svgwidth{6em}
\begingroup%
  \makeatletter%
  \providecommand\color[2][]{%
    \errmessage{(Inkscape) Color is used for the text in Inkscape, but the package 'color.sty' is not loaded}%
    \renewcommand\color[2][]{}%
  }%
  \providecommand\transparent[1]{%
    \errmessage{(Inkscape) Transparency is used (non-zero) for the text in Inkscape, but the package 'transparent.sty' is not loaded}%
    \renewcommand\transparent[1]{}%
  }%
  \providecommand\rotatebox[2]{#2}%
  \newcommand*\fsize{\dimexpr\f@size pt\relax}%
  \newcommand*\lineheight[1]{\fontsize{\fsize}{#1\fsize}\selectfont}%
  \ifx\svgwidth\undefined%
    \setlength{\unitlength}{478.07437359bp}%
    \ifx\svgscale\undefined%
      \relax%
    \else%
      \setlength{\unitlength}{\unitlength * \real{\svgscale}}%
    \fi%
  \else%
    \setlength{\unitlength}{\svgwidth}%
  \fi%
  \global\let\svgwidth\undefined%
  \global\let\svgscale\undefined%
  \makeatother%
  \begin{picture}(1,1.11166689)%
    \lineheight{1}%
    \setlength\tabcolsep{0pt}%
    \put(0,0){\includegraphics[width=\unitlength,page=1]{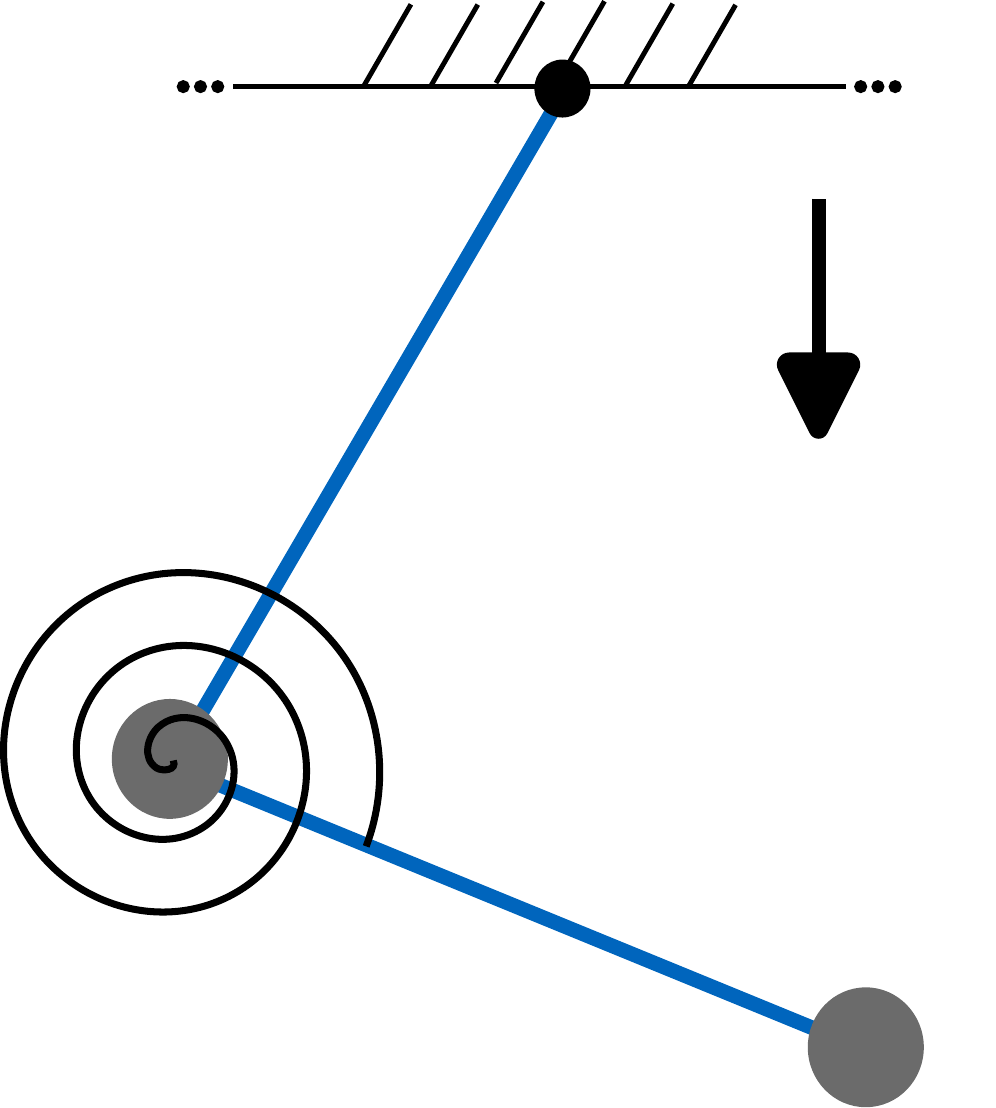}}%
    \put(0.86322285,0.80339654){\color[rgb]{0,0,0}\makebox(0,0)[lt]{\lineheight{1.25}\smash{\begin{tabular}[t]{l}$g$\end{tabular}}}}%
    \put(0.42119468,0.37053664){\color[rgb]{0,0,0}\makebox(0,0)[lt]{\lineheight{1.25}\smash{\begin{tabular}[t]{l}$k$\end{tabular}}}}%
  \end{picture}%
\endgroup%

          \vspace{1.3cm}
          \subcaption{$V_{a}$}\label{fig:Va}
        \end{subfigure}
      \end{tabular}
    \caption{Equilibrium configurations of the three tested potentials: $V_{s1}$ with nonlinear terms from gravity and $V_{s2}$ with only quadratic spring terms are symmetric under $(S_2,\tau_2)$, in the appropriate coordinates. $V_a$ is obtained from $V_{s1}$ by shifting the equilibrium of the spring to $\nicefrac{\pi}{2}$ such that $V_a(\bar{q} + q) \neq V_a(\bar{q} - q)$ for the minimum $\bar{q}$ of $V_a$, making $V_a$ not symmetric.}
    \label{fig:potentials}
\end{figure}

Consider the following potential functions, see also Fig.~\ref{fig:potentials}: 
\begin{align}
    V_{s1}(q) &= \nicefrac{1}{2}k(q_2)^2 + V_g(q), \nonumber\\
    V_{s2}(q) &= \nicefrac{1}{2}k (q_1)^2 + \nicefrac{1}{2}k(q_2-\nicefrac{\pi}{2})^2  \,, \nonumber \\
    V_{a}(q) &= \nicefrac{1}{2}k(q_2-\nicefrac{\pi}{2})^2 + V_g(q)\,. \nonumber
\end{align}
where $V_g(q) = - dmg \big(2\cos(q_1)+\cos(q_1+q_2)\big)$ is the gravitational potential corresponding to the double pendulum and all other terms are quadratic potentials corresponding to linear springs. The essential differences are that $V_{s1}$ and $V_{s2}$ fulfill condition (i) of Theorem~\ref{Theorem:1}, i.e., $V_{s1}(\qeq_{s1} + q) = V_{s1}(\qeq_{s1}-q)$ with $\qeq_{s1} = (0,0)$ and $V_{s2}(\qeq_{s2} + q) = V_{s2}(\qeq_{s2} - q)$ with $\qeq_{s2} = (0,\nicefrac{\pi}{2})$, while $V_{a}(q)$ does not verify this condition in the given coordinates\footnote{$V_a$ may locally fulfill condition (i) of Theorem~\ref{Theorem:1} in other coordinates, this is discussed in Sec.~\ref{sssec:Examples-2DoF-Masses}.}. Further, $V_{s2}$ results in a linear force and $V_{s1},\, V_a$ do not.

Parameters are chosen as $m = 0.4\text{kg}, I = \nicefrac{1}{12} \text{kgm}^2, d = 1\text{m}, k = 10\frac{\text{Nm}}{\text{rad}}, g = 9.81\text{m}/\text{s}^2$\footnote{The units are correct for the double pendulum. For the coupled masses, $k$ should have units of $\frac{\text{N}}{\text{m}}$ and $V_g$ should be interpreted as a nonlinear potential $V_g(q) = -k_2 \big(2\cos(\frac{q_1}{d})+\cos(\frac{q_1+q_2}{d})\big)$ with $k_2 = 3.924\text{N}$}. \\ 
\subsubsection{Procedure}
By different pairings of the above metric tensors and potential functions, one obtains a variety of systems with dynamics governed by~\eqref{eq:hamiltonian_ode_q},~\eqref{eq:hamiltonian_ode_p}, subsequently denoted by a tuple $(M, V)$. 

For these systems, the conditions for Theorem~\ref{thm:Eigenmanifold_Condition}, Corollarly~\ref{cr:n-Eigenmanifolds} and Theorem~\ref{Theorem:1} are checked. 

To confirm predictions, the LSMs of these systems are computed using numerical continuation strategies~\cite{Sachtler2024}. The results are summarized in Figs.~\ref{fig:allplots},~\ref{s4f6} and in Tab.~\ref{s4tab}, and explained in the following.

\begin{table*}[t]
\caption{\small Summary of tested $2\--$DoF examples \label{s4tab}}
\centering
 \begin{tabular}{c|c|c|c|c|c}
 \textbf{System} & \textbf{System type} & \textbf{Potential} & \textbf{Inertia} & \textbf{Thm. \ref{Theorem:1} satisfied?} & \textbf{Figure} \\\hline
 \vspace{-0.1cm}
 $(M_{DP}, V_{s1})$ & non-Euclidean & Symmetric & Symmetric  & \cmark & Fig.~\ref{s4f1}, \ref{s4f1_2} \\ 
  \vspace{-0.1cm}
 $(M_{DP}, V_{a})$  & non-Euclidean & Asymmetric & Symmetric & \xmark & Fig.~\ref{s4f2}, \ref{s4f2_2} \\
 \vspace{-0.1cm}
 $(M_{DP}, V_{s2})$ & non-Euclidean & Symmetric & Asymmetric & \xmark & Fig.~\ref{s4f3}, \ref{s4f3_2} \\
 \vspace{-0.1cm}
 $(M_{C}, V_{s1})$  & Euclidean     & Symmetric & Symmetric & \cmark & Fig.~\ref{s4f4}, \ref{s4f4_2} \\
 $(M_{C}, V_{a})$   & Euclidean     & Asymmetric & Symmetric & \xmark & Fig.~\ref{s4f5}, \ref{s4f5_2} \\
\end{tabular}
\end{table*}

\subsubsection{Euclidean system: Coupled masses}\label{sssec:Examples-2DoF-Masses}
Here, we focus on the systems $(M_C, V_{s1})$ and $(M_C, V_{a})$. The system $(M_C, V_{s2})$ is linear, and is omitted. The conditions for Theorem~\ref{thm:Eigenmanifold_Condition} are fulfilled for both systems: for details on conditions (i) to (iii), refer to the supplementary code~\cite{data-Link}. 

\paragraph{Rosenberg manifolds}
For $(M_C, V_{s1})$, the equilibrium is $\qeq_{s1} = (0,0)$ and the eigenvalues of $M_C(\qeq_{s1})^{-1}\text{Hess}(V_{s1})(\qeq_{s1})$ are $\lambda_1^2 = 21.95\,,\;\lambda_2^2 = 42.29$, with the ratio $\frac{\lambda_2}{\lambda_1} = 1.39$ non-integer, so the non-resonance condition of Theorem~\ref{thm:Eigenmanifold_Condition} holds. By Corollary~\ref{cr:n-Eigenmanifolds} there are $2$ Eigenmanifolds, since also the inverse $\frac{\lambda_1}{\lambda_2}$ is non-integer. Also the conditions for Theorem~\ref{Theorem:1} hold: $V_{s1}(q) = V_{s1}(-q)$, $M_C(q) = M_C(-q)$ holds trivially, and $\qeq_{s1}$ is the unique minimum of $V_{s1}(q)$. Hence, the expected Eigenmanifolds are Rosenberg manifolds. 

This is confirmed in Figs.~\ref{s4f4},~\ref{s4f4_2}: the oscillations (in blue) are projected to configuration space and show as curves without self-intersection and with open ends, i.e., the embedding property (2) of Theorem~\ref{thm:Eigenmanifold_Condition} holds. The shown eigenmodes also pass through the equilibrium configuration $\qeq_{s1}$, making them Rosenberg modes (cf. Def.~\ref{def:rosenberg_mode}).  Fig.~\ref{s4f6} shows this for the collection of modes: the minimum potential energy $V$ along a given mode together is plotted against the energy of the starting point of that mode. This minimum potential is zero when the corresponding mode passes through the equilibrium, and shows a flat line at zero, indicating that the collection of modes is a Rosenberg manifold (cf. Def.~\ref{def:rosenberg_geo}). Videos in the supplementary data~\cite{data-Link} visualize this for the full family of oscillation. 

\paragraph{Eigenmanifolds}
For $(M_C, V_{a})$, equivalent computations predict $2$ Eigenmanifolds starting at the equilibrium $\qeq_{a} = (-0.39,1.28)$. However, the conditions for Theorem~\ref{Theorem:1} do not hold in the given coordinates, such that Rosenberg manifolds are not expected. 

This is confirmed in Figs.~\ref{s4f5},~\ref{s4f5_2}: the results show eigenmodes that do not pass through $\qeq_a$, highlighting that Eigenmanifolds (falling outside Rosenberg's definition) also exist on Euclidean configuration manifolds. Videos in the supplementary data~\cite{data-Link} show that although most modes do not pass through $\qeq_a$, isolated modes do pass through $\qeq_a$. These isolated modes are also recognized in Fig.~\ref{s4f6}: the plotted minimum potential is zero when the corresponding mode passes through the equilibrium, such that the observed isolated zeros identify isolated geometric Rosenberg modes within Eigenmanifolds. This detail does not contradict the presented theory, since isolated Rosenberg modes do not make a Rosenberg manifold. The example shows that non-Rosenberg Eigenmanifolds exist also in Euclidean systems. It also shows, in hindsight, that no other coordinate system exists in which $V_a$ and $M_C$ fulfill condition $(i)$ of Theorem~\ref{Theorem:1}\footnote{This is explained as follows: a coordinate change resulting in a symmetric $V_a$ would also transform $M_C$ into a non-constant inertia tensor, in particular one for which condition (i) does not hold.}.


\subsubsection{Non-Euclidean system: Double pendulum}

%
We focus on the systems $(M_{DP},V_{s1})$, $(M_{DP},V_{s2})$ and $(M_{DP},V_{a})$. The conditions of Theorem~\ref{thm:Eigenmanifold_Condition} are fulfilled by all three systems. Eigenvalues of all systems 
are non-resonant, such that Corollary~\ref{cr:n-Eigenmanifolds} guarantees $2$ Eigenmanifolds. Computations in the supplementary code~\cite{data-Link} show this in detail. Particular systems highlight aspects of the coordinate-free conditions of Theorem~\ref{Theorem:1}, and the ensemble confirms the presence of Eigenmanifolds and Rosenberg manifolds on non-Euclidean configuration manifolds up to high energies.

\begin{figure*}[t]
\centering
    \begin{subfigure}{.495\textwidth}
    \centering
      \def\svgwidth{30em}\scriptsize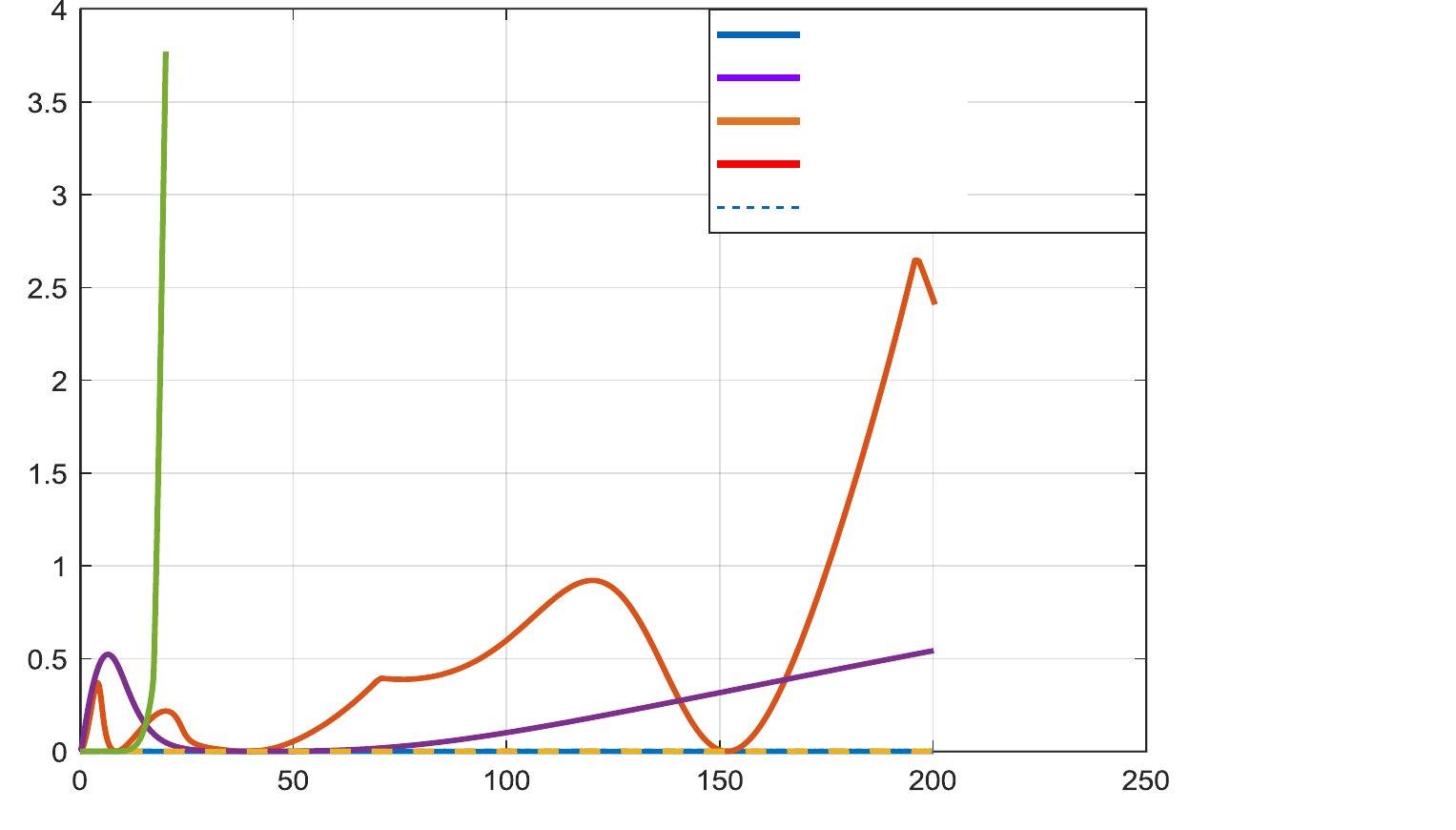
      \subcaption{Generators $R_1$}
      \label{s4f6a}
    \end{subfigure}
    \begin{subfigure}{.495\textwidth}
        \centering
      \def\svgwidth{30em}\scriptsize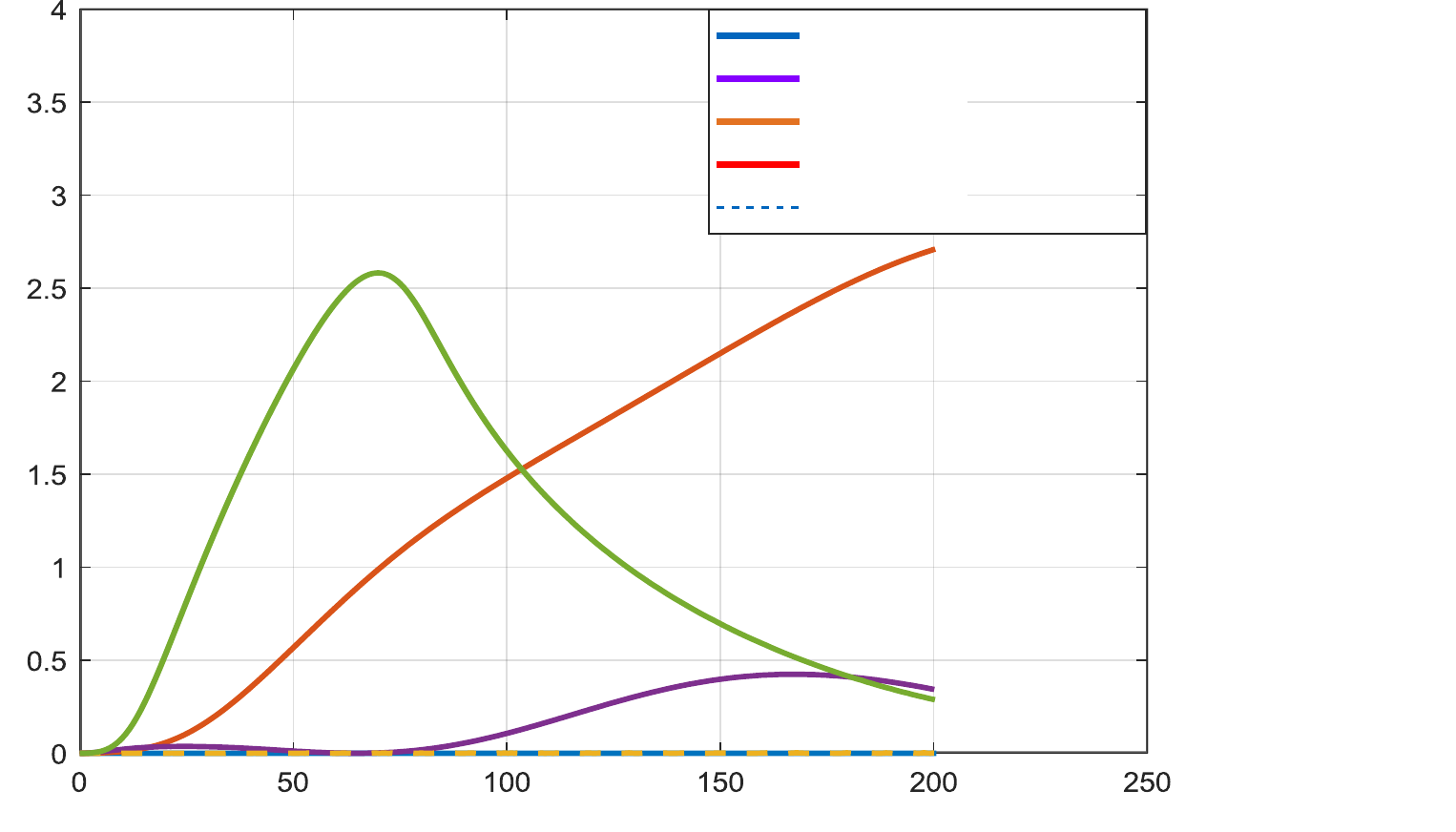
      \subcaption{Generators $R_2$}
      \label{s4f6b}
    \end{subfigure}
\caption{\small Minimal potential energy for modes of the example systems, split by modes collected on first and second generator. }\label{s4f6}
\end{figure*}

\paragraph{Rosenberg manifolds}
The system $(M_{DP},V_{s1})$ satisfies the conditions of Theorem~\ref{Theorem:1}: $V_{s1}(q) = V_{s1}(-q)$, $M_{DP}(q) = M_{DP}(-q)$, and the minimum of $V_{s1}(q)$ is unique. Example trajectories and generators of the predicted Rosenberg manifolds are shown in Figs.~\ref{s4f1} and~\ref{s4f1_2}, where the oscillation is seen to pass through the equilibrium. Videos in the supplementary data~\cite{data-Link} show this for the full family, but also Fig.~\ref{s4f6} shows that the minimum potential along all modes in the family is zero, such that all modes on the family pass through the equilibrium configuration.

\paragraph{Eigenmanifolds}
The systems $(M_{DP},V_{s2})$ and $(M_{DP},V_{a})$ fail to satisfy Theorem~\ref{Theorem:1} in different ways.
%

For the system $(M_{DP},V_{s2})$, the potential $V_{s2}(\qeq_{s2}+q') = V_{s2}(\qeq_{s2}-q')$ fulfills the symmetry condition in coordinates $q' = q - \qeq_{s2}$, and it has a unique minimum. However, in these new coordinates $M_{DP}(q') \neq M_{DP}(-q')$, and Theorem~\ref{Theorem:1} is not fuliflled. This is equivalent to the example $(M_C,V_a)$, in which no coordinates were found such that $V$ and $M$ are both symmetric. The families of oscillations in Figs.~\ref{s4f3} and~\ref{s4f3_2} are Eigenmanifolds, but not Rosenberg manifolds, confirming the insight. This is confirmed by Fig.~\ref{s4f6}, and videos in the supplementary data~\cite{data-Link}, which only show isolated Rosenberg modes.

For the system $(M_{DP}, V_{a})$, the potential does not fulfill the symmetry condition. Analogous to previous observations, Figs.~\ref{s4f2} \&~\ref{s4f2_2} confirm the presence of Eigenmanifolds, but not Rosenberg manifolds. 
%
%
%
%
\begin{figure*}[t]
\begin{tabular}{ccc}
    \begin{subfigure}{.32\textwidth}
        \captionsetup{
        justification=raggedright,
        singlelinecheck=false
        }
      \def\svgwidth{12em}\scriptsize
\begingroup%
  \makeatletter%
  \providecommand\color[2][]{%
    \errmessage{(Inkscape) Color is used for the text in Inkscape, but the package 'color.sty' is not loaded}%
    \renewcommand\color[2][]{}%
  }%
  \providecommand\transparent[1]{%
    \errmessage{(Inkscape) Transparency is used (non-zero) for the text in Inkscape, but the package 'transparent.sty' is not loaded}%
    \renewcommand\transparent[1]{}%
  }%
  \providecommand\rotatebox[2]{#2}%
  \newcommand*\fsize{\dimexpr\f@size pt\relax}%
  \newcommand*\lineheight[1]{\fontsize{\fsize}{#1\fsize}\selectfont}%
  \ifx\svgwidth\undefined%
    \setlength{\unitlength}{495.38807991bp}%
    \ifx\svgscale\undefined%
      \relax%
    \else%
      \setlength{\unitlength}{\unitlength * \real{\svgscale}}%
    \fi%
  \else%
    \setlength{\unitlength}{\svgwidth}%
  \fi%
  \global\let\svgwidth\undefined%
  \global\let\svgscale\undefined%
  \makeatother%
  \begin{picture}(1,0.79467866)%
    \lineheight{1}%
    \setlength\tabcolsep{0pt}%
    \put(0.02507977,0.34921501){\color[rgb]{0,0,0}\rotatebox{90}{\makebox(0,0)[lt]{\lineheight{1.25}\smash{\begin{tabular}[t]{l}$q_2$ in rad\end{tabular}}}}}%
    \put(0,0){\includegraphics[width=\unitlength,page=1]{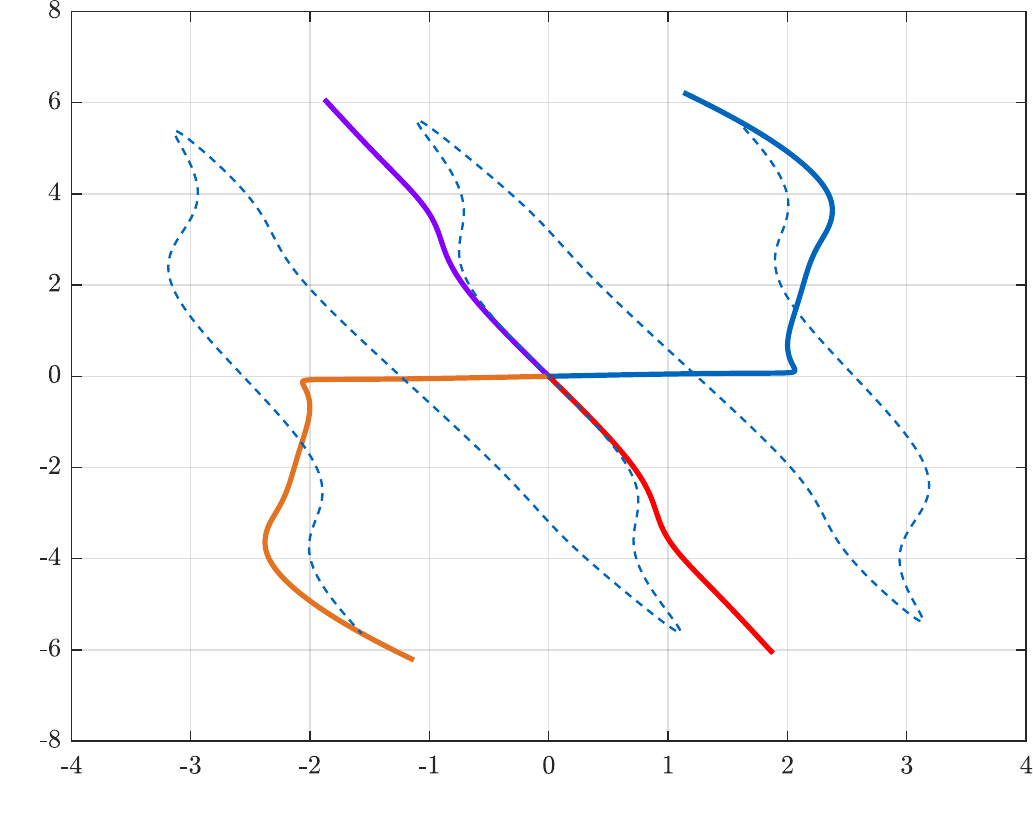}}%
    \put(0.44219468,0.00778504){\color[rgb]{0,0,0}\makebox(0,0)[lt]{\lineheight{1.25}\smash{\begin{tabular}[t]{l}$q_1$ in rad\end{tabular}}}}%
  \end{picture}%
\endgroup%

      \subcaption{$(M_{DP}, V_{s1}), R_1$}
      \label{s4f1}
    \end{subfigure}&
    \hspace{-0.3cm}
    \begin{subfigure}{.32\textwidth}
        \captionsetup{
        justification=raggedright,
        singlelinecheck=false
        }
      \def\svgwidth{13em}\scriptsize
\begingroup%
  \makeatletter%
  \providecommand\color[2][]{%
    \errmessage{(Inkscape) Color is used for the text in Inkscape, but the package 'color.sty' is not loaded}%
    \renewcommand\color[2][]{}%
  }%
  \providecommand\transparent[1]{%
    \errmessage{(Inkscape) Transparency is used (non-zero) for the text in Inkscape, but the package 'transparent.sty' is not loaded}%
    \renewcommand\transparent[1]{}%
  }%
  \providecommand\rotatebox[2]{#2}%
  \newcommand*\fsize{\dimexpr\f@size pt\relax}%
  \newcommand*\lineheight[1]{\fontsize{\fsize}{#1\fsize}\selectfont}%
  \ifx\svgwidth\undefined%
    \setlength{\unitlength}{606.0796253bp}%
    \ifx\svgscale\undefined%
      \relax%
    \else%
      \setlength{\unitlength}{\unitlength * \real{\svgscale}}%
    \fi%
  \else%
    \setlength{\unitlength}{\svgwidth}%
  \fi%
  \global\let\svgwidth\undefined%
  \global\let\svgscale\undefined%
  \makeatother%
  \begin{picture}(1,0.71789843)%
    \lineheight{1}%
    \setlength\tabcolsep{0pt}%
    \put(0,0){\includegraphics[width=\unitlength,page=1]{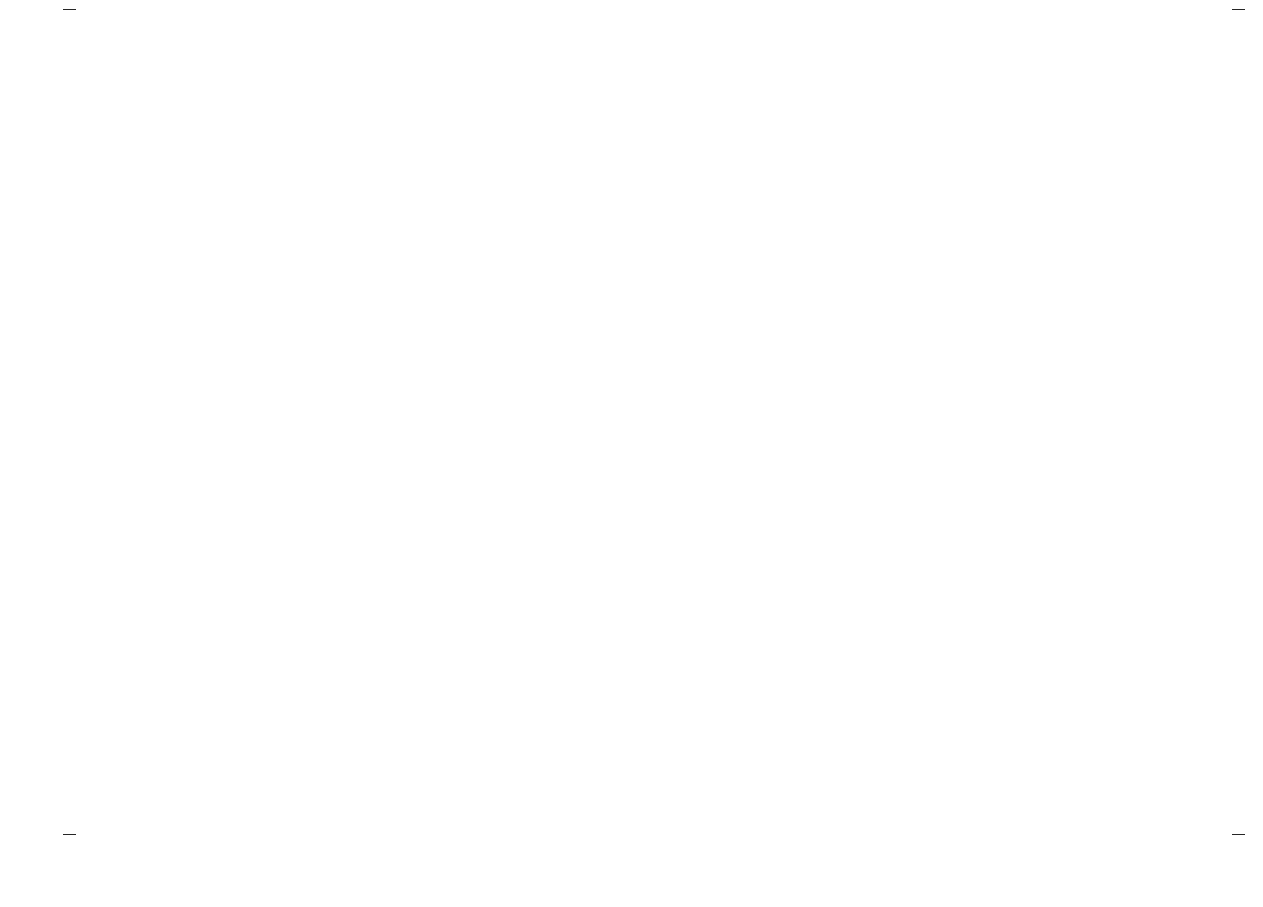}}%
    \put(0.02049932,0.30761436){\color[rgb]{0,0,0}\rotatebox{90}{\makebox(0,0)[lt]{\lineheight{0}\smash{\begin{tabular}[t]{l}$q_2$ in rad\end{tabular}}}}}%
    \put(0.44340464,0.00636322){\color[rgb]{0,0,0}\makebox(0,0)[lt]{\lineheight{0}\smash{\begin{tabular}[t]{l}$q_1$ in rad\end{tabular}}}}%
    \put(0,0){\includegraphics[width=\unitlength,page=2]{DP_Sym_2_svg-tex.pdf}}%
  \end{picture}%
\endgroup%

      \subcaption{$(M_{DP}, V_{s1}), R_2$}
      \label{s4f1_2}
    \end{subfigure}&
    \begin{subfigure}{.33\textwidth}
        \captionsetup{
        justification=raggedright,
        singlelinecheck=false
        }
      \def\svgwidth{12em}\scriptsize
\begingroup%
  \makeatletter%
  \providecommand\color[2][]{%
    \errmessage{(Inkscape) Color is used for the text in Inkscape, but the package 'color.sty' is not loaded}%
    \renewcommand\color[2][]{}%
  }%
  \providecommand\transparent[1]{%
    \errmessage{(Inkscape) Transparency is used (non-zero) for the text in Inkscape, but the package 'transparent.sty' is not loaded}%
    \renewcommand\transparent[1]{}%
  }%
  \providecommand\rotatebox[2]{#2}%
  \newcommand*\fsize{\dimexpr\f@size pt\relax}%
  \newcommand*\lineheight[1]{\fontsize{\fsize}{#1\fsize}\selectfont}%
  \ifx\svgwidth\undefined%
    \setlength{\unitlength}{285.60945753bp}%
    \ifx\svgscale\undefined%
      \relax%
    \else%
      \setlength{\unitlength}{\unitlength * \real{\svgscale}}%
    \fi%
  \else%
    \setlength{\unitlength}{\svgwidth}%
  \fi%
  \global\let\svgwidth\undefined%
  \global\let\svgscale\undefined%
  \makeatother%
  \begin{picture}(1,0.5136665)%
    \lineheight{1}%
    \setlength\tabcolsep{0pt}%
    \put(0,0){\includegraphics[width=\unitlength,page=1]{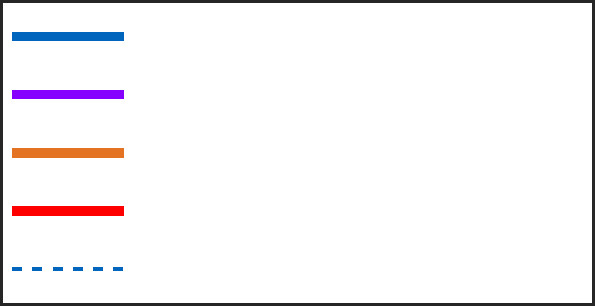}}%
    \put(0.2311515,0.43071758){\color[rgb]{0,0,0}\makebox(0,0)[lt]{\lineheight{1.25}\smash{\begin{tabular}[t]{l}Generator $R_{11}$\end{tabular}}}}%
    \put(0.23140533,0.33796545){\color[rgb]{0,0,0}\makebox(0,0)[lt]{\lineheight{1.25}\smash{\begin{tabular}[t]{l}Generator $R_{21}$\end{tabular}}}}%
    \put(0.23154508,0.2447169){\color[rgb]{0,0,0}\makebox(0,0)[lt]{\lineheight{1.25}\smash{\begin{tabular}[t]{l}Generator $R_{12}$\end{tabular}}}}%
    \put(0.23279003,0.14499427){\color[rgb]{0,0,0}\makebox(0,0)[lt]{\lineheight{1.25}\smash{\begin{tabular}[t]{l}Generator $R_{22}$\end{tabular}}}}%
    \put(0.23748092,0.04652552){\color[rgb]{0,0,0}\makebox(0,0)[lt]{\lineheight{1.25}\smash{\begin{tabular}[t]{l}Example Trajectory\end{tabular}}}}%
  \end{picture}%
\endgroup%

      \vspace{1cm}
      \caption{Legend}
    \end{subfigure}\\
    \vspace{.33cm}&&\\
    \end{tabular}
    \begin{tabular}{cccc}
    \begin{subfigure}{.23\textwidth}
      \def\svgwidth{11em}\scriptsize
\begingroup%
  \makeatletter%
  \providecommand\color[2][]{%
    \errmessage{(Inkscape) Color is used for the text in Inkscape, but the package 'color.sty' is not loaded}%
    \renewcommand\color[2][]{}%
  }%
  \providecommand\transparent[1]{%
    \errmessage{(Inkscape) Transparency is used (non-zero) for the text in Inkscape, but the package 'transparent.sty' is not loaded}%
    \renewcommand\transparent[1]{}%
  }%
  \providecommand\rotatebox[2]{#2}%
  \newcommand*\fsize{\dimexpr\f@size pt\relax}%
  \newcommand*\lineheight[1]{\fontsize{\fsize}{#1\fsize}\selectfont}%
  \ifx\svgwidth\undefined%
    \setlength{\unitlength}{579.12866241bp}%
    \ifx\svgscale\undefined%
      \relax%
    \else%
      \setlength{\unitlength}{\unitlength * \real{\svgscale}}%
    \fi%
  \else%
    \setlength{\unitlength}{\svgwidth}%
  \fi%
  \global\let\svgwidth\undefined%
  \global\let\svgscale\undefined%
  \makeatother%
  \begin{picture}(1,0.72516944)%
    \lineheight{1}%
    \setlength\tabcolsep{0pt}%
    \put(0,0){\includegraphics[width=\unitlength,page=1]{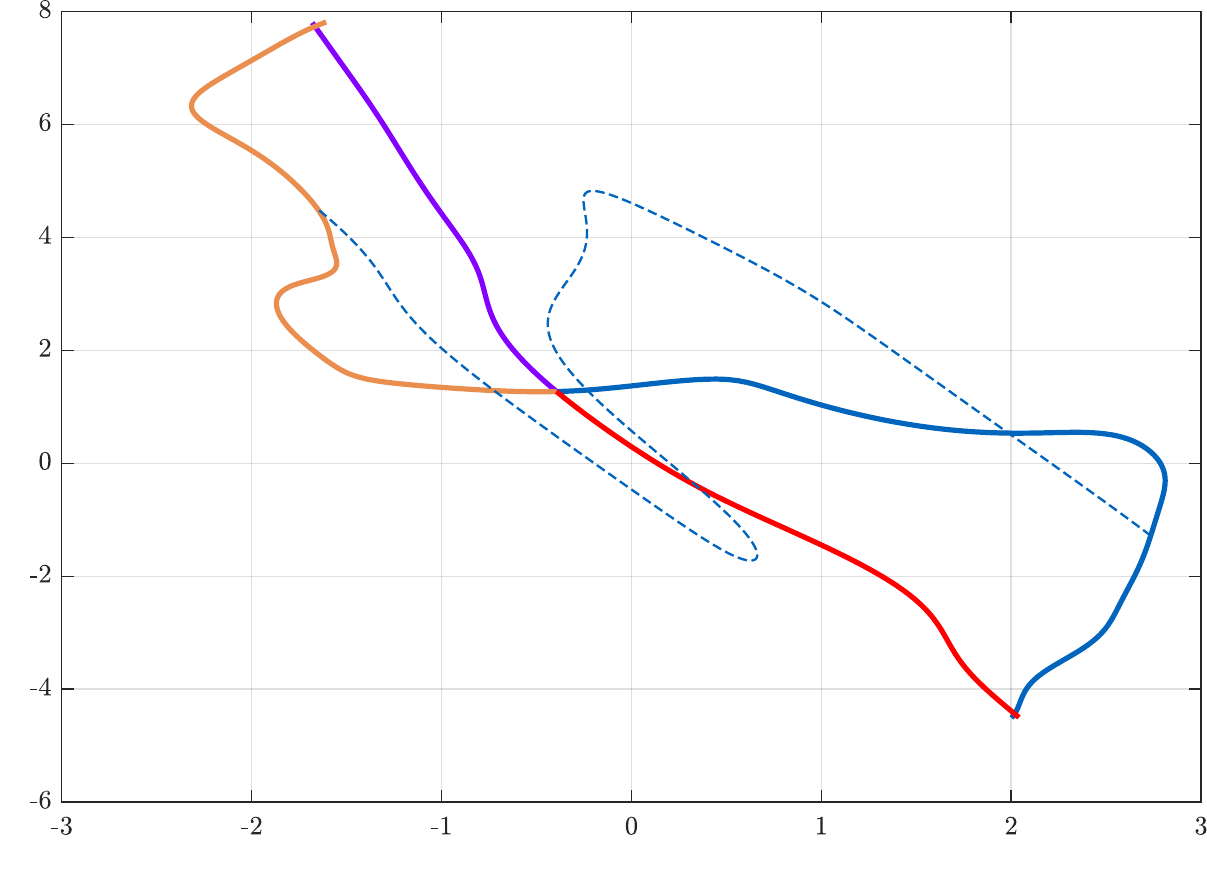}}%
    \put(0.0214533,0.2670085){\color[rgb]{0,0,0}\rotatebox{90}{\makebox(0,0)[lt]{\lineheight{1.25}\smash{\begin{tabular}[t]{l}$q_2$ in rad\end{tabular}}}}}%
    \put(0.45539793,0.00665935){\color[rgb]{0,0,0}\makebox(0,0)[lt]{\lineheight{1.25}\smash{\begin{tabular}[t]{l}$q_1$ in rad\end{tabular}}}}%
  \end{picture}%
\endgroup%

      \subcaption{$(M_{DP}, V_{a}), R_1$}
      \label{s4f2}
    \end{subfigure}&
    \hspace{0.3cm}
    \begin{subfigure}{.23\textwidth}
      \def\svgwidth{11em}\scriptsize
\begingroup%
  \makeatletter%
  \providecommand\color[2][]{%
    \errmessage{(Inkscape) Color is used for the text in Inkscape, but the package 'color.sty' is not loaded}%
    \renewcommand\color[2][]{}%
  }%
  \providecommand\transparent[1]{%
    \errmessage{(Inkscape) Transparency is used (non-zero) for the text in Inkscape, but the package 'transparent.sty' is not loaded}%
    \renewcommand\transparent[1]{}%
  }%
  \providecommand\rotatebox[2]{#2}%
  \newcommand*\fsize{\dimexpr\f@size pt\relax}%
  \newcommand*\lineheight[1]{\fontsize{\fsize}{#1\fsize}\selectfont}%
  \ifx\svgwidth\undefined%
    \setlength{\unitlength}{600.12950592bp}%
    \ifx\svgscale\undefined%
      \relax%
    \else%
      \setlength{\unitlength}{\unitlength * \real{\svgscale}}%
    \fi%
  \else%
    \setlength{\unitlength}{\svgwidth}%
  \fi%
  \global\let\svgwidth\undefined%
  \global\let\svgscale\undefined%
  \makeatother%
  \begin{picture}(1,0.72620476)%
    \lineheight{1}%
    \setlength\tabcolsep{0pt}%
    \put(0,0){\includegraphics[width=\unitlength,page=1]{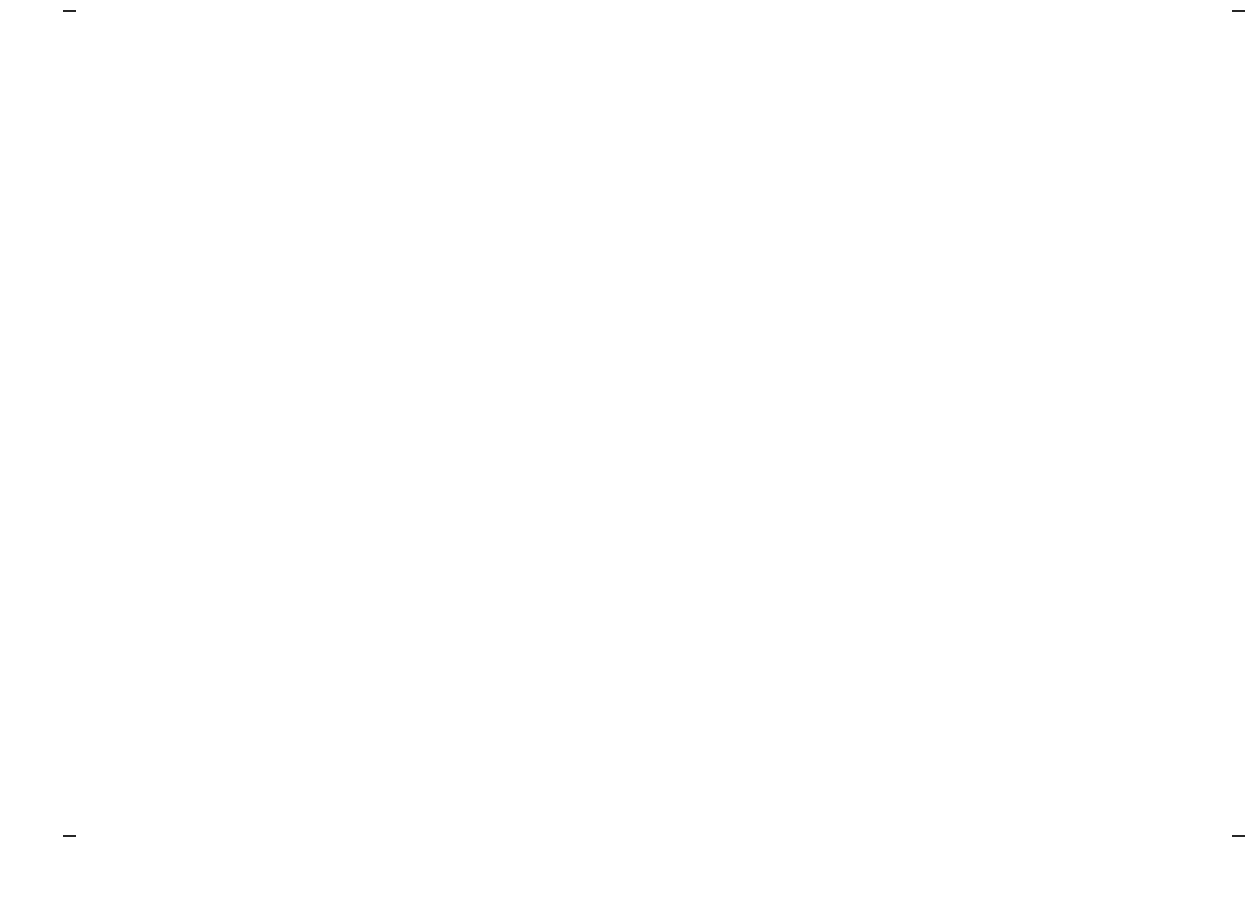}}%
    \put(0.02070256,0.31066427){\color[rgb]{0,0,0}\rotatebox{90}{\makebox(0,0)[lt]{\lineheight{0}\smash{\begin{tabular}[t]{l}$q_2$ in rad\end{tabular}}}}}%
    \put(0.44780087,0.00642631){\color[rgb]{0,0,0}\makebox(0,0)[lt]{\lineheight{0}\smash{\begin{tabular}[t]{l}$q_1$ in rad\end{tabular}}}}%
    \put(0,0){\includegraphics[width=\unitlength,page=2]{DP_As_2_svg-tex.pdf}}%
  \end{picture}%
\endgroup%

      \subcaption{$(M_{DP}, V_{a}), R_2$}
      \label{s4f2_2}
    \end{subfigure}&
    \hspace{0.3cm}
    \begin{subfigure}{.23\textwidth}
      \def\svgwidth{10em}\scriptsize
\begingroup%
  \makeatletter%
  \providecommand\color[2][]{%
    \errmessage{(Inkscape) Color is used for the text in Inkscape, but the package 'color.sty' is not loaded}%
    \renewcommand\color[2][]{}%
  }%
  \providecommand\transparent[1]{%
    \errmessage{(Inkscape) Transparency is used (non-zero) for the text in Inkscape, but the package 'transparent.sty' is not loaded}%
    \renewcommand\transparent[1]{}%
  }%
  \providecommand\rotatebox[2]{#2}%
  \newcommand*\fsize{\dimexpr\f@size pt\relax}%
  \newcommand*\lineheight[1]{\fontsize{\fsize}{#1\fsize}\selectfont}%
  \ifx\svgwidth\undefined%
    \setlength{\unitlength}{487.56343491bp}%
    \ifx\svgscale\undefined%
      \relax%
    \else%
      \setlength{\unitlength}{\unitlength * \real{\svgscale}}%
    \fi%
  \else%
    \setlength{\unitlength}{\svgwidth}%
  \fi%
  \global\let\svgwidth\undefined%
  \global\let\svgscale\undefined%
  \makeatother%
  \begin{picture}(1,0.81906467)%
    \lineheight{1}%
    \setlength\tabcolsep{0pt}%
    \put(0,0){\includegraphics[width=\unitlength,page=1]{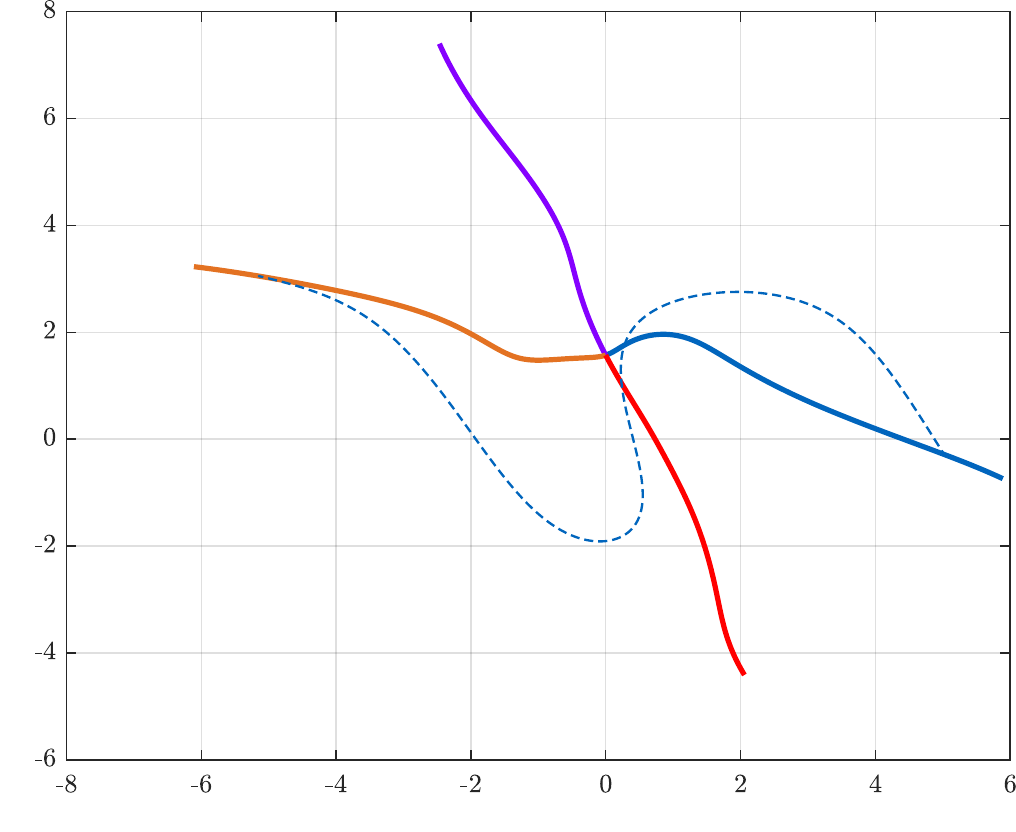}}%
    \put(0.45253318,0.00790998){\color[rgb]{0,0,0}\makebox(0,0)[lt]{\lineheight{1.25}\smash{\begin{tabular}[t]{l}$q_1$ in rad\end{tabular}}}}%
    \put(0.02548226,0.28521962){\color[rgb]{0,0,0}\rotatebox{90}{\makebox(0,0)[lt]{\lineheight{1.25}\smash{\begin{tabular}[t]{l}$q_2$ in rad\end{tabular}}}}}%
  \end{picture}%
\endgroup%

      \subcaption{$(M_{DP}, V_{s2}), R_1$}
      \label{s4f3}
    \end{subfigure}&
    \begin{subfigure}{.23\textwidth}
      \def\svgwidth{11em}\scriptsize
\begingroup%
  \makeatletter%
  \providecommand\color[2][]{%
    \errmessage{(Inkscape) Color is used for the text in Inkscape, but the package 'color.sty' is not loaded}%
    \renewcommand\color[2][]{}%
  }%
  \providecommand\transparent[1]{%
    \errmessage{(Inkscape) Transparency is used (non-zero) for the text in Inkscape, but the package 'transparent.sty' is not loaded}%
    \renewcommand\transparent[1]{}%
  }%
  \providecommand\rotatebox[2]{#2}%
  \newcommand*\fsize{\dimexpr\f@size pt\relax}%
  \newcommand*\lineheight[1]{\fontsize{\fsize}{#1\fsize}\selectfont}%
  \ifx\svgwidth\undefined%
    \setlength{\unitlength}{600.12950892bp}%
    \ifx\svgscale\undefined%
      \relax%
    \else%
      \setlength{\unitlength}{\unitlength * \real{\svgscale}}%
    \fi%
  \else%
    \setlength{\unitlength}{\svgwidth}%
  \fi%
  \global\let\svgwidth\undefined%
  \global\let\svgscale\undefined%
  \makeatother%
  \begin{picture}(1,0.72620475)%
    \lineheight{1}%
    \setlength\tabcolsep{0pt}%
    \put(0,0){\includegraphics[width=\unitlength,page=1]{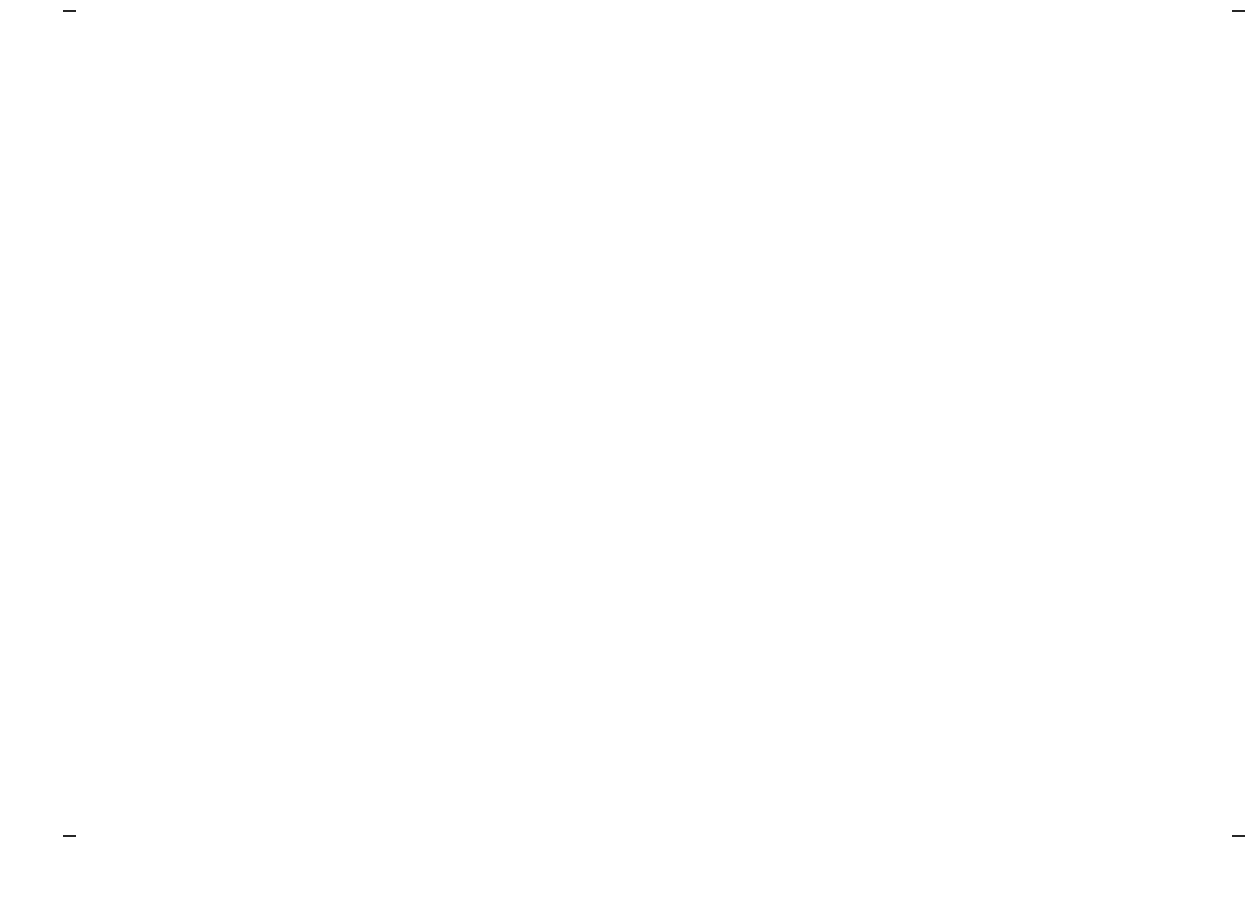}}%
    \put(0.02070256,0.31066427){\color[rgb]{0,0,0}\rotatebox{90}{\makebox(0,0)[lt]{\lineheight{0}\smash{\begin{tabular}[t]{l}$q_2$ in rad\end{tabular}}}}}%
    \put(0.44780087,0.00642631){\color[rgb]{0,0,0}\makebox(0,0)[lt]{\lineheight{0}\smash{\begin{tabular}[t]{l}$q_1$ in rad\end{tabular}}}}%
    \put(0,0){\includegraphics[width=\unitlength,page=2]{DP_AsM_2_svg-tex.pdf}}%
  \end{picture}%
\endgroup%

      \subcaption{$(M_{DP}, V_{s2}), R_2$}
      \label{s4f3_2}
    \end{subfigure}\\
    \vspace{.33cm}&&\\
    \begin{subfigure}{.23\textwidth}
      \def\svgwidth{10em}\scriptsize
\begingroup%
  \makeatletter%
  \providecommand\color[2][]{%
    \errmessage{(Inkscape) Color is used for the text in Inkscape, but the package 'color.sty' is not loaded}%
    \renewcommand\color[2][]{}%
  }%
  \providecommand\transparent[1]{%
    \errmessage{(Inkscape) Transparency is used (non-zero) for the text in Inkscape, but the package 'transparent.sty' is not loaded}%
    \renewcommand\transparent[1]{}%
  }%
  \providecommand\rotatebox[2]{#2}%
  \newcommand*\fsize{\dimexpr\f@size pt\relax}%
  \newcommand*\lineheight[1]{\fontsize{\fsize}{#1\fsize}\selectfont}%
  \ifx\svgwidth\undefined%
    \setlength{\unitlength}{600.18879741bp}%
    \ifx\svgscale\undefined%
      \relax%
    \else%
      \setlength{\unitlength}{\unitlength * \real{\svgscale}}%
    \fi%
  \else%
    \setlength{\unitlength}{\svgwidth}%
  \fi%
  \global\let\svgwidth\undefined%
  \global\let\svgscale\undefined%
  \makeatother%
  \begin{picture}(1,0.72620949)%
    \lineheight{1}%
    \setlength\tabcolsep{0pt}%
    \put(0,0){\includegraphics[width=\unitlength,page=1]{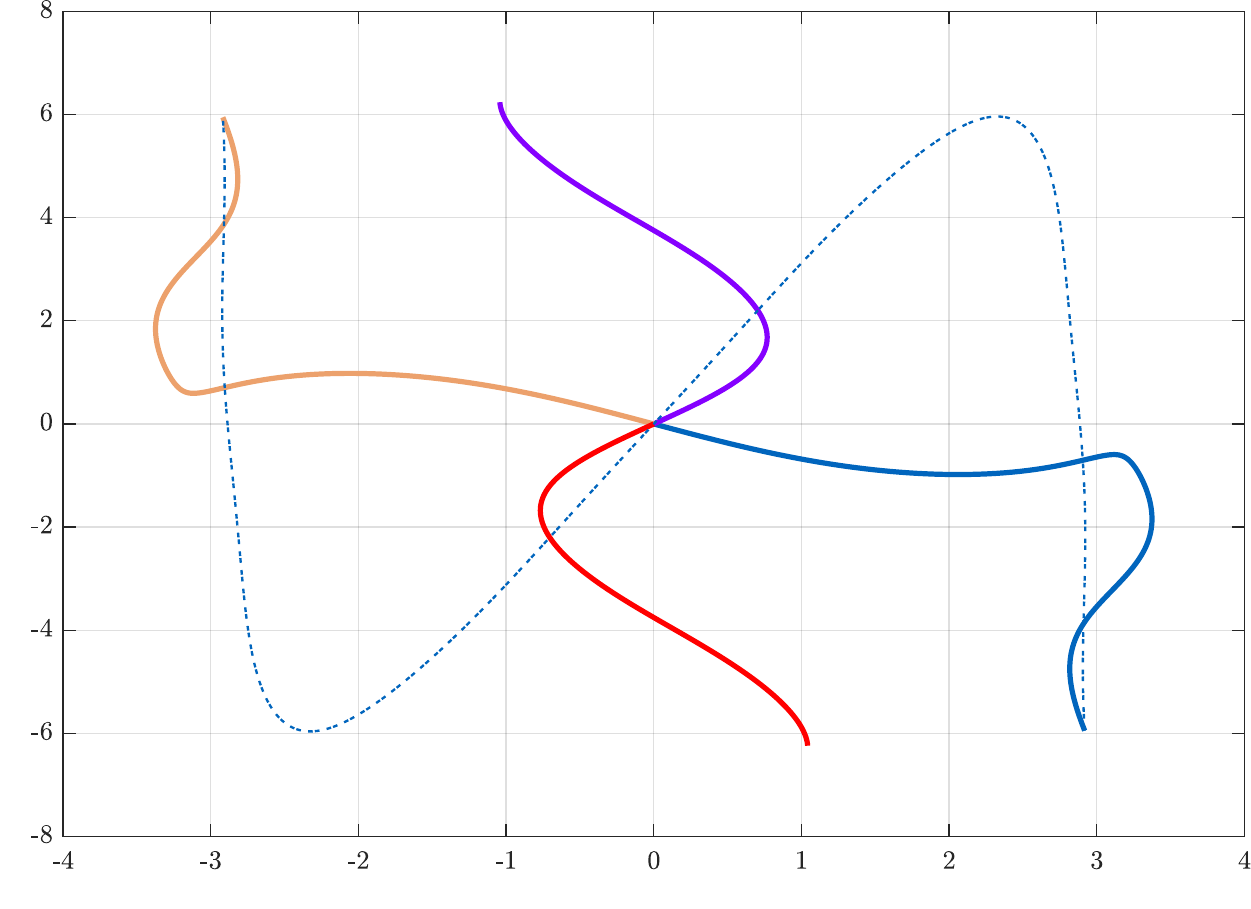}}%
    \put(0.02070052,0.31063358){\color[rgb]{0,0,0}\rotatebox{90}{\makebox(0,0)[lt]{\lineheight{1.25}\smash{\begin{tabular}[t]{l}$q_2$ in m\end{tabular}}}}}%
    \put(0.44775664,0.00642568){\color[rgb]{0,0,0}\makebox(0,0)[lt]{\lineheight{1.25}\smash{\begin{tabular}[t]{l}$q_1$ in m\end{tabular}}}}%
  \end{picture}%
\endgroup%

      \subcaption{$(M_{C}, V_{s1}), R_1$}
      \label{s4f4}
    \end{subfigure}&
    \begin{subfigure}{.23\textwidth}
      \def\svgwidth{10em}\scriptsize
\begingroup%
  \makeatletter%
  \providecommand\color[2][]{%
    \errmessage{(Inkscape) Color is used for the text in Inkscape, but the package 'color.sty' is not loaded}%
    \renewcommand\color[2][]{}%
  }%
  \providecommand\transparent[1]{%
    \errmessage{(Inkscape) Transparency is used (non-zero) for the text in Inkscape, but the package 'transparent.sty' is not loaded}%
    \renewcommand\transparent[1]{}%
  }%
  \providecommand\rotatebox[2]{#2}%
  \newcommand*\fsize{\dimexpr\f@size pt\relax}%
  \newcommand*\lineheight[1]{\fontsize{\fsize}{#1\fsize}\selectfont}%
  \ifx\svgwidth\undefined%
    \setlength{\unitlength}{600.93652398bp}%
    \ifx\svgscale\undefined%
      \relax%
    \else%
      \setlength{\unitlength}{\unitlength * \real{\svgscale}}%
    \fi%
  \else%
    \setlength{\unitlength}{\svgwidth}%
  \fi%
  \global\let\svgwidth\undefined%
  \global\let\svgscale\undefined%
  \makeatother%
  \begin{picture}(1,0.72522952)%
    \lineheight{1}%
    \setlength\tabcolsep{0pt}%
    \put(0,0){\includegraphics[width=\unitlength,page=1]{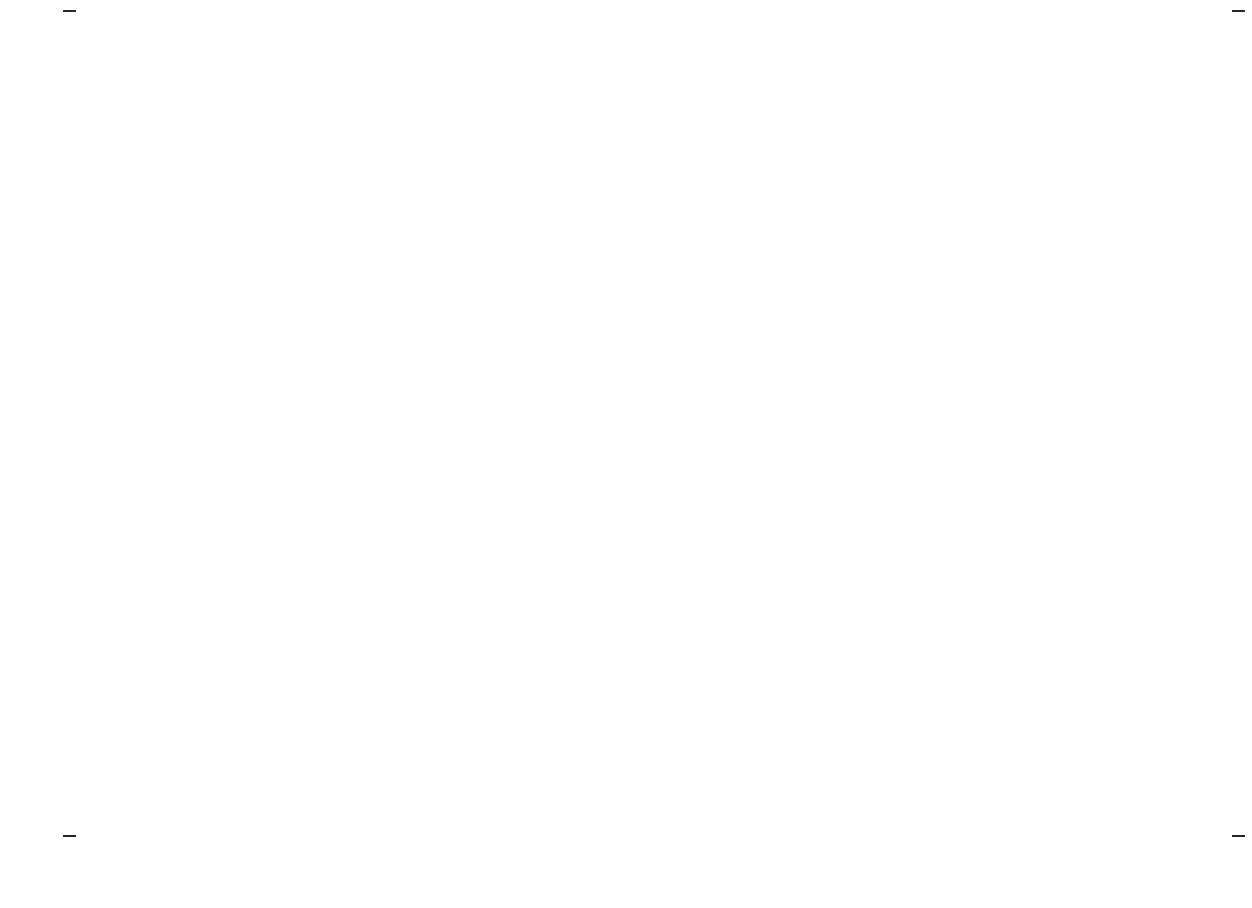}}%
    \put(0.02067476,0.31024707){\color[rgb]{0,0,0}\rotatebox{90}{\makebox(0,0)[lt]{\lineheight{0}\smash{\begin{tabular}[t]{l}$q_2$ in m\end{tabular}}}}}%
    \put(0.44719951,0.00641768){\color[rgb]{0,0,0}\makebox(0,0)[lt]{\lineheight{0}\smash{\begin{tabular}[t]{l}$q_1$ in m\end{tabular}}}}%
    \put(0,0){\includegraphics[width=\unitlength,page=2]{2M_Sym_2_svg-tex.pdf}}%
  \end{picture}%
\endgroup%

      \subcaption{$(M_{C}, V_{s1}), R_2$}
      \label{s4f4_2}
    \end{subfigure}&
    \begin{subfigure}{.23\textwidth}
      \def\svgwidth{10em}\scriptsize
\begingroup%
  \makeatletter%
  \providecommand\color[2][]{%
    \errmessage{(Inkscape) Color is used for the text in Inkscape, but the package 'color.sty' is not loaded}%
    \renewcommand\color[2][]{}%
  }%
  \providecommand\transparent[1]{%
    \errmessage{(Inkscape) Transparency is used (non-zero) for the text in Inkscape, but the package 'transparent.sty' is not loaded}%
    \renewcommand\transparent[1]{}%
  }%
  \providecommand\rotatebox[2]{#2}%
  \newcommand*\fsize{\dimexpr\f@size pt\relax}%
  \newcommand*\lineheight[1]{\fontsize{\fsize}{#1\fsize}\selectfont}%
  \ifx\svgwidth\undefined%
    \setlength{\unitlength}{599.88375302bp}%
    \ifx\svgscale\undefined%
      \relax%
    \else%
      \setlength{\unitlength}{\unitlength * \real{\svgscale}}%
    \fi%
  \else%
    \setlength{\unitlength}{\svgwidth}%
  \fi%
  \global\let\svgwidth\undefined%
  \global\let\svgscale\undefined%
  \makeatother%
  \begin{picture}(1,0.72618918)%
    \lineheight{1}%
    \setlength\tabcolsep{0pt}%
    \put(0,0){\includegraphics[width=\unitlength,page=1]{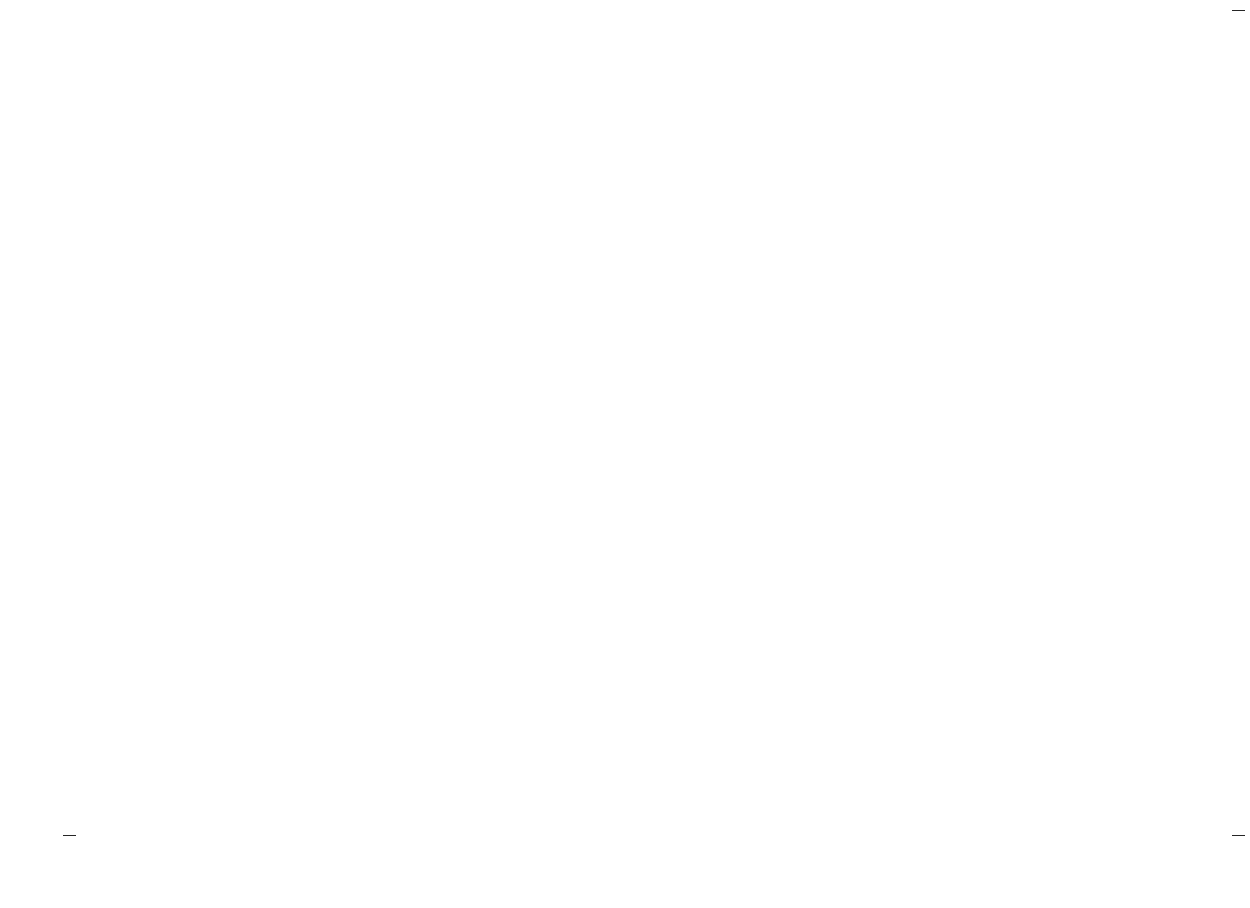}}%
    \put(0.02071105,0.31079154){\color[rgb]{0,0,0}\rotatebox{90}{\makebox(0,0)[lt]{\lineheight{0}\smash{\begin{tabular}[t]{l}$q_2$ in m\end{tabular}}}}}%
    \put(0.44798432,0.00642894){\color[rgb]{0,0,0}\makebox(0,0)[lt]{\lineheight{0}\smash{\begin{tabular}[t]{l}$q_1$ in m\end{tabular}}}}%
    \put(0,0){\includegraphics[width=\unitlength,page=2]{2M_ASym_1_svg-tex.pdf}}%
  \end{picture}%
\endgroup%

      \subcaption{$(M_{C}, V_{a}), R_1$}
      \label{s4f5}
    \end{subfigure}&
    \begin{subfigure}{.23\textwidth}
      \def\svgwidth{10em}\scriptsize
\begingroup%
  \makeatletter%
  \providecommand\color[2][]{%
    \errmessage{(Inkscape) Color is used for the text in Inkscape, but the package 'color.sty' is not loaded}%
    \renewcommand\color[2][]{}%
  }%
  \providecommand\transparent[1]{%
    \errmessage{(Inkscape) Transparency is used (non-zero) for the text in Inkscape, but the package 'transparent.sty' is not loaded}%
    \renewcommand\transparent[1]{}%
  }%
  \providecommand\rotatebox[2]{#2}%
  \newcommand*\fsize{\dimexpr\f@size pt\relax}%
  \newcommand*\lineheight[1]{\fontsize{\fsize}{#1\fsize}\selectfont}%
  \ifx\svgwidth\undefined%
    \setlength{\unitlength}{596.01407241bp}%
    \ifx\svgscale\undefined%
      \relax%
    \else%
      \setlength{\unitlength}{\unitlength * \real{\svgscale}}%
    \fi%
  \else%
    \setlength{\unitlength}{\svgwidth}%
  \fi%
  \global\let\svgwidth\undefined%
  \global\let\svgscale\undefined%
  \makeatother%
  \begin{picture}(1,0.71722663)%
    \lineheight{1}%
    \setlength\tabcolsep{0pt}%
    \put(0,0){\includegraphics[width=\unitlength,page=1]{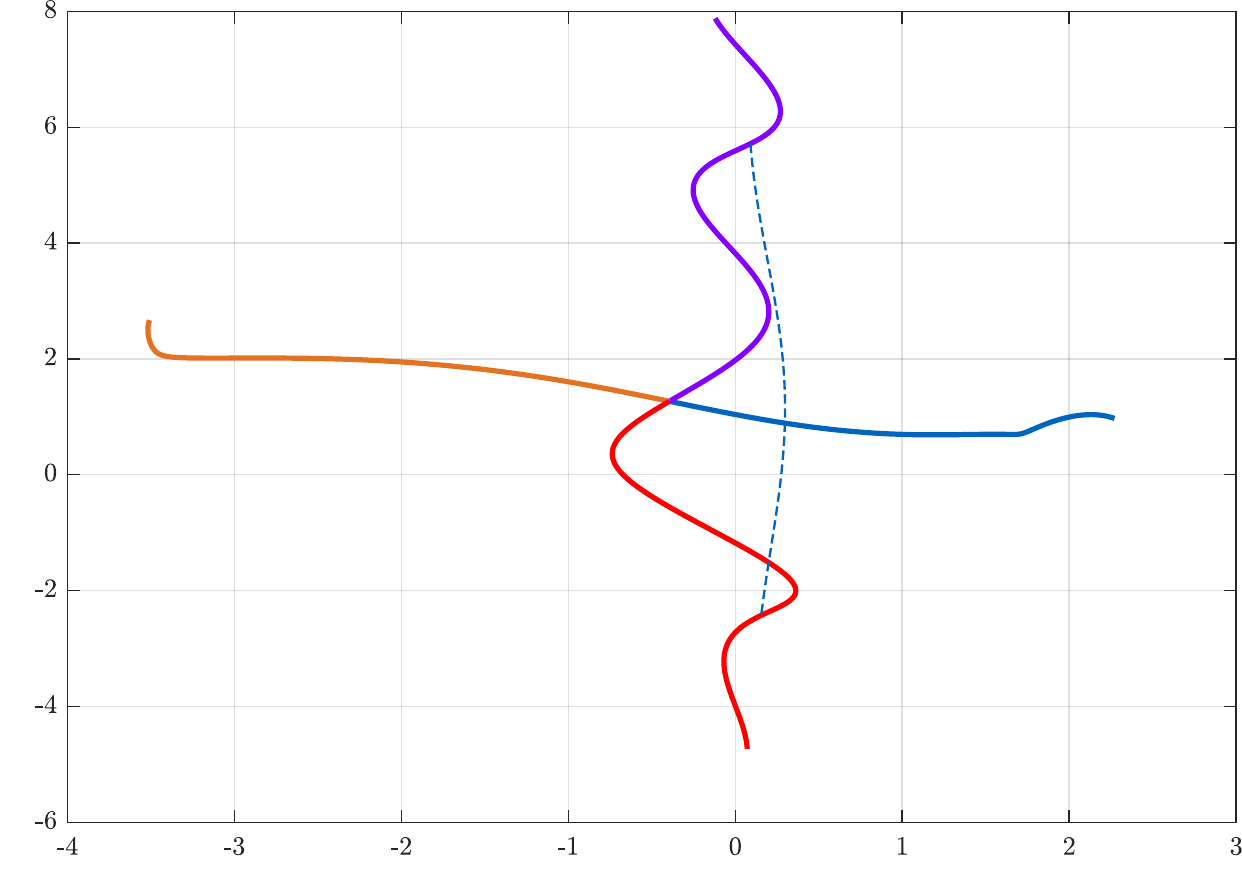}}%
    \put(0.02084551,0.30443212){\color[rgb]{0,0,0}\rotatebox{90}{\makebox(0,0)[lt]{\lineheight{1.25}\smash{\begin{tabular}[t]{l}$q_2$ in m\end{tabular}}}}}%
    \put(0.46538209,0.00647068){\color[rgb]{0,0,0}\makebox(0,0)[lt]{\lineheight{1.25}\smash{\begin{tabular}[t]{l}$q_1$ in m\end{tabular}}}}%
  \end{picture}%
\endgroup%

      \subcaption{$(M_{C}, V_{a}), R_2$}
      \label{s4f5_2}
    \end{subfigure}
  \end{tabular}

  \caption{\small Results of numerical continuation for different example systems generated by combining different potential functions and inertia tensors. The solid lines show computed generators of the system. The pair of blue and orange solid lines show the two generators associated to the first Eigenmanifold and the pair of purple and red the two associated to the second Eigenmanifold. The dashed blue line shows an example modal oscillation of the system for one energy. The conditions stated by Theorem~\ref{Theorem:1} are satisfied for the cases (a,b) and (h,i).
  }\label{fig:allplots}
\end{figure*}
%
%
%
%
%
%
%
%
\subsection{5-DoF System}
\begin{figure}
\begin{tabular}{cc}
\begin{subfigure}[b]{0.1\textwidth}
  \centering
  \def\svgwidth{\columnwidth}\tiny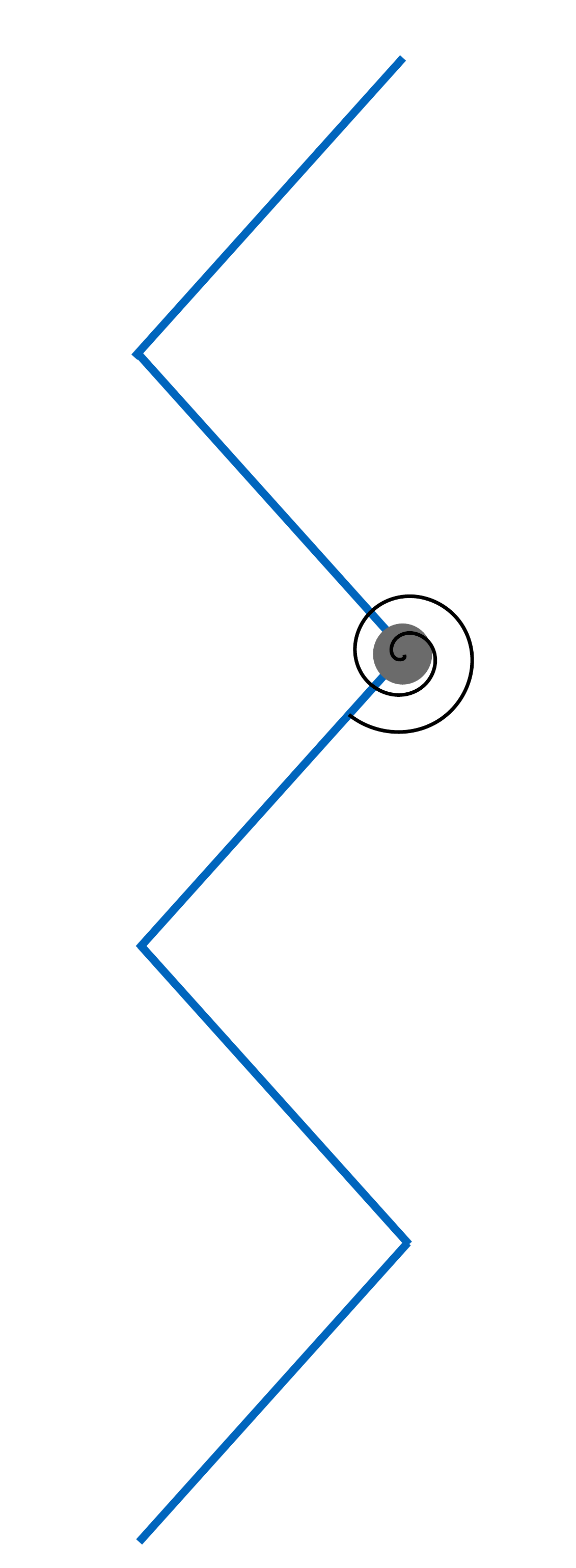
  \vspace{3cm}
  \subcaption{Quintuple elastic pendulum with gravity}
  \label{fig:qp}
\end{subfigure}
&
\begin{subfigure}[b]{0.28\textwidth}
  \def\svgwidth{\textwidth}\tiny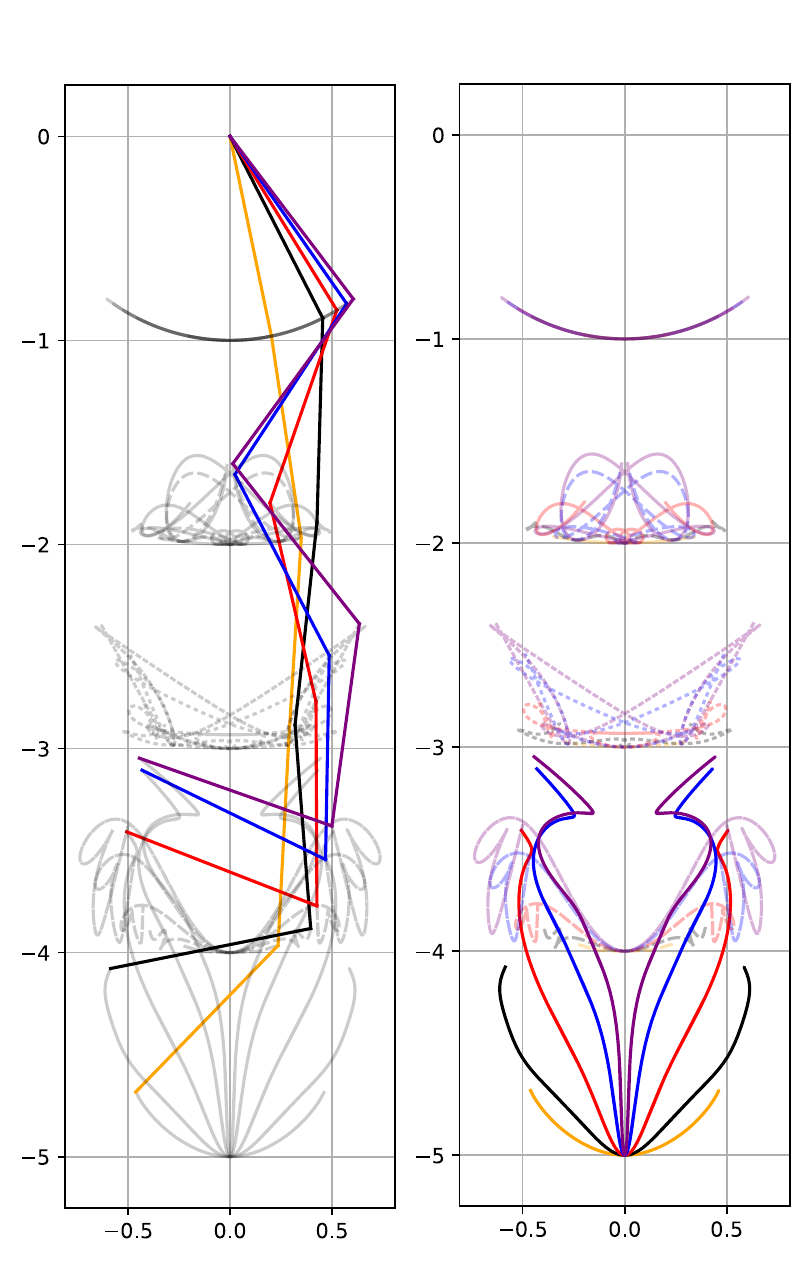
  \subcaption{Modes of different energy levels on the fourth Rosenberg manifold, in Cartesian space. The left shows a number of initial configurations, and the right highlights the modal oscillations of the tip.}
  \label{fig:qp_gen_4}
\end{subfigure}
\end{tabular}
\caption{Quintuple pendulum and various Rosenberg modes.}
\end{figure}

A quintuple pendulum is studied, see Figure~\ref{fig:qp}. The system allows application of the non-resonance condition when there are more than $2$ families of oscillations, and highlights that the symmetry conditions of Theorem~\ref{Theorem:1} allows finding families of geometric Rosenberg modes in a more complex system.

The system is described as follows: for the parameters, all masses are set to $m = 0.4\mathrm{kg}$ and link lengths to $l = 1.0 \mathrm{m}$. The coordinates $q$ are given by joint-angles in radians and the equilibrium configuration $q=0$ corresponds to the straight-down reference configuration. For the potential, a diagonal stiffness matrix $K = 20 I_5 \frac{\mathrm{Nm}}{\mathrm{rad}}$ determines a spring term $V_K(q) = \frac{1}{2} q^\top K q$ and $g = -9.81 \frac{\mathrm{m}}{\mathrm{s}^2}$ a gravitational term $V_g(q)$, which are summed
\begin{equation}
  V(q) = V_K(q) + V_g(q)\,.
\end{equation}
The full inertia tensor $M(q)$ is relegated to the supplementary code~\cite{data-Link}, but we remark that the conditions of Theorem~\ref{Theorem:1} are fulfilled: $M(q) = M(-q)$ and $V(q) = V(-q)$, and the minimum $\qeq = (0,0,0,0,0)$ of $V(q)$ is unique.\\
The conditions of Theorem~\ref{thm:Eigenmanifold_Condition} are fulfilled, as is shown in the supplementary code~\cite{data-Link}. Focussing on Corollary~\ref{cr:n-Eigenmanifolds}, the eigenvalues of $M(\qeq)^{-1}\frac{\partial^2 V}{\partial q \partial q}(\qeq)$ are
\begin{align}
   \lambda_1^2 = -2.5957 \,,\; \lambda_2^2 = &-14.83\,, \lambda_3^2 = -37.05\,,\\ 
  \;\lambda_4^2 = -299.2911&\,,\;  \lambda_5^2 = -441.48 \,. \nonumber
\end{align}
These are mutually non-resonant: the closest-to-integer ratios are $\frac{\lambda_5}{\lambda_1} \approx 13.04$ and $\frac{\lambda_5}{\lambda_4}\approx 1.21$.  
Thus, Corollary~\ref{cr:n-Eigenmanifolds} predicts the existence of 5 weak Eigenmanifolds. By Theorem~\ref{Theorem:1} these are also Rosenberg manifolds.

\begin{rem}
  The ratio $\frac{\lambda_1+\lambda_3}{2\lambda_2} = 0.9996$ is near a combination resonance. This interferes with existence and uniqueness of spectral submanifolds~\cite{Haller2016}, which exist in the presence of dissipation. For beams, it can also cause external forcing to excite multiple modes~\cite{Avramov2002} and energy transfer between modes~\cite{Ruzziconi2026}.
  However, the combination resonance is not relevant for the existence and uniqueness of Eigenmanifolds, where there is no dissipation. In practical application, i.e., with small dissipative terms, non-existence of SSMs may reduce energy-efficiency when target trajectories no longer exactly correspond to solutions. For low-dimensional systems, the effects noted by~\cite{Haller2016} should be limited if an external controller compensates for dissipative terms, and those of~\cite{Avramov2002,Ruzziconi2026} if the external controller stabilizes specific oscillations. 
\end{rem}

\begin{figure*}
    \begin{tabular}{c}
    \begin{subfigure}[b]{\textwidth}
      \def\svgwidth{1\textwidth}\tiny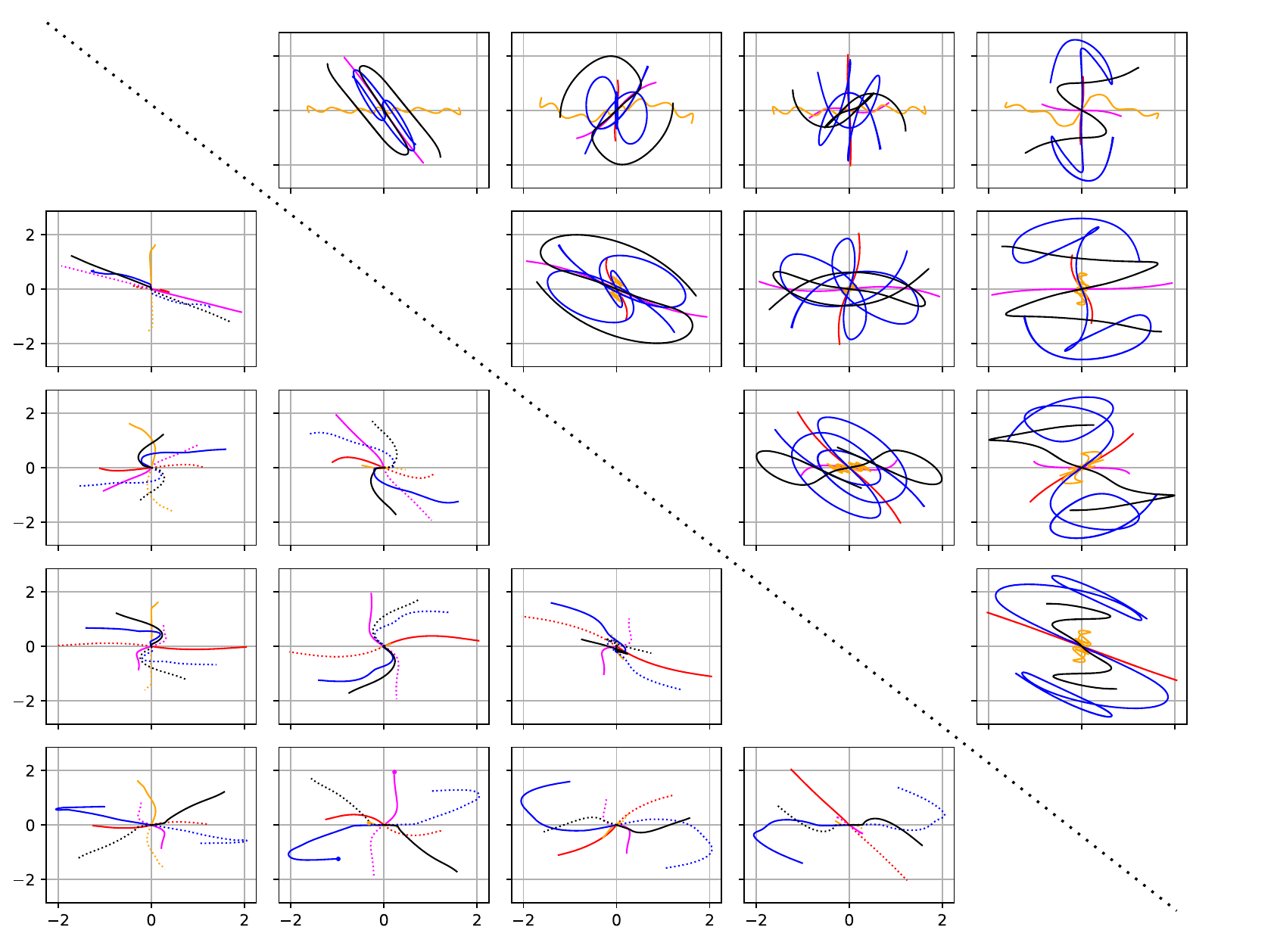
      \subcaption{Generators and high energy modal oscillations of the quintuple pendulum shown in configuration space. The lower triangular matrix of plots shows the generators projected onto different $q_i q_j$-planes. The upper triangular matrix of plots shows the corresponding modal oscillation for the highest energy on the color-matching generator. The system satisfies Theorem~\ref{Theorem:1}, so all modal oscillations pass through the equilibrium.}
      \label{fig:qp_all}
    \end{subfigure}
    \\
    \begin{subfigure}[b]{\textwidth}
      \hspace{1cm}
      \def\svgwidth{\textwidth}\scriptsize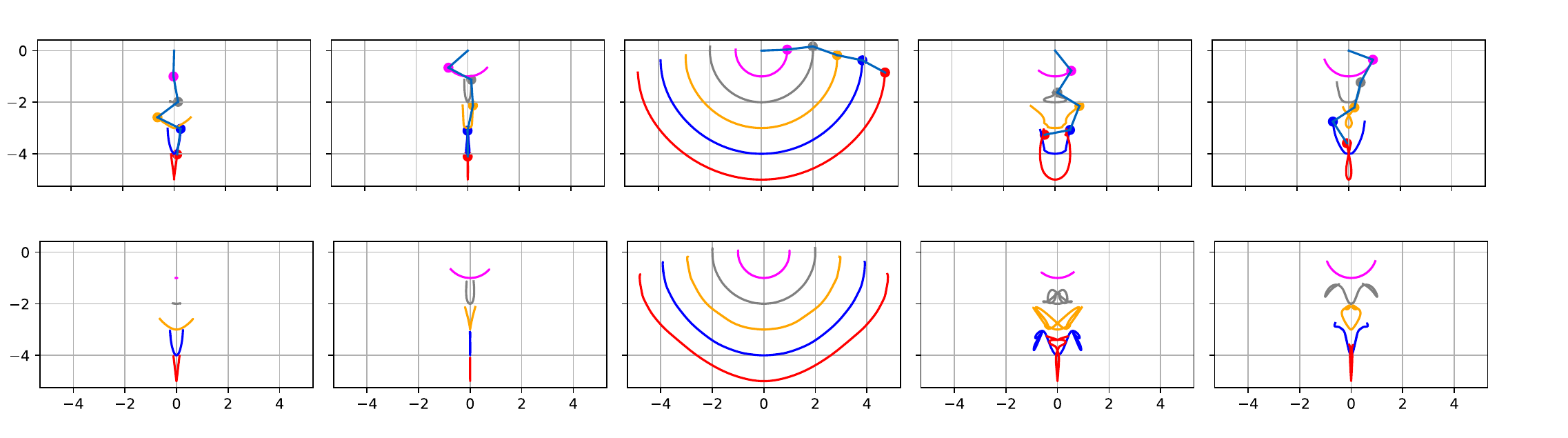
      \subcaption{Cartesian version of Fig.~\ref{fig:qp_all}. The top row shows the generators as Cartesian paths of the joints and the bottom row shows one modal oscillation. All modal oscillation pass through the equilibrium, in which all links of the pendulum point straight down.}
      \label{fig:qp_cart}
    \end{subfigure}
    \end{tabular}
    \caption{Generators of a quintuple pendulum.}
\end{figure*}

%
%
%
To numerically show this, the Eigenmanifolds of the quintuple pendulum are computed up to an energy level of $E_{max} = 100\mathrm{J}$.
The lower triangular matrix in Fig.~\ref{fig:qp_all} shows the projections of the corresponding generators onto $q_i q_j$-planes. 
For each generator, the configuration of maximal potential energy is indicated by dots in the figure. When these configurations are taken as initial configuration for simulating the pendulum. The resulting oscillations pass through the equilibrium configuration and classify as geometric Rosenberg modes. Projections of these trajectories are shown in the upper triangular matrix in Fig.~\ref{fig:qp_all}.
%
Only the generators and highest energy modal oscillation are shown in Fig.~\ref{fig:qp_cart}, but we report that the modal oscillation for all energy levels pass through the equilibrium configuration.
Various modal oscillations on the fourth generator are shown in Fig.~\ref{fig:qp_gen_4}.
Again, all modal oscillation pass through a configuration where all links point straight downwards. This is seen more clearly in the supplementary videos. We note that although $\frac{\lambda_5}{\lambda_1} \approx 13.04$ is close to resonance, this had no effect on the numerical results. 
%
%
\section{Conclusion}\label{ch-EM1:sec:conclusions} Robust, geometric properties of solutions to conservative mechanical systems were investigated, starting from Lyapunov subcenter manifold (LSM) theory.
For conservative mechanical systems, results showed that a time-symmetry guarantees that unique LSMs have the properties of weak Eigenmanifolds: these were defined to be families of oscillations where all orbits are periodic brake trajectories. When the configuration space is two~-~dimensional they are proven to be Eigenmanifolds: as opposed to weak Eigenmanifolds, the configuration trajectory of collected oscillations does not self-intersect. The generator, defined as the collection of brake points on a given Eigenmanifold, was proven to be a uniquely defined, connected 1D submanifold. Finally, conditions on the Riemannian inertia tensor and potential energy were presented that resulted in Eigenmanifolds to fulfill the stronger properties of Rosenberg manifolds: families in which all oscillations pass through the equilibrium. 
Numerical results confirm the validity of the presented Theorems, and present a first step towards stronger properties of nonlinear normal modes as practically feasible target trajectories. 









\bibliographystyle{plainnat}        

\appendix

\section{Proofs}\label{ap:additional_proofs}


\subsection{Time-symmetric Lyapunov subcenter manifolds}\label{ap-ssec:time-symmetry}

We prove Theorem~\ref{thm:Eigenmanifold_Condition} in steps, treating the properties 1 through 4 in order.

\subsubsection{Preliminary: brake points in time-symmetric Hamiltonian systems}
We begin by repeating a well-known theorem on orbits in time-symmetric Hamiltonian systems.

\begin{thm}[Theorem 4.1,~\cite{Lamb1998}]\label{thm:time-symmetry-trajectory}
    {Given Hamiltonian dynamics with $H(\state,P) = H(\state,-P)$.
    \begin{enumerate}[label= (\alph*)]
        \item An evolution $\stateH(t) = (\state(t),P(t))$ 
        of~\eqref{eq:abstract_hamiltonian_vector_field} satisfies 
        \begin{equation}\label{eq:time-symmetry-trajectory}
            \stateH(\mathbb{R}) = {\sigma}_1(\stateH(\mathbb{R}))\,.
        \end{equation}
        if and only if there is $t'$ such that $P(t') = 0$. 
        \item Given $T > 0$, a  $T$-periodic orbit $(\state(t),P(t)) =(\state(t+T),P(t+T))$
        encounters either 0 or 2 distinct points with $P(t)=0$. 
    \end{enumerate}
    }
    \end{thm}
    
    \begin{pf}
        {This directly follows from Theorem 4.1a in~\cite{Lamb1998} with $\text{Fix}(\sigma_1) = \Q \times 0 \subset T^*\Q$, \textit{i.e.}, the fixed set of the symmetry $\sigma_1$ is the zero section of $T^*\Q$.}
    \end{pf}
    By Theorem~\ref{thm:time-symmetry-trajectory}a, any trajectory $(\state(t),P(t))$ with $P(0) = 0$ satisfies
    \begin{equation}\label{eq:backward_ev}
        \state(t) = \state(-t)\,.
    \end{equation}
    \textit{i.e.}, the backwards and forwards evolution of the configuration from a starting point $(\state(0),0)$ are identical. Braking points appear to reflect $\state(t)$. If the orbit is also periodic, $\state(t)$ that satisfies~\eqref{eq:backward_ev} must reflects at a second braking point. Theorem~\ref{thm:time-symmetry-trajectory}b shows that there can be no further reflections in between. 

    \subsubsection{Period between brake points}\label{ap-ssec:period}
    We show that the period and the time between brake points are related: given a period $T_\stateH$, the time to go from one brake point to another is $\Delta t = T_\stateH /2$.
    \begin{lem}\label{ch-EM1:Lemma:2}
    Given an evolution $\stateH(t) = (\state(t),P(t))$ of~\eqref{eq:abstract_hamiltonian_vector_field}. If there are
    two distinct points $(\state_1, 0)$, $(\state_2, 0) \in \stateH(\mathbb{R})$ and denoting $\Delta t = \min \{ t > 0 \, | \, (\state_2, 0) = \Psi^t_{X_H} (\state_1, 0) \}$, then the evolution is periodic with period $T_\stateH = 2\Delta t$. 
    \end{lem}
    
    \begin{pf}
    Let $\stateH(t)$ be such that $\stateH(0) = (\state_1, 0)$ and $\stateH(\Delta t) = (\state_2,0)$. We show that 
    $\stateH(-\Delta t) = \stateH(\Delta t)$, such that the evolution must be periodic with $T_\stateH = 2\Delta t$.  
    Let 
    \begin{equation}\label{ch-EM1:Lemma:2:eq-time-symmetry-period}
        \bar{\stateH}(t) = {\sigma}_1 \stateH(\tau_1t) = \big(\state(-t),-P(-t)\big)
    \end{equation}
    Substituting directly into~\eqref{ch-EM1:Lemma:2:eq-time-symmetry-period}:
    \begin{align}
        \bar{\stateH}(-\Delta t) = (\state(\Delta_t),-P(\Delta_t)) = (\state_2, 0) \,.
    \end{align}
    However, $\bar{\stateH}(0) = (\state_1,0)$ and it follows from uniqueness of solutions that $\stateH(t) = \bar{\stateH}(t)$, so $\stateH(-\Delta t) = \bar{\stateH}(-\Delta t) = \stateH(\Delta t)$, as required. 
    \end{pf}

    \subsubsection{Immersion and brake point properties}\label{ap-ssec:immersion-property}

    \begin{lem}\label{lemma:immersion-time-symmetric-modes}
        Given $T > 0$, a $T$-periodic orbit $\stateH(t)$ there is an immersion $\phi_\stateH:[0,1]\rightarrow\Q$ satisfying~\eqref{eq:immersion-eigenmode} if and only if $\stateH(t)$ satisfies the symmetry
            \begin{equation}
                \stateH(\mathbb{R}) = {\sigma}_1(\stateH(\mathbb{R}))\,.
            \end{equation}
    \end{lem}
    
    \begin{pf}
        If: given the symmetry, Theorem~\ref{thm:time-symmetry-trajectory}a guarantees a point with $P(t')=0$. Then periodicity and Theorem~\ref{thm:time-symmetry-trajectory}b, Lemma~\ref{ch-EM1:Lemma:2} guarantee that $P(t'+T/2) =0$. Then the segment $\state([t',t'+T/2]) \subset \Q$ is such that
        \begin{equation}
            \state([t',t'+T/2]) =  \state([t'+T/2,t'+T])\,.
        \end{equation}
        The map
        \begin{equation}
            \phi_\stateH(s) = \pi_{\Q}(\stateH(t'+sT/2))
        \end{equation}
        is an immersion mapping the interval $[0,1]$ to $\state([t',t'+T/2])$, where $\frac{\ext}{\ext s}\phi_\stateH(s) \neq 0$ for $s\in (0,1)$ since $P \neq 0$, and at $s=0,s=1$, we have that $\phi_\stateH(0) = \state(t')$ and $\phi_\stateH(1) = \state(t'+T/2)$ as required.
        Since, $\state([t',t'+T/2]) =  \state([t'+T/2,t'+T])$, this extends to 
        \begin{equation}
            \phi_\stateH(s) = \pi(\stateH(\mathbb{R}))\,,
        \end{equation}
        and $\stateH(t)$ is a weak geometric eigenmode.
        
        Only if: given a periodic orbit that for which $\phi_\stateH(s)$ can be found, there are necessarily two points with $P(t)=0$, and the symmetry is satisfied as a consequence of Theorem~\ref{thm:time-symmetry-trajectory}. 
    \end{pf}

    \begin{thm}[Time-symmetric modes]\label{thm:time-symmetric-modes-in-time-symmetric-LSMs}
        Given the conditions of Theorem~\ref{thm:LSM_mechanics_geometric}, all $\stateH(t)$ in the unique LSM $\EigN \subset T^*\Q$ are such that
        \begin{equation}
            \stateH(\mathbb{R}) = {\sigma}_1(\stateH(\mathbb{R}))\,.
        \end{equation}
    \end{thm}

    \begin{pf}
        The symmetry $(S_1,\tau_1)$ satisfies the conditions of Theorem~\ref{thm:symmetry_propagation_system_to_LSM}:
        \begin{equation}
            S_1((0,0)) = (0,0)\,,
        \end{equation}
        and the eigenspace $E_0 := D\oplus M(\eqstate)D$ is such that 
        \begin{equation}
            E_0 = \frac{\partial S_1}{\partial \coordC} E_0 = \begin{bmatrix}
                I & 0 \\ 0 & -I
            \end{bmatrix} E_0\,.
        \end{equation}
        Thus, the LSM $\EigN$ is such that 
        \begin{align}
            \EigN &= \sigma_1(\EigN)\,, \\
            \stateH(\R) &= \sigma_1\big(\stateH(\R)\big)\,.
        \end{align}
    \end{pf}
    Then the immersion property of Theorem~\ref{thm:Eigenmanifold_Condition} is a direct consequence of Theorem~\ref{thm:time-symmetric-modes-in-time-symmetric-LSMs} and Lemma~\ref{lemma:immersion-time-symmetric-modes}. The brake point property of Theorem~\ref{thm:Eigenmanifold_Condition} is a direct consequence of Theorem~\ref{thm:time-symmetric-modes-in-time-symmetric-LSMs} and Lemma~\ref{ch-EM1:Lemma:2}.

\subsubsection{Embedding Property}\label{ap-ssec:embedding-property}
    \begin{thm}[Embedding Property]
        Given the conditions of Theorem~\ref{thm:LSM_mechanics_geometric} and the unique LSM $\EigN\subset T^*\Q$, there is a neighborhood $\EigN'\subset \EigN$ containing the equilibrium where there is an embedding $\phi_\stateH:[0,1]\rightarrow\Q$ satisfying~\eqref{eq:immersion-eigenmode}. When $\dim\Q = 2$, then $\EigN' = \EigN$.
    \end{thm}

    \begin{pf}
        By the immersion property~\ref{thm:Eigenmanifold_Condition}, an immersion $\phi_\stateH$ is already guaranteed. To be an embedding, $\phi_\stateH$ must further be injective, \textit{i.e.}, $\state(t)$ should not self-intersect. 
        The limiting behavior of $\stateH(t) = \Psi^t_{X_H}(\stateH_0) \in \mathcal{M}_0$ as $\stateH_0 = (\state,0)$ approaches $(\eqstate,0)$ 
        is that of a linear oscillation, which necessarily does not self-intersect. Hence, there must also be a neighborhood $\EigN'_0 \subset \EigN$ where $\phi_\stateH$ is an embedding.

        We now cover the case $\dim \Q = 2$. One of two scenarios must occur in transitioning from $\EigN'$ to $\EigN$ (see also Figures~\ref{fig:1st_order_intersection} and~\ref{fig:2nd_order_intersection}). 
            \begin{figure}[h]
            \centering
            \begin{tabular}{c c}
                \begin{subfigure}{.45\columnwidth}
                    \def\svgwidth{1.4\textwidth}\scriptsize    
\begingroup%
  \makeatletter%
  \providecommand\color[2][]{%
    \errmessage{(Inkscape) Color is used for the text in Inkscape, but the package 'color.sty' is not loaded}%
    \renewcommand\color[2][]{}%
  }%
  \providecommand\transparent[1]{%
    \errmessage{(Inkscape) Transparency is used (non-zero) for the text in Inkscape, but the package 'transparent.sty' is not loaded}%
    \renewcommand\transparent[1]{}%
  }%
  \providecommand\rotatebox[2]{#2}%
  \newcommand*\fsize{\dimexpr\f@size pt\relax}%
  \newcommand*\lineheight[1]{\fontsize{\fsize}{#1\fsize}\selectfont}%
  \ifx\svgwidth\undefined%
    \setlength{\unitlength}{969.79191197bp}%
    \ifx\svgscale\undefined%
      \relax%
    \else%
      \setlength{\unitlength}{\unitlength * \real{\svgscale}}%
    \fi%
  \else%
    \setlength{\unitlength}{\svgwidth}%
  \fi%
  \global\let\svgwidth\undefined%
  \global\let\svgscale\undefined%
  \makeatother%
  \begin{picture}(1,0.58144748)%
    \lineheight{1}%
    \setlength\tabcolsep{0pt}%
    \put(0.62501786,0.00454651){\color[rgb]{0,0,0}\makebox(0,0)[lt]{\lineheight{1.25}\smash{\begin{tabular}[t]{l}$\Q$\end{tabular}}}}%
    \put(0,0){\includegraphics[width=\unitlength,page=1]{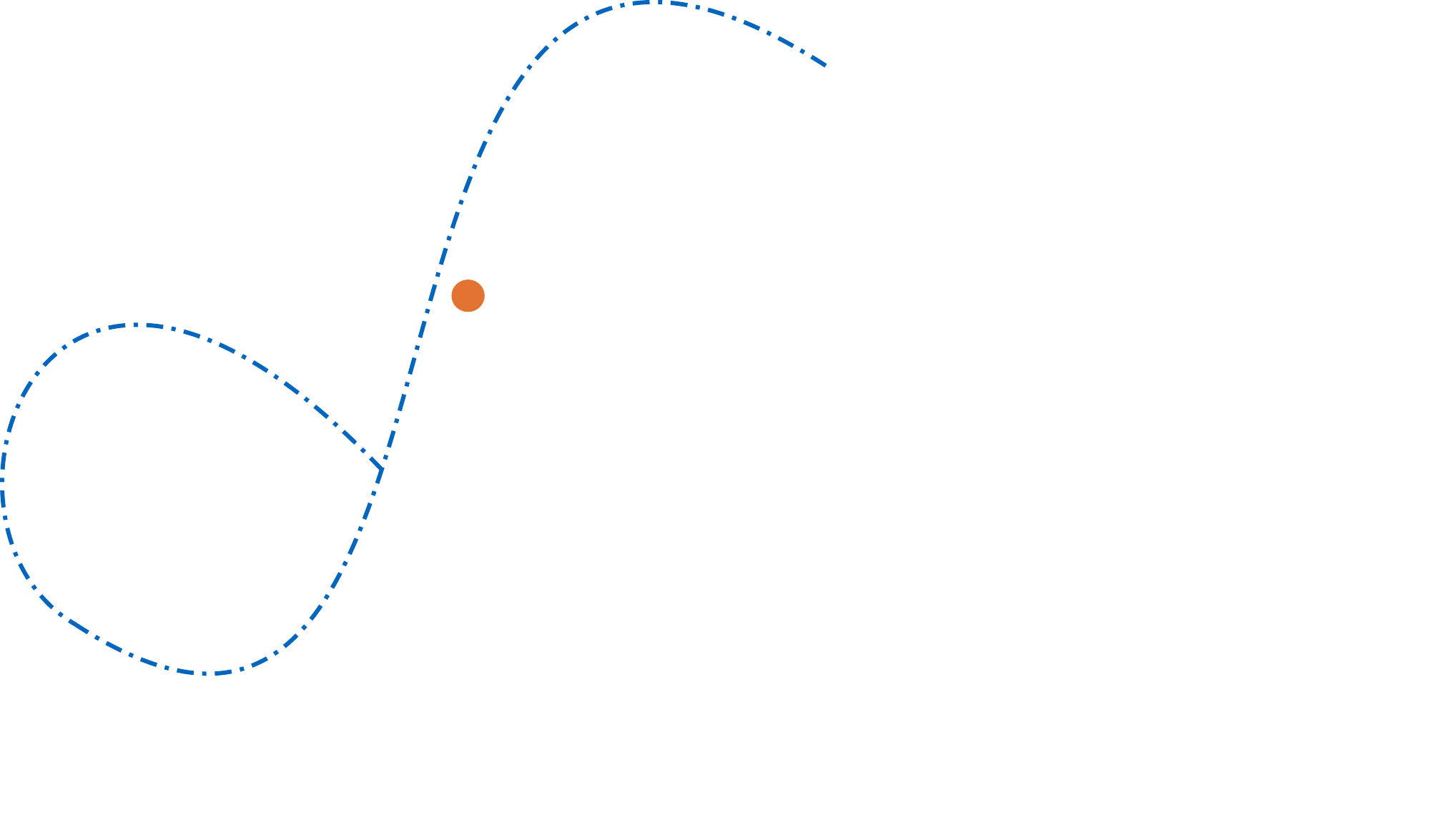}}%
    \put(0.33378264,0.32509504){\color[rgb]{0,0,0}\makebox(0,0)[lt]{\lineheight{1.25}\smash{\begin{tabular}[t]{l}$\eqstate$\end{tabular}}}}%
    \put(0,0){\includegraphics[width=\unitlength,page=2]{intersection_1st_order_svg-tex.pdf}}%
    \put(0.27282944,0.21065077){\color[rgb]{0,0,0}\makebox(0,0)[lt]{\lineheight{1.25}\smash{\begin{tabular}[t]{l}$\left(\state_2, 0\right)$\end{tabular}}}}%
    \put(0.63215587,0.54329727){\color[rgb]{0,0,0}\makebox(0,0)[lt]{\lineheight{1.25}\smash{\begin{tabular}[t]{l}$\left(\state_1, 0\right)$\end{tabular}}}}%
    \put(0.04084009,0.07953985){\color[rgb]{0,0,0}\makebox(0,0)[lt]{\lineheight{1.25}\smash{\begin{tabular}[t]{l}$\left(\state(t), P(t)\right)$\end{tabular}}}}%
  \end{picture}%
\endgroup%

                    \subcaption{First order self-intersection of a geometric eigenmode.}\label{fig:1st_order_intersection}
                \end{subfigure}
                &
                \begin{subfigure}{.45\columnwidth}
                    \def\svgwidth{\textwidth}\scriptsize    
\begingroup%
  \makeatletter%
  \providecommand\color[2][]{%
    \errmessage{(Inkscape) Color is used for the text in Inkscape, but the package 'color.sty' is not loaded}%
    \renewcommand\color[2][]{}%
  }%
  \providecommand\transparent[1]{%
    \errmessage{(Inkscape) Transparency is used (non-zero) for the text in Inkscape, but the package 'transparent.sty' is not loaded}%
    \renewcommand\transparent[1]{}%
  }%
  \providecommand\rotatebox[2]{#2}%
  \newcommand*\fsize{\dimexpr\f@size pt\relax}%
  \newcommand*\lineheight[1]{\fontsize{\fsize}{#1\fsize}\selectfont}%
  \ifx\svgwidth\undefined%
    \setlength{\unitlength}{747.65433123bp}%
    \ifx\svgscale\undefined%
      \relax%
    \else%
      \setlength{\unitlength}{\unitlength * \real{\svgscale}}%
    \fi%
  \else%
    \setlength{\unitlength}{\svgwidth}%
  \fi%
  \global\let\svgwidth\undefined%
  \global\let\svgscale\undefined%
  \makeatother%
  \begin{picture}(1,0.81668627)%
    \lineheight{1}%
    \setlength\tabcolsep{0pt}%
    \put(0.9071709,0.00589733){\color[rgb]{0,0,0}\makebox(0,0)[lt]{\lineheight{1.25}\smash{\begin{tabular}[t]{l}$\Q$\end{tabular}}}}%
    \put(0,0){\includegraphics[width=\unitlength,page=1]{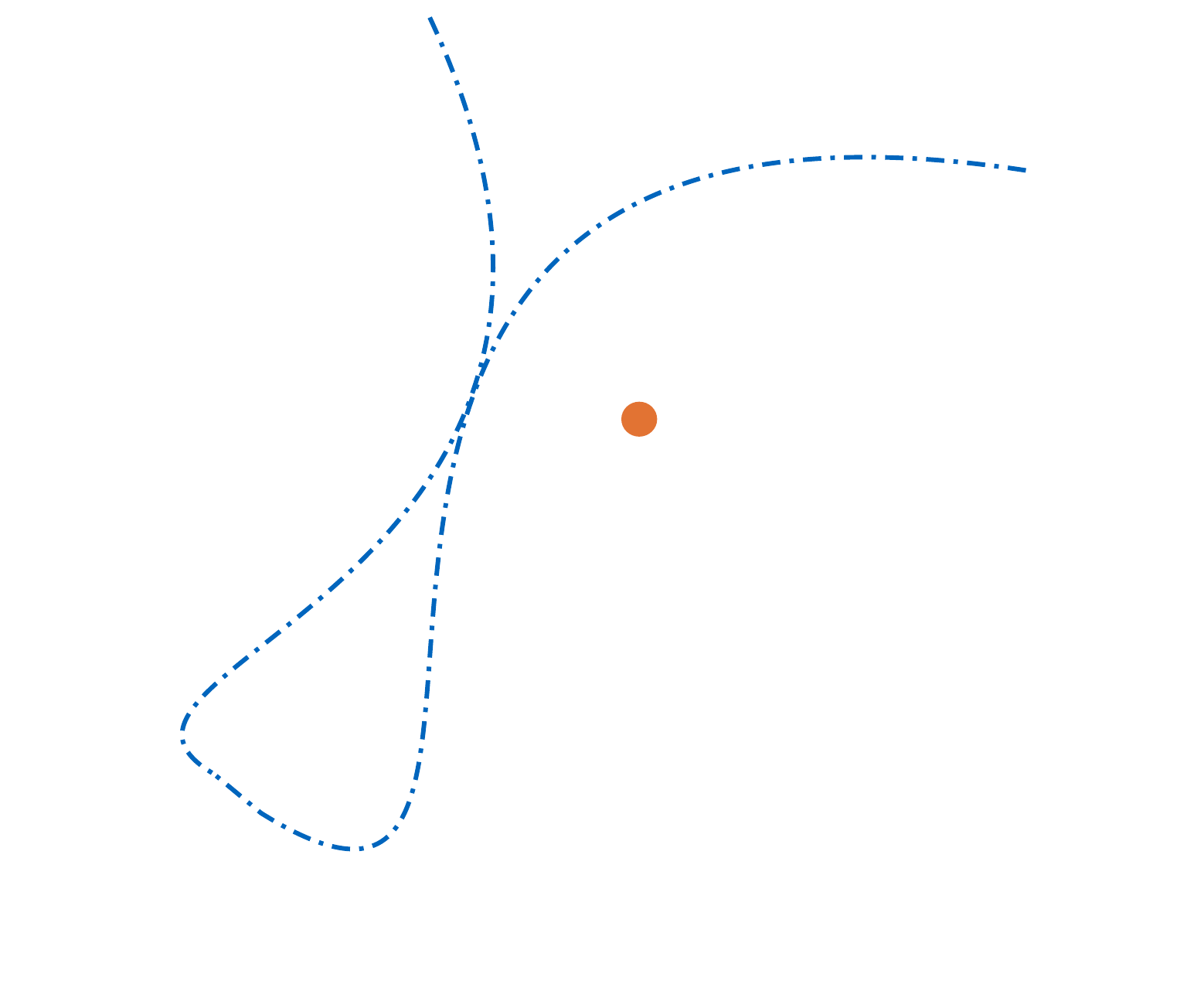}}%
    \put(0.54373651,0.40162219){\color[rgb]{0,0,0}\makebox(0,0)[lt]{\lineheight{1.25}\smash{\begin{tabular}[t]{l}$\eqstate$\end{tabular}}}}%
    \put(0,0){\includegraphics[width=\unitlength,page=2]{intersection_2nd_order_svg-tex.pdf}}%
    \put(-0.00333072,0.39665526){\color[rgb]{0,0,0}\makebox(0,0)[lt]{\lineheight{1.25}\smash{\begin{tabular}[t]{l}$\left(\state, P\right)$\end{tabular}}}}%
    \put(0.16375701,0.08310944){\color[rgb]{0,0,0}\makebox(0,0)[lt]{\lineheight{1.25}\smash{\begin{tabular}[t]{l}$\left(\state(t), P(t)\right)$\end{tabular}}}}%
  \end{picture}%
\endgroup%

                    \subcaption{Higher order self-intersection of a geometric eigenmode. }\label{fig:2nd_order_intersection}
                \end{subfigure}
            \end{tabular}
            \caption{\small Possible self-intersections of Eigenmodes along a generator, assuming a two dimensional configuration space. In transitioning from a geometric eigenmode to a weak geometric eigenmode, either \textbf{(a)} a first order intersection takes places or \textbf{(b)} a higher order intersection takes place. 
            }\label{fig:apf1}
            \end{figure}\\
            Scenario a, a first order intersection where either one of the end-points $\state_1$ or $\state_2$ passes through a $\state_3$ contained within $\pi_{\Q}\big(\stateH(\R)\big)$. A first order intersection would imply e.g., $V(\state_3) = V(\state_2)$ and by conservation of energy the momentum at $\state_3$ must be $P_3 = 0$. By Theorem~\ref{thm:time-symmetry-trajectory}b, such a third point with $P = 0$ on the same mode $\stateH(t)$ is not possible. \\
            Scenario b concerns a higher order intersection where $\stateH(\R)$ is tangent to itself at some point $(\state,P) \in \stateH(\R)$ that is not an end-point. Uniqueness of solutions forbids this transition, since the system must have a unique solution through any $(\state,P)$. \\
            Since the transition from embedded to immersed-but-not-embedded modes is impossible for $\dim \Q = 2$, it must be that all modes are embedded and $\EigN' =\EigN$ for this case.
        \end{pf}
        
        \begin{rem}
            The proof breaks for $\dim \Q = 2$ down in higher dimensions, since other first order intersections of points $(\state,P) \in \stateH(\R)$ that are end-points become possible. However, we expect that self-intersections become less ``likely'' with increasing $n$, being most common for $\text{dim}(\Q) = 3$. 
        \end{rem} 

\subsubsection{Generator Property}\label{proof:generator}
To show that the generator property holds, we prove the following theorem: 
\begin{thm}[Existence of generator]
    Given a unique LSM $\EigN\subset T^*\Q$ that satisfies the immersion property in Theorem~\ref{thm:Eigenmanifold_Condition}. Then the set $\mathcal{R} := \{(\state,0) \,|\, (\state,0) \in \EigN\}$ is a connected, 1D submanifold containing the equilibrium.
\end{thm}

\begin{pf}
We first show that $\mathcal{R} = \{(\state,0) \,|\, (\state,0) \in \EigN\}$ is a regular 1-dimensional submanifold of $\EigN$. 
We use the regular level set theorem \citep[Corr. 5.24]{Lee2012}, which states that $\mathcal{R} \subset \EigN$ is a closed, embedded submanifold of dimension 1 if there is a function $\alpha:\EigN \rightarrow \R$ such that $\mathcal{R} = \alpha^{-1}(0)$ and $\text{Im}(\text{d}\alpha) = T_0\R$. 
To this end, note that 
\begin{equation}
    \alpha(\state,P) = M^{-1}(\text{d}V(\state),P) + M^{-1}(P,P)\,. 
\end{equation}
is as required: points on the generator satisfy that $P = 0$, which corresponds to $\alpha(\state,0) = 0$. For $\text{Im}(\text{d}\alpha) = T_0\R$, note that $\EigN$ is tangent to the flow $X_H$, and apply $\text{d}\alpha(\state,P)$ to the vector $X_H(\state,P)$ at $(\state,0) \in \EigN$:
\begin{align}
    \text{d}\alpha(X_H) &= M^{-1}(\text{d}V(\state),\dot{P}) \\
    &= M^{-1}(\text{d}V(\state),\text{d}V(\state)) \neq 0 \nonumber \,.
\end{align}
Thus, it follows from the regular level set theorem~\cite{Lee2012} that $\mathcal{R}$ is a closed, embedded 1D submanifold. This holds where $\text{d}V(\state) \neq 0$, \textit{i.e.}, at all points other than the equilibrium. The equilibrium itself may be included as a limit point, since geometric eigenmodes arbitrarily close to $(\eqstate,0)$ may be found by definition of $\EigN$.\\
%
%
From the closedness of $\mathcal{R}$ it can be shown that it is connected: Consider coordinates $X:\EigN\rightarrow\R^2$, such that orbits are mapped to lines of constant radius, and angle $\theta$ ranges from $0$ to $2\pi$ over a full oscillation. Then the generator is 
parameterized as 
\begin{equation}
    \mathcal{R} = \big(r(\coordC),\theta(\coordC)\big)\,.
\end{equation}
Let the coordinates be such that the antipodal brake points guaranteed by Lemma~\ref{ch-EM1:Lemma:2} are
\begin{equation}
    \sigma_1(\mathcal{R}) = \big(r(\coordC),-\theta(\coordC)\big)\,.
\end{equation}
Assume, for the sake of contradiction, that $\mathcal{R}$ is not connected but consists of connected, closed submanifolds $\mathcal{R}_1,\mathcal{R}_2\subset \mathcal{R}$. Let $\bar{r}_1$ be the radius on the boundary of $\mathcal{R}_1$. By closure, $\bar{r}_1$ must also lie on the boundary of $\mathcal{R}_2$.   
Yet, if $\mathcal{R}_{0,i}$, $\mathcal{R}_{0,j}$, $\sigma_1(\mathcal{R}_{0,i})$ and $\sigma_1(\mathcal{R}_{0,j})$ are distinct for $\bar{r}_1$, there would have to be geometric eigenmodes with four brake points. This contradicts Theorem~\ref{thm:time-symmetry-trajectory}b, so the generator must be connected. 
\end{pf}

\subsection{Spatial symmetry}\label{ch-EM1:ap:symmetry} 

We show that a set of equivariant coordinates exists in which $\bar{q} = 0$ and $\varphi:\Q\rightarrow\Q$ is the map $q\mapsto -q$ and $S_2(q,p) = (-q,-p)$. Afterwards, we use this set of coordinates to prove the symmetry of the equations.

\subsubsection{Equivariant coordinates}\label{ap:equivariant_coordinates}
Given an diffeomorphism $\varphi:\Q\rightarrow\Q$, that leaves an equilibrium $\eqstate \in \Q$ fixed, satisfies that $\varphi_*(\bar{q}) = -\text{id}_{T_{\bar{q}}\Q}$ and $\varphi\circ\varphi = \text{id}_{\Q}$.\\
Pick any chart $(U,X)$ with chart-region $\eqstate \in U\subset \Q$ and (diffeomorphic) chart-map $X:\Q\rightarrow\R^n$. We want to find a chart $(V',Y)$ such that $S_2(q,p) = (-q,-p)$. For this it is sufficient that $Y\circ\varphi\circ Y^{-1}(q) = - q$, or equivalently
\begin{equation}
      Y\circ\varphi = - Y\,.
\end{equation}
Define a candidate chart $(V,Y)$ by  
\begin{align}
    V &:= \{\state \in U \, | \, \varphi(\state) \in U\}\,, \\
    Y &:= \frac{1}{2} (X - X \circ \varphi)\,.
\end{align}
Then the coordinate description of the symmetry is $Y \circ \varphi(\state) = -Y(\state)$, as desired. Also note that $V$ is well defined since $\varphi\circ\varphi = \text{id}_{\Q}$.\\
It remains to be shown that $Y$ is a valid chart-map, \textit{i.e.}, that $Y$ is a local diffeomorphism. This follows from the inverse function theorem \citep[Theorem 1.3.12]{FoundationsOfMechanics} by showing that $Y_*(\eqstate)$ is full rank:
\begin{equation}
    Y_*(\eqstate) = \frac{1}{2} (X_*(\eqstate) - X_*(\eqstate)\varphi_*(\eqstate)) = X_*(\eqstate)\,. 
\end{equation}
Here, the first step uses that $\varphi(\eqstate) = \eqstate$ and the second step uses that $\varphi_*(\eqstate) = -\text{id}_{T_{\eqstate}\Q}$.
This is indeed full rank, since $X$ is itself a local diffeomorphism. Thus, a neighborhood $V'\subset V$ of $\eqstate$ exists such that $(V',Y)$ is a valid chart with the desired properties.

\subsubsection{Proof of spatial symmetry}

The spatial symmetry is stated in $\varphi$-equivariant coordinates (see also Appendix~\ref{ap:equivariant_coordinates}).

\begin{thm}[Spatial symmetry]\label{thm:spatial-symmetry}
    If $M = \varphi^* M$ and $V = V \circ \varphi$ hold for $\varphi$ satisfying~\eqref{eq:eq-sym-1} to~\eqref{eq:eq-sym-3}, then the symmetry
    \begin{equation}\label{ap-eq:eq_sym}
        \sigma_2((\state,P)) = (\varphi(\state), (\varphi^{-1})^*P)\,, \quad \tau_2 = 1\,. 
    \end{equation}
    satisfies
    \begin{equation}
        X_H\circ\sigma = \tau \sigma_* X_H \,. \label{ap-eq:symmetry_condition}
    \end{equation}
\end{thm}

\begin{pf}
    We use $\varphi$-equivariant coordinates. Then $M = \varphi^*M$ and $V = V\circ\varphi$ read, respectively:
    \begin{align}
        M(q) &= M(-q)\,, \\
        V(q) &= V(-q)\,,
    \end{align} 
    I.e., $V(q)$ and the component-functions $M(q)$ are even functions.
    The dynamics~\eqref{eq:abstract_hamiltonian_vector_field} read
    \begin{equation}
        f\big(\coordC \big) = 
        \begin{pmatrix} 
                M(q)^{-1}p \\
    -\frac{\partial}{\partial q} (p^\top M(q)^{-1} p) - \frac{\partial V}{\partial q}(q)
            \end{pmatrix}\,.
    \end{equation}
    And condition~\ref{ap-eq:symmetry_condition} reads
    \begin{equation}
        \frac{\partial S_2}{\partial \coordC} f\big(\coordC\big) = f\big(S_2(\coordC)\big)\,.
    \end{equation}
    Note that the gradient of an even function $h(q) = h(-q)$ is odd, i.e., $\frac{\partial h}{\partial q}(q) = -\frac{\partial h}{\partial q}(-q)$. Thus 
    the right-hand side is
    \begin{align}
        f\big(S_2(\coordC)\big) &= f(-q,-p) \nonumber \\ &= 
        \begin{pmatrix} 
                - M(-q)^{-1}p \\
                -1/2 \frac{\partial}{\partial s} (p^\top M(s)^{-1} p)_{s=-q} - \frac{\partial V}{\partial s}_{s=-q}\,.
        \end{pmatrix} \nonumber \\ &= 
        \begin{pmatrix} 
                - M(q)^{-1}p \\
                \frac{\partial}{\partial q} (p^\top M(q)^{-1} p) + \frac{\partial V}{\partial q}\,.
        \end{pmatrix}\,.
    \end{align}
    For the left-hand side we get
    \begin{align}
        \tau_1 \frac{\partial S_2}{\partial \coordC} F\big(\coordC\big) &= \begin{bmatrix}
            -I & 0 \\ 0 & -I
        \end{bmatrix} F(\coordC) \nonumber \\  &=         
        \begin{pmatrix} 
                - M(q)^{-1}p \\
    \frac{\partial}{\partial q} (p^\top M(q)^{-1} p) + \frac{\partial V}{\partial q}(q)\,.
        \end{pmatrix} \,.
    \end{align}
    Both sides are the same, so the symmetry condition~\ref{ap-eq:symmetry_condition} holds. This completes the proof.
\end{pf}

In the same $\varphi$-equivariant coordinates, the symmetry means that if a trajectory $\big(q(t),p(t)\big)$ is a solution of the dynamics~\eqref{eq:hamiltonian_ode_q},~\eqref{eq:hamiltonian_ode_p} then also $\big(-q(t),-p(t)\big)$ is a solution.

\subsection{Time- and spatially symmetric Lyapunov subcenter manifolds}\label{proof:rosenberg_manifold}
In this section, the properties guaranteed by Theorem~\ref{Theorem:1} are proven in order.  

\begin{thm}[Spatially symmetric modes]\label{thm:spatially-symmetric-modes-in-spatially-symmetric-LSMs}
    Given the conditions of Theorem~\ref{thm:LSM_mechanics_geometric}, and if $\varphi:\Q\rightarrow\Q$ satisfies properties~\eqref{eq:eq-sym-1} to~\eqref{eq:eq-sym-3}, and $V = V\circ\varphi$, $M = \varphi^*M$, then then all $\stateH(t)$ in the unique LSM $\EigN \subset T^*\Q$ are such that
    \begin{equation}
        \stateH(\mathbb{R}) = {\sigma}_2(\stateH(\mathbb{R}))\,.
    \end{equation}
\end{thm}

\begin{pf}
    Analogous to the proof of~\eqref{thm:time-symmetric-modes-in-time-symmetric-LSMs}, the symmetry $(S_2,\tau_2)$ satisfies the conditions of Theorem~\ref{thm:symmetry_propagation_system_to_LSM}. With 
    \begin{equation}
        S_2((0,0)) = (0,0)\,,
    \end{equation}
    and eigenspace $E_0 := D\oplus M(\eqstate)D$ satisfying 
    \begin{equation}
        E_0 = \frac{\partial S_2}{\partial \coordC} E_0 = 
        \begin{bmatrix}
            -I & 0 \\ 0 & -I
        \end{bmatrix} E_0\,.
    \end{equation}
    Thus, the LSM $\EigN$ is such that 
    \begin{align}
        \EigN &= \sigma_2(\EigN)\,, \\
        \stateH(\R) &= \sigma_2\big(\stateH(\R)\big)\,.
    \end{align}
\end{pf}

\begin{lem}\label{Lemma:4}
Given an evolution $\stateH(t) = (\state(t),P(t))$ of~\eqref{eq:abstract_hamiltonian_vector_field}. If there is a pair of points $(\state,P), (\varphi(\state),{\varphi^{-1}}^*P) \in \stateH(\R)$ 
and denoting $\Delta t = \min \{ t > 0 \, | \, (\varphi(\state),{\varphi^{-1}}^*P) = \Psi^t_{X_H} ((\state,P)) \}$,  
then the evolution is periodic with period $T = 2\Delta t$ and  
\begin{align}\label{s2p2:eq1}
    \state(t + T/2) &= \varphi\big(\state(t)\big)\,,\\ 
    P(t + T/2) &= {\varphi^{-1}}^*P(t)\,. 
\end{align}
\end{lem}

\begin{pf}
Let $\stateH_1(t) = (\state(t),P(t))$ be a solution to~\eqref{eq:abstract_hamiltonian_vector_field},~\eqref{eq:abstract_hamiltonian}. Via~\eqref{eq:eq_sym}
\begin{equation}
    \stateH_2(t) = (\varphi(\state)(t),{\varphi^{-1}}^*P(t))\,,
\end{equation}
must also be a solution. 
By hypothesis, the points $\stateH_1(0) = (\state,P)$ and $\stateH_2(0) = (\varphi(\state),{\varphi^{-1}}^*P)$ lie on the same trajectory, and
\begin{equation}\label{eq:x1_eq_x2}
    \stateH_1(\Delta t) = \stateH_2(0)\,.
\end{equation}
Similarly $\stateH_2(\Delta t) = \stateH_1(0)$, such that 
$\stateH_1(2\Delta t) = \stateH_2(\Delta t) = \stateH_1(0)$, i.e., $\stateH(t)$ is periodic and $T = 2\Delta t$ is the period. 
It also follows from~\eqref{eq:x1_eq_x2} that $\stateH_1(t+T/2) = \stateH_2(t)$, which yields 
\begin{align}
        \state(t + T/2) &= \varphi\big(\state(t)\big)\\
        P(t + T/2) &= {\varphi^{-1}}^*P(t)\,,
\end{align} 
which concludes the proof. 
\end{pf}

\begin{lem}[Time to equilibrium]\label{Lemma:5}
Given an evolution $\stateH(t) = (\state(t),P(t))$ of~\eqref{eq:abstract_hamiltonian_vector_field}. If there is a pair of points $(\state,0), (\varphi(\state),0) \in \stateH(\R)$ and $\Delta t = \min \{ t > 0 \, | \, (\varphi(\state),0) = \Psi^t_{X_H} ((\state,0)) \}$. Then $\stateH(t)$ is periodic with period $T = 2\Delta t$ and there is $\widetilde{\state} \in \Q$ satisfying $\widetilde{\state} = \varphi(\widetilde{\state})$ and $P\in T^*_{\widetilde{\state}}\Q$, such that $(\widetilde{\state},P) = \Psi^{T/4}_{X_H} ((\state,0))$.
\end{lem}

\begin{pf}
Since there are two points with $P = 0$, Lemma~\ref{ch-EM1:Lemma:2} applies and $\stateH(t)$ is periodic with period $T = 2\Delta t$. We set $\stateH(0) = (\state,0)$ and $\stateH(T/2) = (\varphi(\state),0)$.
Equation~\eqref{eq:backward_ev} applies, and a particular case is
\begin{equation}\label{Proof:Lemma:5, eq:1}
    \state(T/4) = \state(-T/4)\,.
\end{equation} 
Points $(\state,0), (\varphi(\state),0) \in \stateH(\R)$ are of the form required for Lemma~\ref{Lemma:4}, so equation~\eqref{s2p2:eq1} holds and a particular case is
\begin{equation}\label{Proof:Lemma:5, eq:3}
    \state(T/4) = \varphi\big(\state(-T/4)\big)\,. 
\end{equation}
Combining~\eqref{Proof:Lemma:5, eq:1} and~\eqref{Proof:Lemma:5, eq:3} yields
\begin{equation}
    \state(-T/4) = \varphi\big(\state(-T/4)\big) \,,
\end{equation} 
i.e., $\state(t)$ is guaranteed to move through $\widetilde{\state} = \state(-T/4)$, which is a fixed point of $\varphi:\Q\rightarrow\Q$\,. The statement $(\widetilde{\state},P) = \Psi^{T/4}_{X_H} ((\state,0))$ follows directly. 
\end{pf}

Property 2.\ of Theorem~\ref{Theorem:1} is a direct consequence of Theorem~\ref{thm:spatially-symmetric-modes-in-spatially-symmetric-LSMs} and Lemma~\ref{Lemma:5}.

\begin{lem}[Equilibrium Property]\label{lemma:equilibrium-property}
    Given $T > 0$. If a $T$-periodic orbit $\stateH(t)$ satisfies 
        \begin{equation}
            \stateH(\mathbb{R}) = {\sigma}_1(\stateH(\mathbb{R})) = {\sigma}_2(\stateH(\mathbb{R})) \,,
        \end{equation}
    and if the fixed point $\eqstate$ of ${\sigma}_1,{\sigma}_2$ is the only equilbrium such that $V(\eqstate) < V(\state)$ for $(\state,0) \in \stateH(\R)$. Then 
    \begin{equation}
        \eqstate \in \pi\big(\stateH(\R)\big)\,.
    \end{equation}
\end{lem}

\begin{pf}
    By Lemma~\ref{Lemma:5}, $\stateH(t)$ that satisfies $\stateH(\R) = \sigma_1\big(\stateH(\R)\big) = \sigma_2\big(\stateH(\R)\big)$ necessarily contains a fixed point $\widetilde{\state} = \varphi(\widetilde{\state})$, and by~\cite{Moehlis2007} such $\widetilde{\state}$ is necessarily an Equilibrium. If $\eqstate$ is the only equilibrium such that $V(\eqstate) < V(\state)$ for $(\state,0) \in \stateH(\R)$, then it necessarily holds that $\eqstate = \widetilde{\state}$ and $\stateH(t)$ passes through the equilibrium $\eqstate$ for some $P$. 
\end{pf}
Property 1.\ is equivalent to Lemma~\ref{lemma:equilibrium-property}.

\section{Additional Material}\label{ap:additional_material}

\subsection{Canonical charts on $T^*\Q$}\label{ap:canonical_chart}
See this section in the supplementary material~\cite{arxiV-Link}.
Define coordinates on $\Q$ by a chart $(U,X)$ with $U\subset \Q$ and $X:U\rightarrow \mathbb{R}^n$, assigning coordinates to $\state \in \Q$ by
\begin{equation}
    (q_1,\hdots,q_n) := X(\state)\,.
\end{equation}
At the point $\state$, this induces a fibre-wise chart $(T_Q^*\Q,{X^{-1}}^*)$ via the pullback ${X^{-1}}^*:T_Q^*\Q \rightarrow \mathbb{R}^n$ that assigns to $P \in T_Q^*\Q$ the coordinates 
\begin{equation}
      (p_1,\hdots,p_n) = {X^{-1}}^*(P) \,.
\end{equation}
A canonical chart on $T^*U \subset T^*\Q$ then assigns to $(\state,P) \in T^*U$ the coordinates $(q_1,\hdots,q_n, p_1, \hdots, p_n)$.

\subsection{The geometric Hessian of a function}

Denote by $\text{d}:C^\infty(\M) \rightarrow \Gamma(T^*\M)$ the differential of a function and by $\nabla_{\cdot} \cdot: \X(\M)\times \Gamma(T^p_q\M)\rightarrow \Gamma(T^p_q\M)$ the covariant derivative w.r.t. a given connection on $\M$. 
The Hessian of $V \in C^2(\M,\R)$ is defined as
\begin{equation}
    \text{Hess}(V)(X,Y) := \big(\nabla_X\text{d}V\big)(Y)\,.
\end{equation}
In a local chart we denote $\text{d}V = \frac{\partial V}{\partial \coordC^j}d\coordC^j$, write $X = X^i\frac{\partial}{\partial \coordC^j}$, $Y = Y^i\frac{\partial}{\partial \coordC^j}$, and the choice of connection corresponds to a choice of Christoffel Symbols $\Gamma^k_{ij}$. Using Einstein summation notation, the Hessian reads
\begin{equation}
    \text{Hess}(V)(X,Y) = X^j \frac{\partial^2 V}{\partial \coordC^i\partial \coordC^j} Y^i - X^j \frac{\partial V}{\partial \coordC^k}\Gamma^k_{ij} Y^i\,.
\end{equation}
At points $m \in \M$ where $\text{d}V = 0$, we have that 
\begin{equation}
    \text{Hess}(V) = \frac{\partial^2 V}{\partial \coordC^i\partial \coordC^j} \text{d}\coordC^i\text{d}\coordC^j \,,
\end{equation}
i.e., $\text{Hess}(V)(X,Y)$ is a $(0,2)$-tensor at critical points of $V$, and is then also independent of the choice of Christoffel symbols $\Gamma^k_{ij}$.

\end{document}